\documentclass[msom,nonblindrev]{informs3} 
\OneAndAHalfSpacedXI 

\usepackage{endnotes}
\usepackage{algorithm}
\usepackage{algpseudocode}
\usepackage{mathrsfs} 
\let\footnote=\endnote
\let\enotesize=\normalsize
\def\notesname{Endnotes}%
\def\makeenmark{$^{\theenmark}$}
\def\enoteformat{\rightskip0pt\leftskip0pt\parindent=1.75em
  \leavevmode\llap{\theenmark.\enskip}}
\usepackage{amssymb, amsmath, amstext}
\usepackage[dutch,english]{babel}
\usepackage{soul}
\usepackage{mathtools}

\usepackage{tikz}
\usetikzlibrary{decorations.pathreplacing,positioning, arrows.meta}
\usepackage{tabularx,booktabs}%
\usepackage{enumerate}
\usepackage{bm}
\usepackage{nicefrac}
\usepackage{caption}
\usepackage{subcaption}

\usepackage{natbib}
 \bibpunct[, ]{(}{)}{,}{a}{}{,}%
 \def\bibfont{\small}%
 \def\bibsep{\smallskipamount}%
\usepackage{bigstrut}
\usepackage{multirow}
\usepackage{dsfont}
\usepackage{makecell}
\usepackage{enumitem}
\setlist{nolistsep}
\usepackage{authblk}
\usepackage[font=footnotesize,labelfont=bf]{caption}
\usepackage{algorithmicx}
\usepackage{algorithm}
\usepackage{algpseudocode}
\TheoremsNumberedThrough     
\ECRepeatTheorems

\EquationsNumberedThrough    

\def\N{\mathbb{N}}
\def\R{\mathbb{R}}

\def \E {\mathbb{E}}

\def \P {\mathbb{P}}

\def \-> {\rightarrow}

\usepackage{calrsfs}

\def \leqst {\leq_{\text{st}}}
\def \eqst {=_{\text{st}}}

\usepackage{comment}
\usepackage[hidelinks]{hyperref}
\begin{document}


\RUNAUTHOR{Drent, Drent, and Arts}

\RUNTITLE{Condition-Based Production for Stochastically Deteriorating Systems}

\TITLE{Condition-Based Production for Stochastically Deteriorating Systems: Optimal Policies and Learning}

\ARTICLEAUTHORS{%
\AUTHOR{Collin Drent, Melvin Drent}
\AFF{School of Industrial Engineering, Eindhoven University of Technology, Eindhoven, the Netherlands, \\ \EMAIL{
\{c.drent, m.drent\}@tue.nl}}
\AUTHOR{Joachim Arts}
\AFF{Luxembourg Centre for Logistics and Supply Chain Management, University of Luxembourg, Luxembourg-City, Luxembourg, \EMAIL{joachim.arts@uni.lu}}
} 

\ABSTRACT{\textbf{Problem Definition:} Production systems deteriorate stochastically due to usage and may eventually break down, resulting in high maintenance costs at scheduled maintenance moments. This deterioration behavior is affected by the system's production rate. While producing at a higher rate generates more revenue, the system may also deteriorate faster. Production should thus be controlled dynamically to trade-off deterioration and revenue accumulation in between maintenance moments. 
We study systems for which the relation between production and deterioration is known and the same for each system as well as systems for which this relation differs from system to system and needs to be learned on-the-fly.
The decision problem is to find the optimal production policy given planned maintenance moments (operational) and the optimal interval length between such maintenance moments (tactical). 
\textbf{Methodology/Results:} For systems with a known production-deterioration relation, we cast the operational decision problem as a continuous-time Markov decision process and prove that the optimal policy has intuitive monotonic properties. We also present sufficient conditions for the optimality of bang-bang policies and we partially characterize the structure of the optimal interval length, thereby enabling efficient joint optimization of the operational and tactical decision problem.
For systems that exhibit variability in their production-deterioration relations, we propose a Bayesian procedure to learn the unknown deterioration rate under any production policy.
Numerical studies indicate that on average across a wide range of settings (i) condition-based production increases profits by 50\% compared to static production, (ii) integrating condition-based production and maintenance decisions increases profits by 21\% compared to the state-of-the-art sequential approach, and (iii) our Bayesian approach performs close -- especially in the bang-bang regime -- to an Oracle policy that knows each system's production-deterioration relation. 
\textbf{Managerial Implications:}  
Production should be adjusted dynamically based on real-time condition monitoring and the tactical maintenance planning should anticipate and integrate these operational decisions. Our proposed framework assists managers to do so optimally. 
}%



\maketitle

%

\vspace{-20pt}
\section{Introduction}
The modern manufacturing industry relies heavily and increasingly on advanced technical systems in their daily operations. These systems -- for instance lithography systems in the semiconductor industry, wind turbines in wind farms, and assembly equipment in production lines -- deteriorate due to usage and may break down as a result of such deterioration, leading to extremely high costs in the event of prolonged downtime. For example, unplanned downtime of lithography equipment can cost up to 72,000 euros per hour \citep{lamghari2021new}. The global manufacturing industry, which includes the semiconductor industry, faces an estimated annual cost of around \$50 billion due to unplanned system downtimes and their consequences \citep{deloitte2017}. 

Motivated by such cost figures, there has been a considerable amount of research on maintenance optimization, where the main goal is to determine when to preventively maintain a system based on its condition in order to prevent costly failures and their associated downtimes (see \cite{Wang2002} and \cite{de2019review} for excellent overviews of the field). A common assumption in this field is that the underlying deterioration of a system is exogenous; that is, the deterioration is not influenced by factors other than by its own randomness. However, for many advanced technical systems in the manufacturing industry, the deterioration behavior of its components is often affected by the production load applied to the system itself. For example, gearboxes and generators in a wind turbine deteriorate faster under higher rotational speeds \citep{feng2013monitoring, macquart2019stall} which can be adjusted instantly and remotely by controlling the wind turbine's pitch \citep{spinato2009reliability,zhang2015condition}. Other examples include cutting tools in high-speed cutting machinery that wear faster under higher speeds \citep{dolinvsek2001wear}, conveyor belts in production lines that are increasingly likely to fail as their rotational speeds are increased \citep{nourelfath2012integrated}, stamping machines that deteriorate faster as the stamping speed increases \citep{hao2015}, and rapid gravity filters in waterworks that deteriorate faster as the filtration rate is increased \citep{sun2019reliability}. 
There is a caveat to this: Production loads that lead to less deterioration accumulation typically also come with less revenue accumulation. Indeed,  wind turbines generate less wind energy when their rotational speeds are lower. Production lines and cutting tools have less output when conveyor belts and cutting speeds are slowed down, respectively. Production decisions in make-to-stock environments are not driven by order deadlines but by revenue targets and equipment maintenance requirements. It is in these environments that our contribution will apply best.


Systems usually fail when components deteriorate beyond a critical level. As such, failures can be prevented before they occur by making instant maintenance decisions based on the component's deterioration level. (Hence most research has focused on the decision when to preventively maintain a system based on its condition.) However, in practice, maintenance planning has limited flexibility and cannot be done last minute because arranging the logistic support takes time, particularly for remotely located systems (e.g., wind turbines in offshore wind farms). 
Moreover, it is inconvenient and costly to shut down an entire system each time an individual component needs maintenance. Maintenance intervals across a set of assets/components are therefore often fixed to reduce the fixed cost of shutting down an asset and bringing all the maintenance resources together. The classical example of this is the block maintenance policy studied by \cite{barlow1960optimum}, forming the basis for more recent work, such as that of \cite{vanDijkhuizenvanHarten1997}, \cite{ArtsBasten2018}, and \cite{uit2020condition}. For these reasons, in practice, maintenance activities are usually performed at planned maintenance moments that are scheduled periodically and well in advance. At such planned maintenance moments, the entire system is deliberately shut down and components that require maintenance are maintained.

When these planned maintenance moments are scheduled in advance, another strategy that decision makers (DMs) can employ to reduce maintenance costs and increase revenues, is to control the deterioration behavior in between such maintenance moments. Specifically, slowing down the deterioration rate by decreasing the production rate when a component's deterioration is close to its critical level leads to fewer failures between consecutive maintenance moments. On the other hand, if components deteriorate slowly, DMs can increase the production rate, thereby accumulating a higher revenue without running the risk of system failure before the next maintenance moment.

This operational strategy has become a viable option due to the increasing and improved availability of real-time data stemming from ubiquitous sensors installed in modern equipment, integrated in the so-called Internet-of-Things \citep{PWC2014, MCk2017}. It particularly enables  (i) remote accrual of deterioration data, and (ii) continuous remote control of production operations at near-zero costs and in real time \citep{olsen2020industry}. 
This novel idea was recently proposed by \cite{uit2020condition}, who coined the term condition-based production to describe this strategy. The authors consider condition-based production rate decisions for systems with an adjustable production rate that directly affects the deterioration rate. There is an inherent trade-off to this concept: A lower production rate has both a lower deterioration and revenue rate, while a higher production rate has both a higher deterioration and revenue rate. By formulating a continuous-time optimal control problem, \cite{uit2020condition} are able to study this trade-off and establish structural insights for the case that deterioration is deterministic.
As yet, there are no analytical results for stochastic deterioration processes despite significant practical relevance.

In this paper, we study such optimal condition-based production policies. Contrary to the existing literature, we consider stochastically deteriorating systems, and our objective is to characterize the structure of the optimal policy. 
To this end, we model the stochastic deterioration as a Poisson process whose uncontrolled intensity is equal to a certain base rate. 
Modeling deterioration as a Poisson process is a common and empirically validated approach in the literature for systems subject to so-called stochastic shock deterioration, where  systems only deteriorate due to shocks that occur at random points in time \citep[][]{esary1973shock,sobczyk1987stochastic}. The base rate represents the rate at which the system deteriorates under normal, uncontrolled operating conditions. A DM can continuously adjust production which leads to instantaneous changes in the revenue rate and the rate of the Poisson process (proportional to the base rate). We shall discuss this in more detail in \S \ref{ProbForm}, but conceptually, we have that a higher (lower) production rate leads to higher (lower) revenue and deterioration rates. By employing such a condition-based production policy, the DM can steer the deterioration as she nears the planned maintenance moment, at which a maintenance cost depending on the deterioration state is incurred.  

We first consider the case where each newly installed system has the same fixed base rate that is a-priori known to the DM. We model this problem as a continuous-time Markov decision process (MDP) when the planned maintenance interval is fixed and rigorously analyze its Hamilton-Jacobi-Bellman (HJB) equations to characterize the monotonic behavior of the optimal policy. We also provide easily verifiable conditions that are sufficient for a so-called bang-bang policy (production is either set to the maximum or off) to be optimal, thereby reducing the action space of the continuous-time MDP from a compact set with infinitely many actions to a discrete set with only two actions. We then study the case where each newly installed system has a different base rate that is unknown a-priori and must be learned in real time, that is, the population of systems exhibits unit-to-unit variability in their base rates. For this case, we develop an easy-to-implement Bayesian heuristic, that mimics the structure of the optimal policy for the known base rate case. In this scenario, a DM needs to learn the a-priori unknown base rate based on the chosen production policy and the deterioration she observes as a consequence of that policy. 

Finally, we explore the tactical trade-off due to maintenance costs levied at a planned maintenance moment. Specifically, scheduling planned maintenance moments too frequently leads to unnecessary maintenance activities and associated costs, while a very long maintenance interval both increases the risk of a failed system and decreases the revenue accumulation as the DM is then forced to decrease the production rate. This implies that in deciding upon the length of the interval between planned maintenance moments (a tactical decision), the optimal condition-based production policy (operational decisions) needs to be taken into account. We study this trade-off and show that the optimization of the maintenance interval is a unimodal optimization problem.

The main contributions of this paper are fourfold: 
\begin{enumerate} 
\item We are the first to analytically establish the structure of the optimal condition-based production policy under stochastic deterioration. The structure of the optimal policy has intuitive appeal: Decrease production if the system's condition is close to failure, and increase production if the planned maintenance moment is nearing or when the base rate is lower. We also establish sufficient conditions for when the optimal policy has a bang-bang structure. As a byproduct of our analysis, we obtain new submodularity preservation properties that may have merit themselves in future work in the broad field of monotonic optimal control theory.  
\item We show that under the optimal condition-based production policy, the length of the interval between planned maintenance moments can be easily optimized. This allows for the optimal integration of condition-based production on an operational level and the length of planned maintenance intervals on a tactical level. 
\item We complement our theoretical results with an extensive numerical study, which demonstrates that on average across a wide range of practical settings: (i) condition-based production policies lead to 50\% profit increases compared to optimal static policies, and (ii) integrating maintenance decisions and condition-based production to 21\% profit increases compared to the state-of-the-art sequential approach. Further comparative static analyses offer insights into the driving factors behind these profit increases, and consequently highlight practical settings that are associated with higher or lower increases.  
\item For systems that exhibit unit-to-unit variability in their base rates, we propose a Bayesian procedure to tractably learn this unknown rate under any production policy, and we use it to build a heuristic policy. In an extensive simulation study, we show that this heuristic performs close to an Oracle policy that knows the base rate a-priori. Notably, our results highlight a key insight: The heuristic yields particularly good performance in the bang-bang regime, and this performance remains robust against unit-to-unit variability.
\end{enumerate}
The rest of the paper is structured as follows. In \S \ref{litReview}, we review related literature. We present the problem formulation in \S \ref{ProbForm}. In \S \ref{cbpStructural}, we present the main analytical results for populations of systems whose base rates are fixed and known. \S \ref{sec:num} presents a numerical study that underscores the practical value of our theoretical results. In \S\ref{sec:heterogenous}, we study system populations with unit-to-unit variability in their base rates. Finally, \S \ref{sec:conclusion} provides concluding remarks. All proofs are relegated to Appendix \ref{AppendixProofs}. 

\section{Literature Review}\label{litReview}
Studies on production control under uncertainty as well as on maintenance optimization are abundant in literature. \cite{sethi2002optimal} and \cite{MULA2006271} provide excellent overviews of the former field while \cite{ALASWAD201754} and \cite{de2019review} do so for the latter.
The intersection of these two fields has received considerably less attention in the literature. In what follows we discuss studies on this intersection that are most relevant to our work. 
Additionally, since our methodology exhibits similarities with modeling approaches in revenue management, we conclude by discussing and contrasting our work with this literature stream too.

Most of the early works focus on production dependent failures. In such settings, production equipment resides in either a failed or non-failed state, and the probability of failing depends on the production rate. 
\cite{liberopoulos1994production} were among the first to study such a setting. 
They consider a make-to-stock production system where the system's lifetime has an exponential distribution whose rate depends on the production rate. Upon failure, corrective maintenance is needed, which also takes an exponentially distributed amount of time.
\cite{liberopoulos1994production} focus on the cost performance evaluation of so-called threshold polices. Under such a policy, the production rate is governed by thresholds on the inventory level, where higher inventory levels prescribe lower production rates.
Many follow-up studies that depart from this work have appeared since then.
While some authors have focused on providing conditions under which certain threshold policies are optimal, typically dependent on the production-failure relationship \citep[e.g.,][]{hu1995monotonicity, martinelli2010manufacturing}, others have extended the model by incorporating the option to perform preventive maintenance that is cheaper than corrective maintenance \citep[e.g.,][]{boukas1995robust, boukas2001production}, by considering multiple production machines \citep[e.g.,][]{francie2014stochastic}, or by letting the failure rate depend on other variables, for instance the equipment's age \citep[e.g.,][]{zied2011optimal}.
Although there are practical settings where production equipment suddenly ceases to operate from a normal operating state, 
most production equipment deteriorates gradually during production.
For that reason, research has focused on extending the research line of production dependent failures to models in which the deterioration of the equipment can be described by more than one non-failed state, leading to so-called production dependent deterioration.   

Early works on production dependent deterioration can be found in \cite{conrad1987drilling} and \cite{rishel1989controlled}. \cite{conrad1987drilling} study a drilling machine with controllable speed whose condition deteriorates according to a speed-dependent Brownian motion. \cite{rishel1989controlled} studies a similar problem, but considers a wear process that is monotone. Similar to our paper, both works let the deterioration and revenue rate depend on the control variable, but unlike our work, they incorporate a deterioration-dependent cost rate and wish to maximize the expected profit over the system's  random lifespan. These works assume that maintenance can be performed instantaneously at any time. This assumption is often not tenable in practice since maintenance -- especially for remotely located systems -- is usually planned well in advance and therefore permits limited flexibility, if any. 
In our work, we take this into account and consider the case where maintenance can only be performed at pre-specified planned times; we then wish to maximize the expected profit over such a planning horizon. 
Moreover, \cite{conrad1987drilling} and \cite{rishel1989controlled} do not provide analytical results regarding the structure of the optimal policy. Since these early works, researchers have extended production dependent deterioration models with the inclusion of maintenance activities to improve the equipment's condition \citep[e.g.,][]{sloan2000combined,sloan2002using,kazaz2008production,batun2012reassessing,kazaz2013impact,uit2021joint}, but again with few explicit structural results.

The work on production dependent deterioration that is most relevant to ours can be found in \cite{uit2020condition}. Like us, they study how to dynamically adjust the production rate over a finite planning horizon that ends with a maintenance moment, where a higher production rate leads to more revenue, but also to a faster deterioration rate, and maintenance costs at the end of the planning horizon. The authors formulate an optimal control problem and characterize the structure of the optimal condition-based production policy under the assumption that deterioration is deterministic and its parameters fully known. 
When deterioration is stochastic, the authors numerically explore the optimal policy but do not report on structural results. Our work extends the work of \cite{uit2020condition} in three important ways. Firstly, we establish the structural properties of the optimal condition-based production policy for stochastically deteriorating systems; these properties are found to differ significantly from those under deterministic deterioration. Secondly, we study the case in which each system's deterioration process is characterized by an a-priori unknown parameter that must be learned on-the-fly, which gears our model even more towards practice. Indeed, recent practice-based research has shown that deterioration processes of the type studied in this work may exhibit unit-to-unit variability \citep{drent2021cbm}, meaning that the parameters of the deterioration process of each newly installed system differ from previous ones.  Finally, we show that the interval length between planned maintenance moments can be easily optimized under the optimal condition-based production policy. This allows for the optimal integration of condition-based production decisions on an operational level and maintenance decisions on a tactical level, which has not been studied before.  

Methodologically, our modeling approach shares some similarities with approaches to model the canonical dynamic pricing problem in revenue management, where a DM seeks to maximize revenue by pricing a finite number of products within a fixed time frame. In this problem, introduced by \cite{gallego1994optimal}, customers arrive according to a non-homogeneous Poisson process whose intensity depends on the chosen price. They formulate this problem as an intensity control problem, also leading to a set of HJB equations. Their work instigated a large stream of literature in which various extensions are studied, such as multiple products \citep{gallego1997multiproduct}, unknown market responses \citep{araman2009dynamic,farias2010dynamic}, multiple offered prices \citep{aydin2009personalized}, time-varying demand \citep{cao2012optimal}, discounts and rebates \citep{aydin2008pricing,hu2017optimal}, and sales target constraints \citep{du2021sales}. In these works, demand for products is modeled as a controlled Poisson process. As such, they too model their problem as Poisson intensity control problems and subsequently study the corresponding HJB equations to establish structural properties of the optimal policy. 

There are two fundamental differences between our formulation and the ones encountered in revenue management (see \cite{gallego2019revenue} for a complete overview of the field). First, a chosen price impacts both the revenue accumulation and demand process in a different way than a chosen production rate impacts the revenue accumulation and deterioration process in our problem. In dynamic pricing, revenue only accumulates at discrete demand arrival epochs, while in our problem, revenue is accumulated continuously. In the former, choosing a higher price leads to a lower demand rate and potentially more revenue if demand occurs, while in our problem, choosing a higher production rate induces a higher deterioration rate and instantaneously a higher revenue rate. Second, we deal with a boundary condition that models the maintenance cost levied at the end of the planning horizon. This boundary condition needs to be taken into account when adopting a production policy throughout the planning horizon to prevent high maintenance costs at the end. In the revenue management literature, such boundary conditions typically do not exist.

\section{Problem Formulation}\label{ProbForm}
We consider a single-unit production system whose deterioration can be represented by a single parameter. The deterioration process of the system is continuously monitored and described by a Poisson process $Y=(Y_t)_{t\geq 0}$, with $Y_0=0$ because at each planned maintenance moment a new unit is installed and hence has zero deterioration. Poisson deterioration is a common modeling choice in the OM literature, for recent examples, see, e.g., \cite{tian2021optimal} and \cite{li2022after}. The system fails when the deterioration level is equal to a fixed failure level $\xi\in \N$. This is in line with the prevailing maintenance literature and practice, see, e.g., \cite{Wang2002}, \cite{gebraeel2005residual}, \cite{elwanyetal2011}, \cite{uit2020condition}, and \cite{drent2021cbm}, and the case studies therein. We denote the set of all deterioration levels with $\mathscr{X} \triangleq \{0,1,\ldots,\xi\}$. (Discrete deterioration states are common, e.g. \cite{derman1963optimal,kim2013joint,tian2021optimal,li2022after,drent2021cbm} and \cite{zhang2023analytical}.) 

The DM can adjust a univariate production control variable of the system that will both impact the deterioration process and the rate at which the production system generates revenue. We remark that in the remainder of the paper, we will refer to this variable as the production rate, though in practice, it could represent any adjustable setting or control variable of the system that has the said impact. Hence, if the DM selects production rate $s \in \mathscr{S} \triangleq [0,\bar{s}]$ (where $0<\bar{s}<\infty$ is the maximum production rate), then the rate of the Poisson process is instantaneously governed by the base rate, denoted with $\lambda$, multiplied by the deterioration function $f: \mathscr{S} \rightarrow \R_+$, and the revenue rate is instantaneously governed by the revenue function $r: \mathscr{S} \rightarrow \R$. In practice, one can view the base rate as the base line for deterioration under normal operating conditions.

For now we shall assume that the base rate $\lambda$ is fully known to the DM and that each system has the same fixed base rate upon replacement at the scheduled maintenance moment. The results for such systems with fixed base rates also help us to develop and benchmark a heuristic for systems for which the base rates are unknown a-priori to the decision-maker and differ from system to system. We discuss such populations of systems with unit-to-unit variability in detail in \S \ref{sec:heterogenous}. 

Throughout this paper, we use increasing and decreasing in the weak sense. We impose the following assumption on the deterioration and revenue function.
\begin{assumption}\label{assumptionFunctions}
The deterioration function $f(\cdot)$ and revenue function $r(\cdot)$ are both continuous and increasing. Furthermore, $f(0)=0$ and $r(0)\leq 0$.
\end{assumption}
Assumption \ref{assumptionFunctions} together with the compactness of $\mathscr{S}$ also guarantees that  $f(\cdot)$ and $r(\cdot)$ are bounded. The increasing assumption is reasonable in practice as higher production rates generally lead to higher revenue and deterioration rates. Two remarks regarding $r(\cdot)$ are in order here.
\begin{remark}\label{remarkRev}
Our analysis does not require $r(\cdot)$ to be increasing: All results also hold if $r(\cdot)$ is decreasing. However, if $r(\cdot)$ is decreasing then this will result in a practically meaningless policy in which the system is always turned off as turning it on leads to a decrease in revenues.
\end{remark}
\begin{remark}\label{remarkRev2}
In settings where the focus is to maximize production revenues, it is natural to have $r(0)=0$ (no revenue is earned when the system is off), though there are situations where $r(0)<0$. One practical setting is when a DM aims to satisfy a (continuous) demand rate, say $D>0$, with the production  rate $s$. This can be modeled by choosing $r(\cdot)$ such that the downfall $(D-s)^+$ is penalized and the surplus $(s-D)^+$ is rewarded, where $(x)^+\triangleq\max(0,x)$, so that $r(0)$ is negative. One notable example is $r(s)\triangleq b(s-D)^+-p(D-s)^+$, where $b$ $(p)$ is the bonus (penalty) rate. One extension to this case, which we discuss in Appendix \ref{extensionMulti}, is when a single demand rate $D$ can be met by the cumulative production rate of multiple systems.
\end{remark}

Maintenance is only performed at scheduled maintenance moments, and the next maintenance action is scheduled at time $T$ (e.g., maintenance teams visit offshore wind farms semi-annually, or the system undergoes intentional shutdowns yearly). 
The maintenance cost at scheduled maintenance moments is governed by an increasing convex maintenance cost function $c_m : \mathscr{X} \rightarrow \R_+$, which reflects that higher deterioration levels lead to greater adverse effects, and these effects worsen as deterioration progresses. The canonical maintenance cost structure with a single corrective ($c_u$) and single preventive maintenance cost ($c_p<c_u$) belongs to this class of maintenance cost functions. 

When maintenance moments are scheduled, the DM's operational problem is to dynamically adjust the production rate with the goal of maximizing the expected profit (i.e. revenue minus costs). We formulate this optimization problem with time indexed backwards. 
To this end, let $X_t = \{Y_{\tilde{t}}\ |\ 0 \leq \tilde{t}\leq T \}$ be the deterioration level at time $t=T-\tilde{t}$ before the next maintenance moment. Since $Y_{\tilde{t}}$ is the number of deterioration increments up to time $\tilde{t}$, $dY_{\tilde{t}}=1$ if deterioration increases at time $\tilde{t}$. 
Let $\mathscr{U}$ be the set of all admissible non-anticipating production policies satisfying $\int_{k=0}^t d Y_k\leq \xi - X_t$ (almost surely). This constraint ensures that whenever the system fails (i.e. when $Y_t = \xi$), the system is turned off. It can always be satisfied because $s=0\in\mathscr{S}$ with $f(0)=0$ turns the system off. For a given production policy $u\in \mathscr{U}$, let $s^u_{\lambda}(x,k)$ be the induced production rate if the deterioration level is $x$, the time to planned maintenance moment is $k\in [t,0)$ and the base rate is $\lambda$. The expected profit, starting from state $(x,t)$ under a policy $u\in \mathscr{U}$ is then denoted by 
\begin{align}\label{profitUnderPolicy}
J_{\lambda}^{u}(x,t)  \triangleq  \E_u \Bigg[\int_0^t r\big(s^u_{\lambda}(X_k,k)\big) dk -c_m \big( X_0 \big) \ \bigg|\  X_t = x \Bigg],
\end{align}
where $\E_u$ denotes that the expectation is taken with respect to policy $u$. The integral part in \eqref{profitUnderPolicy} represents the expected revenue, while the second part represents the expected maintenance costs under policy $u$. We make the dependence of the total expected profit on the base rate $\lambda$ explicit by using the subscript, thereby establishing the correspondence to a population of systems with the same fixed rate. The goal is to find the optimal condition-based production policy, denoted with $u_{\lambda}^*$, that maximizes the expected profit generated over $[t,0]$, denoted with
$J_{\lambda}^{*}(x,t)  \triangleq  \sup_{u\in\mathscr{U}} J_{\lambda}^{u}(x,t).
$ For brevity, we introduce the shorthand notation $\Delta_{\lambda}(x,t) \triangleq J_{\lambda}^*(x,t) - J_{\lambda}^*(x+1,t)$, which can be interpreted as the marginal value of less deterioration for fixed $t$ and $\lambda$ under the optimal policy. Using this additional notation, we are in the position to present our first result in Lemma \ref{lemmaHJB} below: The HJB equation for the optimal expected profit $J_{\lambda}^*(x,t)$.
\begin{lemma}\label{lemmaHJB}
$J_{\lambda}^*(x,t)$ is the solution of 
\begin{align}\label{HJEqn}
\frac{\partial J_{\lambda}^*(x,t)}{\partial t} =  \begin{cases}
			\max_{s\in \mathscr{S}} \Big[r(s) -  \lambda f(s)\Delta_{\lambda}(x,t) \Big], 
& \text{if $x \in \{0,1,\ldots, \xi-1\},$}\\
            0, & \text{if $x=\xi$},
		 \end{cases}
\end{align}
with boundary conditions $J_{\lambda}^*(x,0) \triangleq  -c_m(x)$ for all $x\in \mathscr{X}$. 
\end{lemma}
The formulation in Lemma \ref{lemmaHJB} is a Markovian intensity control problem; we refer the reader to \cite{bremaud1980} for an excellent treatment of the general theory on this subject. Since $f(\cdot)$ and $r(\cdot)$ are continuous and bounded (by Assumption \ref{assumptionFunctions}), and $\mathscr{S}$ is a compact subset of $\R$, we are guaranteed by Theorem III of \cite{bremaud1980} that there exists a solution to \eqref{HJEqn}, which we denote by $s_{\lambda}^*(x,t)$; that is, $s_{\lambda}^*(x,t)$ is the optimal production rate for state $(x,t)$ when the base rate is $\lambda$. 

\section{Fixed Base Rate: Analytical Results}\label{cbpStructural}
This section contains our analytical results for systems with a known, fixed base rate. We first present structural properties of the optimal expected profit and optimal condition-based production policy for a given maintenance interval length $T$. We then study the tactical problem of scheduling maintenance moments, given that the optimal production policy is used in between such moments.

\subsection{The Optimal Condition-Based Production Policy} 
We first show that the optimal expected profit function $J_{\lambda}^{*}(x,t)$ behaves monotonically in its state variables and the base rate, which we then use to prove that the optimal production rate also has certain monotonicity properties. We illustrate these theoretical results with numerical examples.

Using sample path arguments, one can deduce first-order monotonic properties of $J_{\lambda}^{*}(x,t)$ with respect to $x$, $t$, and $\lambda$ that are intuitive and useful for subsequent analysis. These properties are summarized in the following lemma. 

\begin{lemma}\label{lemmaIncreasing}
The optimal expected profit $J_{\lambda}^*(x,t)$ is decreasing in the deterioration level $x$, increasing in the time to the next maintenance moment $t$, and decreasing in the base rate $\lambda$. 
\end{lemma}

Everything else fixed, Lemma \ref{lemmaIncreasing} states that for a more deteriorated system, the expected profit will be lower than for a system that is less deteriorated. Similarly, everything else fixed, the expected profit is higher when the time to planned maintenance is higher, because there is more time to generate revenue. Finally, a system whose base rate is higher, deteriorates faster under any production rate, which in turn leads to a lower expected profit.   

We now turn our attention to higher-order monotonic properties of the expected profit with respect to its parameters, which are needed to prove structural properties of the optimal policy. 
 \begin{proposition}\label{concaveInX}
The optimal expected profit $J_{\lambda}^*(x,t)$  is concave in $x$ for all $t\geq0$ and $\lambda\geq0$. 
\end{proposition}
The concavity of the optimal expected profit in the deterioration level implies that expected profit drops ever steeper as the system deterioration level is closer to the failure level for any fixed time.
This is intuitive: A system that is approaching the failure threshold has less and less room (in terms of remaining useful condition) to control the deterioration, and hence for generating revenue.     

The next result establishes that the marginal value of less deterioration, increases in the remaining time until the next planned maintenance. 
\begin{proposition}\label{increasingInT}
$\Delta_{\lambda}(x,t)$ is increasing in $t$ for all $x\in\mathscr{X}$ and $\lambda\geq0$. 
\end{proposition}
The intuition behind Proposition \ref{increasingInT} is as follows. Consider a system with a better condition (i.e. with one deterioration level less) than another system. The former generates more revenue than the latter, see Lemma \ref{lemmaIncreasing}. Now, if the scheduled maintenance moment is further away, then for the former system, there is a longer planning horizon to exploit this better condition -- for generating revenue -- than when there is less time to do so. This leads to a higher additional revenue.    

Using a similar intuition one can also reason that the optimal expected profit must be concave in the time to the next scheduled maintenance moment for all deterioration levels and each base rate. This is indeed true and we present this result in the next proposition. 

\begin{proposition}\label{concaveInT}
The optimal expected profit $J_{\lambda}^*(x,t)$  is concave in $t$ for all $x\in\mathscr{X}$ and $\lambda\geq0$. 
\end{proposition}

We now proceed with proving that the optimal expected profit is submodular in $(x, \lambda)$, which subsequently enables us to establish the monotonic behavior of the optimal policy in the base rate $\lambda$ in light of the well-known Topkis' Theorem \citep[see, for instance,][Theorem 2.8.1]{topkis2011supermodularity}. We first state two technical results that relate to preservation of submodularity. These results are preceded by a formal definition of submodularity.

\begin{definition}[Submodularity]\label{submodular}
A bivariate function, $g(x,y)$, is submodular in $(x,y)$ if for all $x^+\geq x^-$ and $y^+\geq y^-$, we have
$\big(g(x^+,y^-)-g(x^-,y^-)  \big)  \geq  \big( g(x^+,y^+) - g(x^-,y^+)  \big).$
\end{definition}

\begin{lemma}\label{submodFunction}
Suppose $g(y)$ is nonnegative, and increasing in $y$, and $f(x,y)$ is nonnegative,  submodular in $(x,y)$, and  decreasing in $x$, then $g(y)f(x,y)$ is submodular in $(x,y)$.
\end{lemma}
\begin{lemma}\label{submodExpect}
Suppose $f(x,\theta)$ is submodular in $(x,\theta)$ and concave in $x$ and the nonnegative random variable $Z(\theta)$ is stochastically increasing in the usual stochastic order in $\theta$, then $\E[f(x+Z(\theta),\theta)]$ is submodular in $(x,\theta)$ and concave in $x$. 
\end{lemma}

Lemma \ref{submodFunction} and \ref{submodExpect} complement other preservation results -- usually needed for establishing structural properties in optimization problems -- that can be found in the literature; for instance a closely related result by \cite{chao2009dynamic} that focuses on $\E[f\left(x,Z(\theta)\right)]$ (in our notation). Using these two lemma's, we are able to establish the submodularity of the  optimal expected profit in $(x, \lambda)$, which is presented in the result below. 

\begin{proposition}\label{submodvalue}
For all $t \geq 0$, the optimal expected profit  $J_{\lambda}(x,t)$ is submodular in $(x,\lambda)$.
\end{proposition}
Definition \ref{submodular} shows that submodularity of a bivariate function $g(x,y)$ in $(x,y)$ implies that the difference $g(x^+,y) - g(x^-,y)$ is decreasing with respect to $y$. Using this definition, Proposition \ref{submodvalue} thus implies that the decrease in the optimal expected profit due to more deterioration, i.e. $J_{\lambda}(x^+,t)-J_{\lambda}(x^-,t)$ with $x^+ \geq x^-$, is decreasing in the base rate. In other words, the decrease in the optimal expected profit due to deterioration accumulation is smaller for systems facing a lower base rate than for systems facing a higher base rate.  

The next result is the main result of this section; it characterizes the monotonic behavior of the optimal production rate.

\begin{theorem}\label{optPolicy} The optimal production rate when the deterioration level is $x$, the time until planned maintenance is $t$ and the base rate is $\lambda$, i.e. $s_{\lambda}^*(x,t)$, is decreasing in $x$, $t$, and $\lambda$.   
\end{theorem}

The intuition behind Theorem \ref{optPolicy} is as follows.
\begin{enumerate}
\item \textbf{Monotonicity in $x$:} Given the same time until planned maintenance, production will be slowed down whenever failure is nearing to avoid high maintenance costs. 
Likewise, a DM will select a higher production rate for a system that is less deteriorated than for a system that is close to failure. In doing so, she increases revenues without running a high risk of failure.
\item \textbf{Monotonicity in $t$:} Given the same deterioration level, a DM is more conservative at the beginning of the planning horizon and chooses a lower production rate than when the scheduled maintenance moment is nearby. Similarly, when the planned maintenance moment is nearing and the deterioration level is kept constant, the DM will exploit the remaining useful condition for accumulating revenue by choosing an increasingly higher production rate.
\item \textbf{Monotonicity in $\lambda$:} If a system deteriorates faster under any production rate through a higher base rate, then a DM will choose lower production rates for all deterioration levels and remaining times than for a system that deteriorates more slowly. This ensures that in the former, high maintenance costs are prevented, while in the latter, more revenue is generated.  
\end{enumerate}
Notably, the structure of the optimal condition-based production policy is inherently different from the optimal policy under deterministic deterioration. For the latter, \cite{uit2020condition} prove that it may be optimal to purposely break down the system when a failure is unavoidable, breaking down any monotonic properties of the optimal policy. This sharply contrasts with our result, where a DM would gradually decrease production when a potential failure is nearing. However, Theorem \ref{optPolicy} does formally confirm some of their reported intuitions -- based on numerical explorations -- of the structure of the optimal policy for the stochastic setting \citep[see][\S 5.3]{uit2020condition}. 

\begin{figure}[htb!]
\begin{center}
\includegraphics[width=0.75\textwidth]{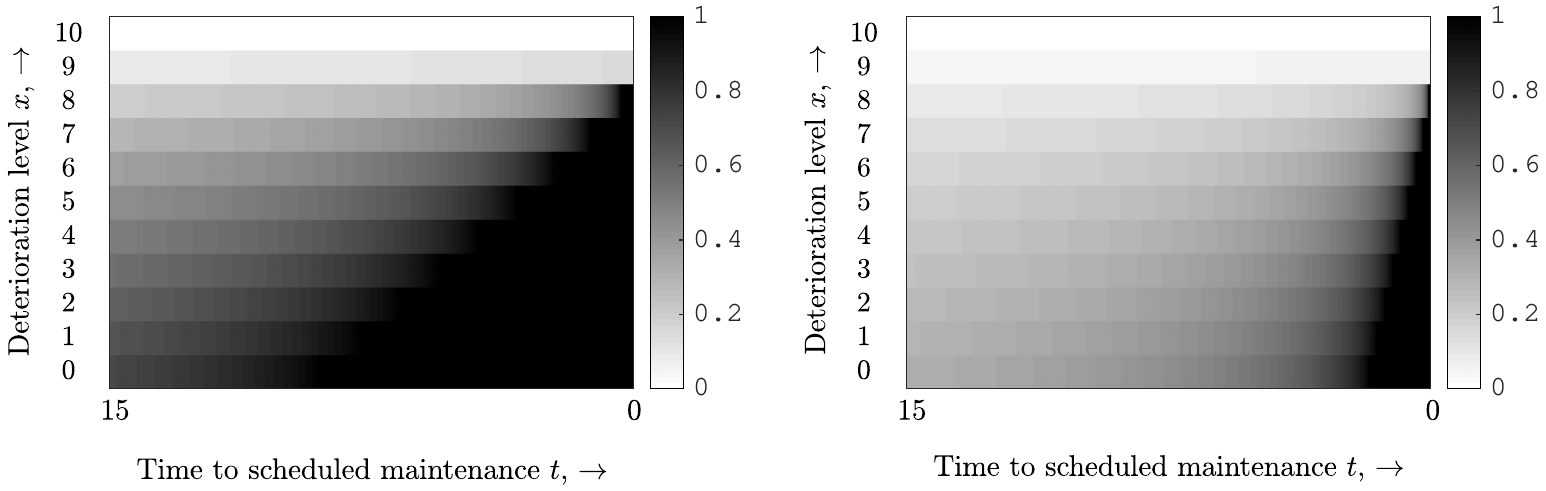}
\caption{Illustration of monotonic behavior of the optimal condition-based production policy with base rate $\lambda = 1$ (left) and base rate $\lambda =4$ (right). The color-bar indicates the optimal production rate.}
\label{optimalPolicyIllus}
\end{center}
\end{figure}

\begin{figure}[htb!]
\begin{center}
\includegraphics[width=0.625\textwidth]{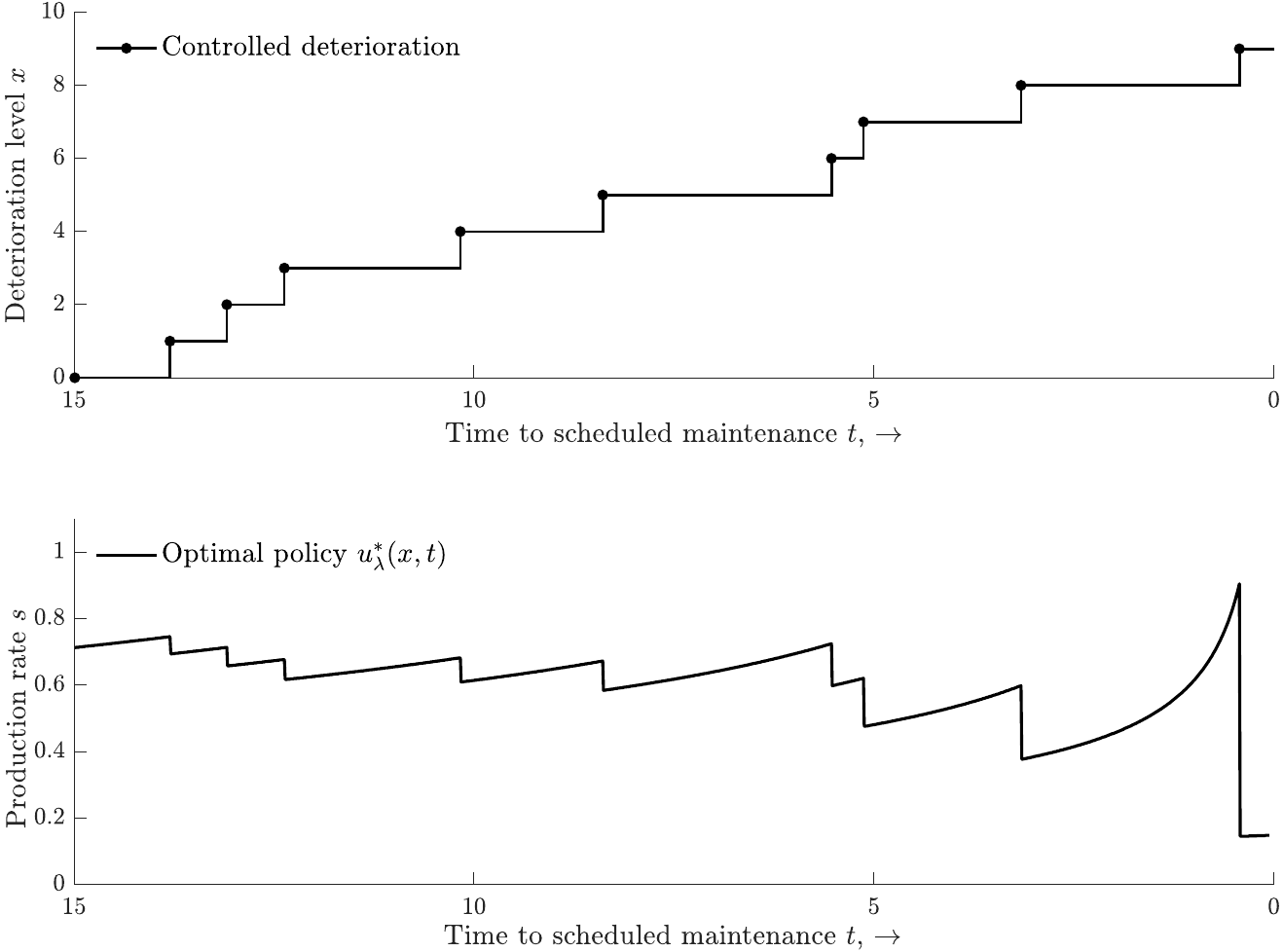}
\caption{Sample paths of the optimal condition-based production policy (bottom) and controlled deterioration process (top), when the base rate $\lambda$ is 1.}
\label{samplepathOptimal}
\end{center}
\end{figure}

The monotonic behavior of the optimal condition-based production policy is illustrated by Figures \ref{optimalPolicyIllus} and \ref{samplepathOptimal}. In both figures, we use $\xi=10$ (failure threshold), $c_m(\xi)=5$ (corrective maintenance cost), and $c_m(\cdot)=1$ for all non-failed states (preventive maintenance cost). The planning horizon $T$ is equal to 15 and the production rate space is equal to $\mathscr{S}=[0,1]$. For the revenue and deterioration functions, we use $r(s)=s^{1/2}$ and $f(s)=s^{2}$. We obtained the optimal policy in these illustrations, and elsewhere in this paper, using the standard approach for solving continuous time dynamic programs, based on the finite difference formulation (see Appendix \ref{appendixProcedure} for more details). 

In the left (right) subfigure of Figure \ref{optimalPolicyIllus} the base rate is 1 (4). In each subfigure, we clearly observe that the optimal production rate increases as the time to scheduled maintenance decreases for a fixed deterioration level (moving right horizontally), and that the optimal production rate decreases as the condition deteriorates (moving up vertically). If we compare the left with the right subfigure, we clearly see that lower rates are prescribed in all states $(x,t)$ in the right subfigure. 

Figure \ref{samplepathOptimal} presents sample paths of the optimal condition-based production policy (bottom) and the corresponding controlled deterioration (top). We see that when the system's condition deteriorates, production is immediately lowered, which  highlights the monotonicity of the optimal production rate in the  deterioration level for fixed time to scheduled maintenance. Also, the gradual rise of the production rate between deterioration increments highlights the monotonicity of the optimal production rate in the time to scheduled maintenance for a fixed deterioration level. 

\subsubsection{Optimality of Bang-Bang Production Policies}
In this subsection, we present sufficient conditions for the optimal policy to be a so-called bang-bang policy. See Definition \ref{defbangbang} below for a formal description of a bang-bang policy. 
\begin{definition}[Bang-bang policy] \label{defbangbang}
    A bang-bang production policy is a policy in which production is either set to the maximum rate $s=\bar{s}$ (full production) or turned off $s=0$ (no production). 
\end{definition}

The next result provides a sufficient condition in terms of $r(\cdot)$, $f(\cdot)$, and their maxima, for a bang-bang policy to be optimal. Note that this condition does not depend on the base rate $\lambda$.
\begin{proposition}\label{bangbangCondition}
The optimal policy has a bang-bang structure if
\begin{align}
\frac{r(s)}{f(s)} \leq \frac{r(\bar{s})}{f(\bar{s})}\mbox{ for all 
 } s\in (0,\bar{s}). \label{bangcondition}
 \end{align}
\end{proposition} 

The next result, which follows from Proposition \ref{bangbangCondition}, provides sufficient conditions only on the functions $r(\cdot)$ and $f(\cdot)$ for a bang-bang policy to be optimal that can be verified without verifying whether \eqref{bangcondition} holds for all  $s\in (0,\bar{s})$.  
\begin{corollary}
\label{extremalConvexConcave}
The optimal policy has a bang-bang structure if any of the following conditions is satisfied:
\begin{enumerate}[label=(\roman*)]
    \item $r(s)=f(s)$ for all $s\in\mathcal{S}$.
    \item $r(s)= a\cdot s^\alpha$ and $f(s)= b \cdot s^\beta$ with constants $\alpha\geq\beta$, $a$, and $b$.
    \item $r(\cdot)$ is convex and $f(\cdot)$ is concave.
\end{enumerate}
\end{corollary}
We remark that Corollary \ref{extremalConvexConcave} $(iii)$ provides theoretical foundation for numerical observations of \citet[p. 804]{uit2020condition}. Their experiments indicate that bang-bang policies -- or, on-off policies in their nomenclature -- have the same expected profit as the optimal policy under stochastic deterioration when the deterioration function is concave and the revenue function linear.

Under a bang-bang policy, the policy itself can be described by a switching curve, denoted with $\delta^{(\lambda,t)}$, such that the system, and hence production, is switched off if the deterioration level is greater than or equal to this curve. The following result establishes the monotonic behavior of this switching curve. This result is an immediate consequence of Theorem \ref{optPolicy}; hence we omit the proof.
\begin{corollary}\label{switchcurve}
For a bang-bang policy, the switching curve $\delta^{(\lambda,t)}$ is non-decreasing in the time to scheduled maintenance $t$, and non-increasing in the base rate $\lambda$. 
\end{corollary}
The optimality of a monotonic bang-bang policy is not only convenient because of its practical appeal but it also significantly decreases the computational burden of computing the optimal policy. Indeed, it effectively reduces the action space of our decision problem from a compact set with infinitely many actions to a discrete set with only two actions.  

\begin{figure}[htb!]
\begin{center}
\includegraphics[width=0.75\textwidth]{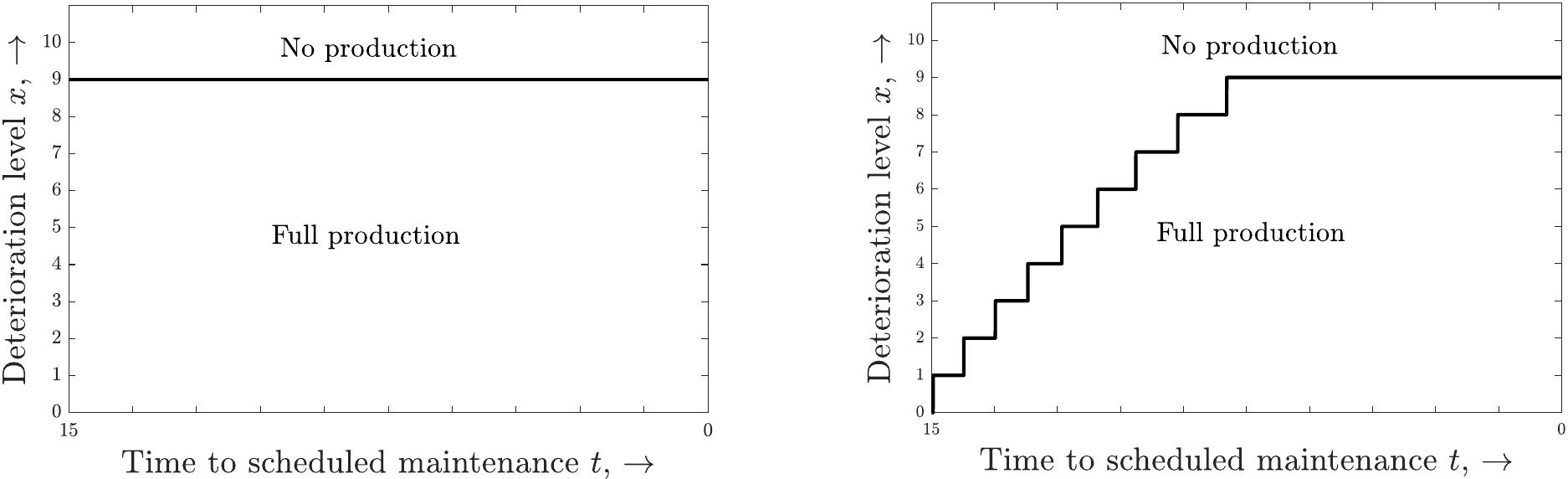}
\caption{Illustration of a bang-bang policy for base rate $\lambda = 1$ (left) and base rate $\lambda =4$ (right). The system is switched off if the deterioration level is greater than or equal to the switching curve (i.e. black line).}
\label{optimalBangIllus}
\end{center}
\end{figure}

Figure \ref{optimalBangIllus} illustrates the bang-bang policy. In this example we use the same values as in Figure \ref{optimalPolicyIllus}, except that we swap the revenue and deterioration function to comply with Corollary \ref{extremalConvexConcave}, i.e. $r(s)=s^{2}$ and $f(s)=s^{1/2}$. In line with Corollary \ref{switchcurve}, we see that the switching curve is non-decreasing in the time to scheduled maintenance, and non-increasing in the base rate.

\begin{remark}
    The optimal policy when multiple production systems seek to meet a joint demand rate is explored numerically in Appendix \ref{extensionMulti}. Here we observe that different production systems synchronize their production rates to meet the demand rate with machines in better health carrying a higher proportion of the total production rate.
\end{remark}

\subsection{Optimal Integration of Maintenance and Condition-Based Production}
\begin{figure}[!htb]
\begin{center}
\begin{subfigure}{.4\textwidth}
    \caption{}
\includegraphics[width=0.9\textwidth,keepaspectratio]{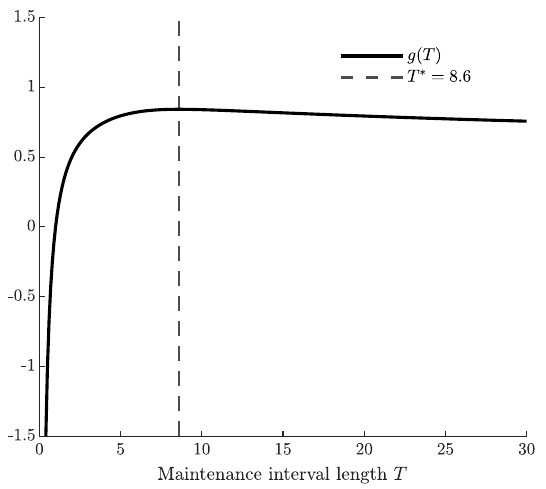}
  \label{fig:averageFigure}
\end{subfigure}
\begin{subfigure}{.4\textwidth}\centering
\caption{}
\includegraphics[width=0.77\textwidth,keepaspectratio]{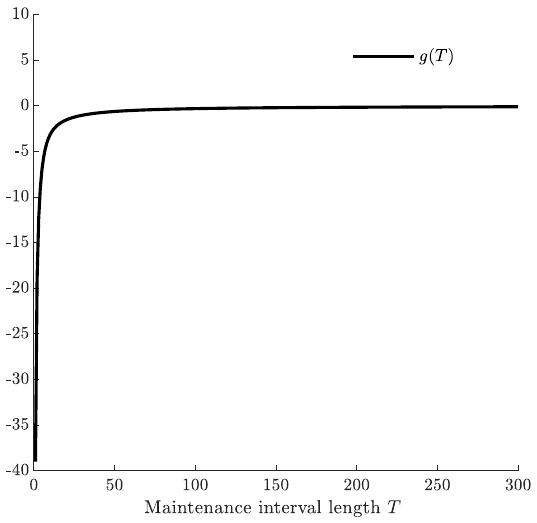}
  \label{fig:averageFigureCounter}
\end{subfigure}%
\end{center}
\caption{Two functions $g(T)$ when the base rate is 1. In Panel (a), the optimal maintenance interval is depicted with the dashed line. Panel (b) illustrates an example with no global maximizer.}
\label{average}
\end{figure}

Theorem \ref{optPolicy} characterizes the structure of the optimal production policy for the operational decision problem given that maintenance is performed every $T$ time units. In practice, a DM also needs to decide upon the value of $T$ in which she deals with a tactical trade-off. Specifically, frequent planned maintenance (small $T$) leads to unnecessary maintenance activities and associated costs, while infrequent maintenance (large $T$) both increases system failure risk and reduces revenue due to necessitated lower production rates. In deciding upon $T$, a DM seeks to maximize the average profit, by incorporating both the operational and tactical trade-off. More formally, let $g(T)$ denote the average profit per time unit when an optimal condition-based production policy $u_\lambda^*$ is followed on the operational level and the planned maintenance moment is scheduled every $T$ time units, that is, $
g(T) \triangleq \frac{1}{T}J_{\lambda}^*(0,T).$
The DM seeks to find $T^*\triangleq\arg\sup_{T\in\R_{>0}} g(T)$.
This decision problem is complex because of its nested structure. The next result simplifies its optimization considerably. 
\begin{theorem}\label{optimalT}
Let $T'$ be a strict local maximizer of $g(T)$, then $T'$ is the strict global maximizer. 
\end{theorem}
Condition-based production on an operational level and maintenance decisions on a tactical level should be made in conjunction rather than in isolation, and Theorem \ref{optimalT} provides us with a way to do so. Specifically, readily available uni-modal optimization methods can be used to efficiently find the optimal length of the maintenance interval.  Figure \ref{fig:averageFigure} illustrates such a uni-modal average profit function $g(T)$. In this example, where we use the same settings as in Figures \ref{optimalPolicyIllus} and \ref{samplepathOptimal}, the optimal maintenance interval length is equal to 8.6 with an average profit of 0.84 per time unit.  

In the next section, we present numerical studies for instances
as they are typically encountered in practice. For all these instances, $g(T)$ has a strict local maximizer and hence, by virtue of Theorem \ref{optimalT}, a unique maximizer. 
Unfortunately, this does not hold in general; we found that $g(T)$ may also be non-decreasing and negative for all $T$. Figure \ref{fig:averageFigureCounter} illustrates this case, where $g(T)$ is negative for all $T$ and with $g(T)$ approaching zero from below as $T$ grows large. In this example, $\lambda=1$ (base rate), $\xi=10$ (failure level), $c_m(\xi)=50$ (corrective maintenance cost), and $c_m(\cdot)=40$ for all non-failed states (preventive maintenance cost). The production rate space is $\mathscr{S}=[0,1]$ and we use $r(s)=s^{2}$ and $f(s)=s^{1/2}$ for the revenue and deterioration functions, respectively. 
Under this set of choices, the system deteriorates relatively faster than it generates revenue for all $s\in\mathscr{S}$, which, in combination with the high maintenance costs at the end of the planning horizon, results in a non-profitable system. In such a case, the potential revenue under any production policy is not sufficient to compensate the high maintenance costs, and choosing a large value of $T$ ensures that at least the maintenance costs are on average very low per time unit. Such non-profitable systems are arguably not the ones we operate in practice.

\section{Value of Condition-Based Production and Maintenance Management}
\label{sec:num}
This section contains two comprehensive numerical studies that highlight the practical value of our theoretical results on populations of systems with the same fixed base rate. The goal of this section is twofold. Firstly, we quantify the value of dynamically adjusting the production based on the system's deterioration and time to scheduled maintenance. Secondly, we investigate the value of integrating maintenance and condition-based production decisions. In both numerical studies, we use the canonical maintenance cost structure with a fixed preventive maintenance cost $c_m(\cdot) \triangleq c_p$ for all non-failed states, and a corrective maintenance cost $c_m(\xi)\triangleq c_u$ for the failed state.

\subsection{Static versus Condition-Based Production}
\label{subsec:staticCon}
We compare the optimal condition-based production policy with the optimal static policy that cannot control production in real time. 
Under this static policy, the DM chooses a fixed production rate, which is operated continuously in between consecutive maintenance moments until either failure or until the maintenance moment, whichever comes first. We call this static policy the fixed-rate policy, denoted by $\pi_{FS}$, and obtain it by solving the following maximization problem:
\begin{align}\label{fixedSetting}
\pi_{FS} \in \argmax_{s \in \mathscr{S}} \Big[ \E\big[ r(s) \cdot \min\{T_{\xi},T\} \big] - c_p - (c_u - c_p)\cdot\P \big[T_{\xi}\leq T \big]  \Big],
\end{align}
where $T_{\xi}$ is an Erlang random variable with shape $\xi$ and rate $\lambda f(s)$, and $T$ is the length of the maintenance interval. The first part of the objective function is the expected accumulated revenue, while the second part is the expected maintenance cost at time $T$. The maximization problem can be easily solved numerically to obtain $\pi_{FS}$. 

In order to assess the value of using the optimal condition-based production policy, $\pi_{CS}$, instead of the static fixed-rate policy, $\pi_{FS}$, we use the relative profit increase, denoted by $\mathcal{R}= 100\cdot (P_{\pi_{CS}}-P_{\pi_{FS}})/P_{\pi_{FS}}$, where $P_\pi$ is the expected profit (over the planning horizon) of policy $\pi\in\{\pi_{CS},\pi_{FS}\}$. 

\subsubsection{Full factorial experiment}
In this section and Appendix \ref{subsec:staticCBPapp}, we report on a full factorial experiment to quantify the value of condition-based production averaged across a wide range of parameter settings. Our test bed consists of 3750 instances, obtained by permuting various values for all input parameters. See Table \ref{tab:testbed1} in Appendix \ref{subsec:staticCBPapp} for all input parameter values of the test bed; here we describe only a subset. We use the functions $f=s^{\gamma}$ and $r=s^{\nu}$ to represent the deterioration and revenue functions, respectively. We use the interval $\mathscr{S}=[0,2]$ for the production rates to choose from. By using different values for $\gamma$ and $\nu$, our test bed consists of concave, linear, and convex deterioration and revenue functions. Furthermore, we use different base rates and failure levels so that we vary the expected time to failure under normal deterioration.

The main results of the full factorial experiment are summarized in
Table \ref{tab:resultsTestbedMain}, which presents the average and standard deviation of $\mathcal{R}$, as well as the average value $P_{\pi_{CS}}$ over all instances.  Detailed results are provided in Appendix \ref{subsec:staticCBPapp}. Table \ref{tab:resultsTestbedMain} indicates that adjusting the production rate based on the system's deterioration and time to maintenance, instead of using the best fixed-rate policy, leads on average to sizeable profit increases across a wide range of parameter settings. 
These profit increases are mainly driven by two mechanisms. First, the condition-based production policy is able to prevent failures by slowing down deterioration when the failure threshold is nearing, thereby saving on high corrective maintenance costs. This is different from the fixed-rate policy, which just continues to operate with the fixed production rate until failure. Second, when the scheduled maintenance moment is nearing and the system is still in good condition, the condition-based production policy will exploit the remaining useful condition to accumulate revenue by setting a higher production rate. Such opportunities are not seized by fixed-rate policies.

\begin{table}[htbp]
\centering
\caption{Main results static versus condition-based production.}
\fontsize{9pt}{9pt}\selectfont
\label{tab:resultsTestbedMain}
\begin{tabular}{ccc}
\toprule
 Mean {$\mathcal{R}$}  & SD {$\mathcal{R}$}  & $P_{\pi_{CS}}$ \\
\midrule
50.42 & 202.46 & 21.38 \\
    \bottomrule
    \end{tabular}%
\end{table}%
\vspace{-8pt} 

While on average condition-based production leads to significant profit increases, Table \ref{tab:resultsTestbed1} in Appendix \ref{subsec:staticCBPapp} also shows that the exact value depends on the specific input parameters. The value seems to depend strongly on the expected time to failure under normal deterioration (through $\lambda$ and $\xi$) and on how production affects both deterioration and revenue accumulation (through $\nu$ and $\gamma$). In the next section, we investigate the impact of these input parameters on the value of condition-based production through comparative static analyses.

\subsubsection{Comparative statics}\label{subsec:cbpparametric}
We define a base case of parameter settings, and then study how parametric changes with respect to this base case affect the relative profit increase of condition-based production $\mathcal{R}$. For the base case, we set the base rate $\lambda=1$, the failure level $\xi=14$, the preventive maintenance cost $c_p=2$, the corrective maintenance cost $c_u=10$, $\gamma=\nu=1$, and the planning horizon $T=10$. Note that these parameters are set such that they coincide with their respective median values of the full factorial experiment above, except for the planning horizon $T$. 

\begin{figure}[!htb]
\begin{subfigure}{.33\textwidth}
  \centering
    \caption{Relative profit increase $\mathcal{R}$.}
\includegraphics[width=0.97\textwidth,keepaspectratio]{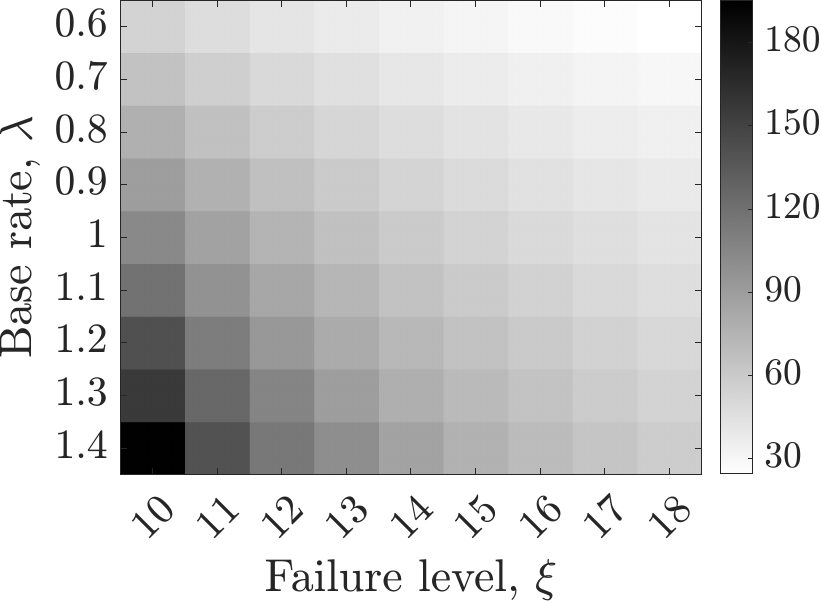}
  \label{fig:lf_profitIncrease}
\end{subfigure}
\begin{subfigure}{.33\textwidth}\centering
\caption{Revenue under $\pi^*_{CS}$.}
\includegraphics[width=0.94\textwidth,keepaspectratio]{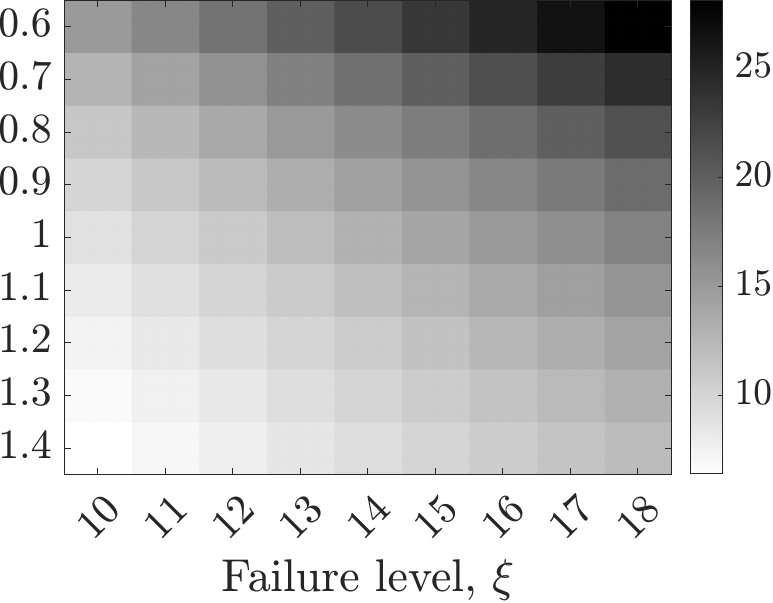}
  \label{fig:lf_revCS}
\end{subfigure}%
\begin{subfigure}{.33\textwidth}
  \centering
    \caption{Revenue under $\pi^*_{FS}$.}
\includegraphics[width=0.94\textwidth,keepaspectratio]{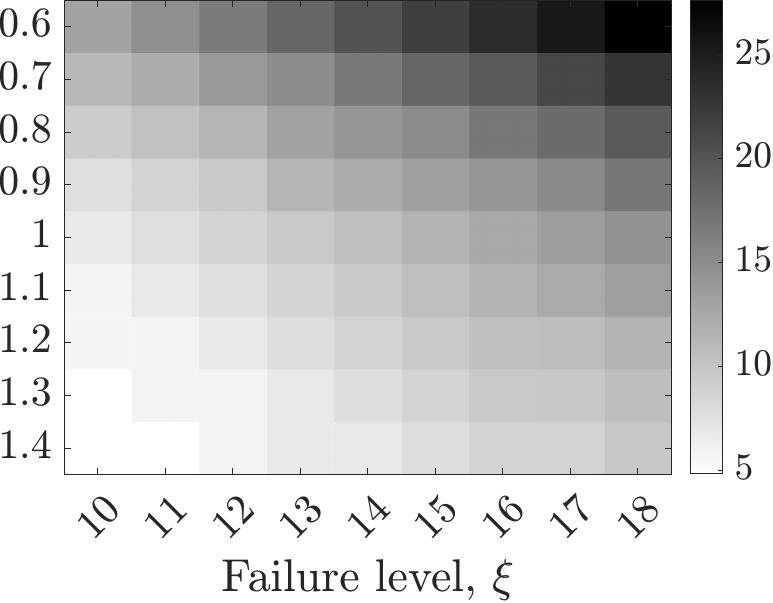}
  \label{fig:lf_revFS}
\end{subfigure}
\par\medskip
\begin{subfigure}{.34\textwidth}\centering
 \caption{Fixed rate under $\pi^*_{FS}$.}
\includegraphics[width=0.97\textwidth,keepaspectratio]{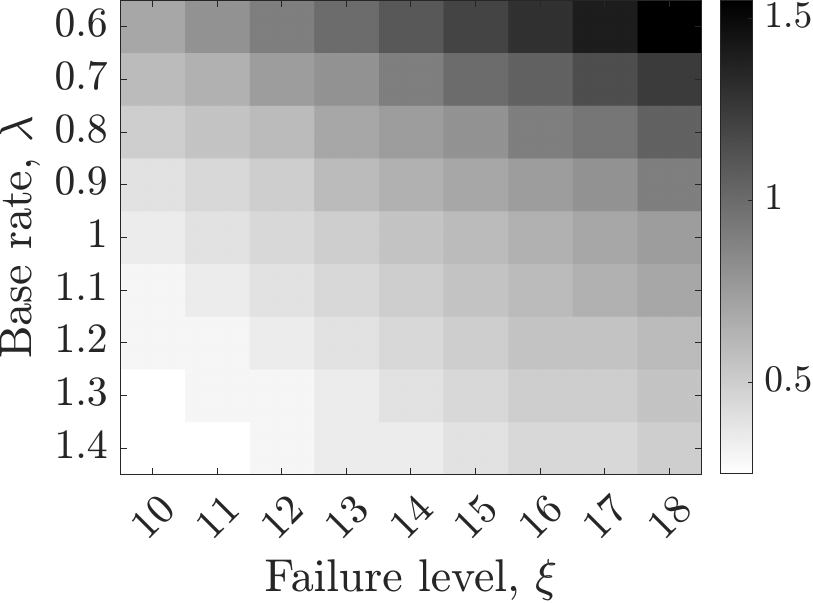}
  \label{fig:lf_rate}
\end{subfigure}%
\begin{subfigure}{.33\textwidth}
  \centering
    \caption{$P\{X_0\geq \xi\}$ under $\pi^*_{FS}$.}
\includegraphics[width=0.94\textwidth,keepaspectratio]{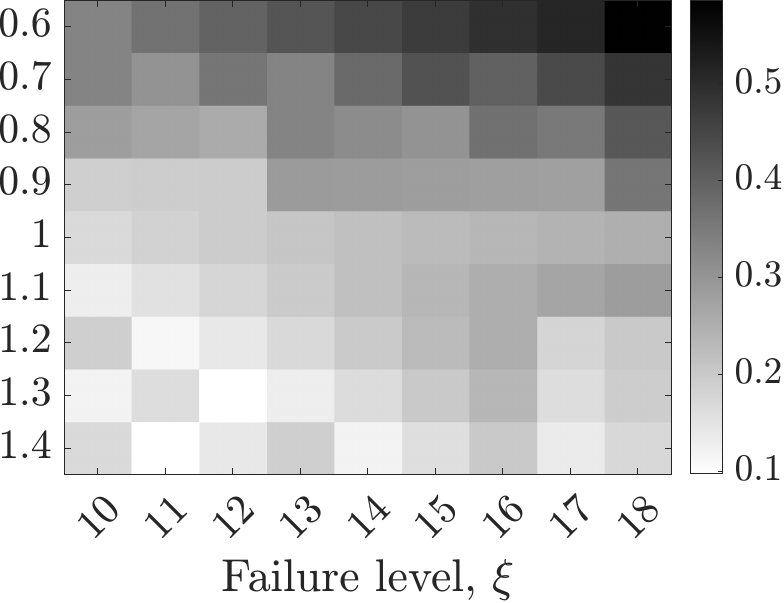}
  \label{fig:lf_probCorrective}
\end{subfigure}
\begin{subfigure}{.32\textwidth}
\centering
 \caption{Cost under $\pi^*_{FS}$.}
\includegraphics[width=0.94\textwidth,keepaspectratio]{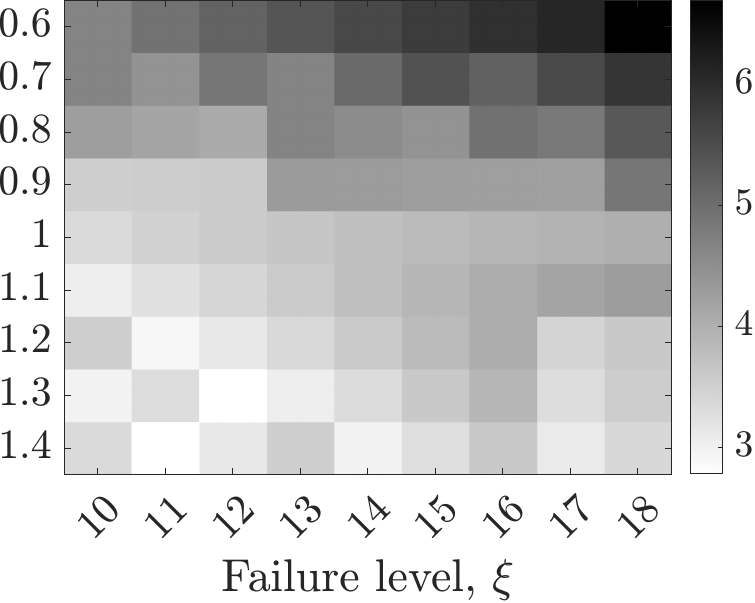}
  \label{fig:lf_costFS}
\end{subfigure}%
\caption{Various output measures for different values of the base rate $\lambda$ and the failure level $\xi$ when the other parameters are held fixed at the base case values. Note that the expected time to failure under normal deterioration increases as the base rate (failure level) decreases (increases).}
\label{lf_parametric}
\end{figure}

We first investigate how the relative value of condition based production is affected by the expected time to failure under normal deterioration. Figure \ref{fig:lf_profitIncrease} illustrates that as the expected time to failure increases (i.e. starting from $(\lambda,\xi)=(1.4,10)$ moving along either axis), the value of optimal condition based production polices over static fixed-rate policies decreases substantially from more than 180\% to less than 30\%. Comparing Figures \ref{fig:lf_revCS} and \ref{fig:lf_revFS}, we see that the optimal condition-based production policy consistently accumulates more revenue than the optimal static fixed-rate production policy. Observe further that this difference diminishes as the expected time to failure increases, while the absolute revenues of both policies increase owing to the increased time available for production. Figure \ref{fig:lf_probCorrective} indicates that the optimal rate of the fixed-rate production policy increases in the expected time to failure. This enables more revenue accumulation, but comes at the expense of a higher probability of corrective maintenance at the end of the planning horizon, as indicated in Figure \ref{fig:lf_probCorrective}, and thus a higher total maintenance cost, as indicated in Figure \ref{fig:lf_costFS}. By contrast, under the optimal condition-based production policy, the deterioration level never reaches the failure threshold for any of the instances considered in this analysis, resulting in only incurring the preventive maintenance cost at the end of the planning horizon. 

Figure \ref{lf_parametric} suggest -- and additional numerical experiments confirm this -- that as the expected time to failure increases, the optimal rate of the fixed-rate policy increases towards the maximum production rate. For sufficiently large expected times to failure, the probability of corrective maintenance goes to zero under the maximum production rate, thus implying that the revenue and cost of both the fixed-rate and the condition-based production policy converge to the same values. That is to say, the relative profit increase of condition-based production over fixed-rate production goes to zero as the expected time to failure grows large. This argument can be made more rigorous by noting that the probability of corrective maintenance $\mathbb{P}[T_\xi\leq T]$ in the expected profit function \eqref{fixedSetting} under the fixed rate $\bar{s}$ ($<\infty$) goes to zero as the expected time to failure $\E[T_\xi]$ grows large, through increasing the failure level $\xi$ and/or decreasing the base rate $\lambda$. 

\begin{figure}[tbh!]
\begin{subfigure}{.33\textwidth}
  \centering
    \caption{Relative profit increase $\mathcal{R}$.}
\includegraphics[width=0.97\textwidth]{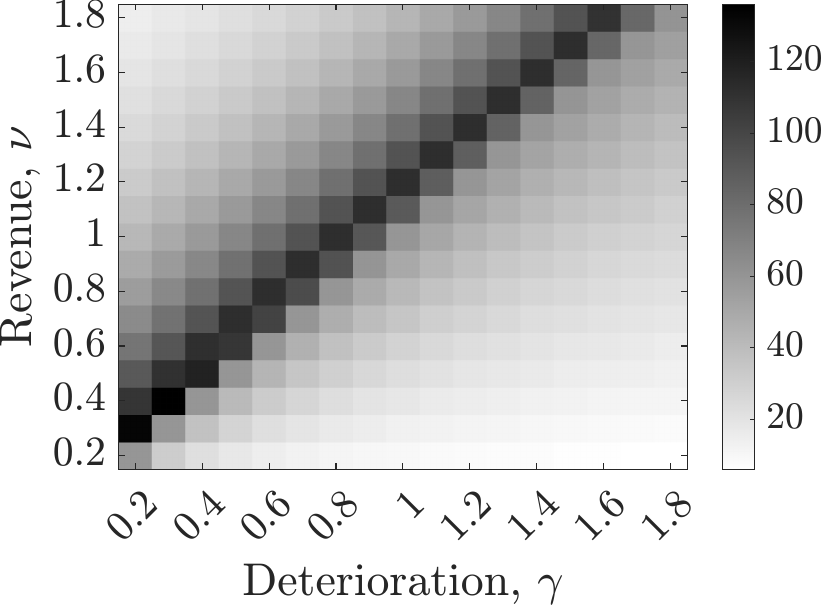}
  \label{fig:pd_profitIncrease}
\end{subfigure}
\begin{subfigure}{.33\textwidth}\centering
\caption{Revenue under $\pi^*_{CS}$.}
\includegraphics[width=0.94\textwidth]{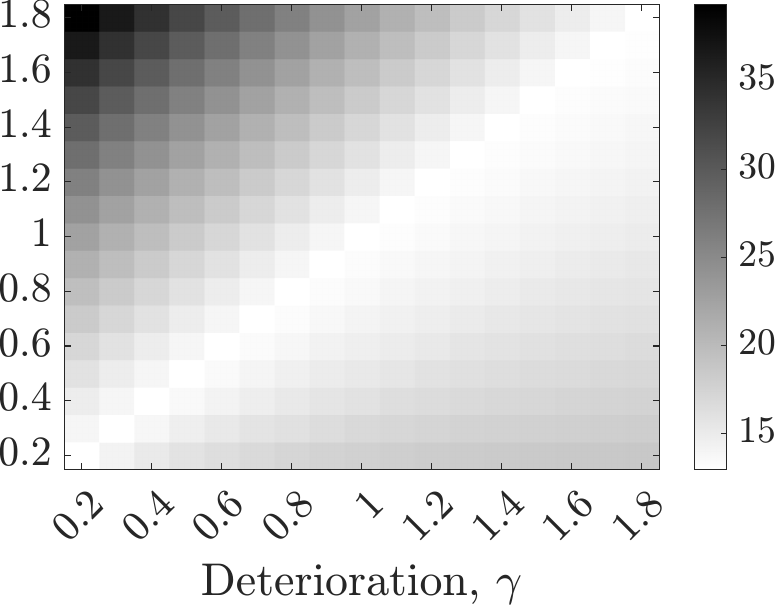}
  \label{fig:pd_revCS}
\end{subfigure}%
\begin{subfigure}{.33\textwidth}
  \centering
    \caption{Revenue under $\pi^*_{FS}$.}
\includegraphics[width=0.94\textwidth]{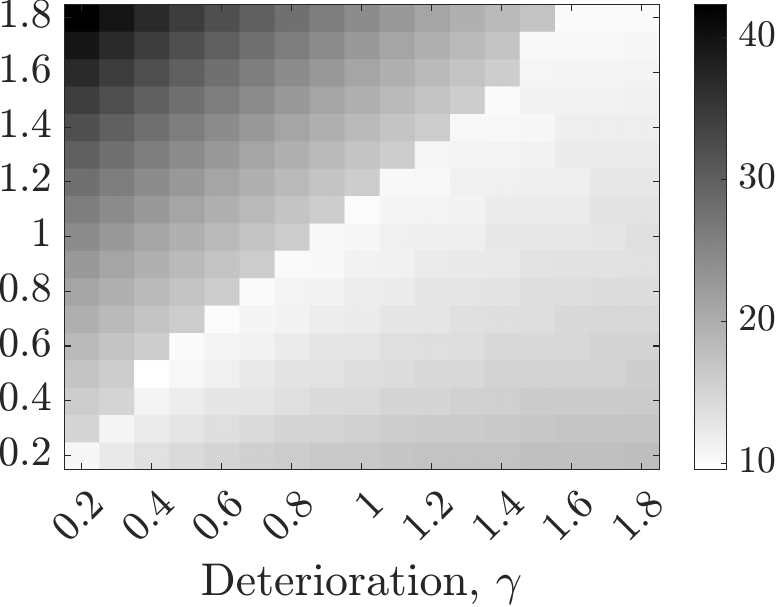}
  \label{fig:pd_revFS}
\end{subfigure}
\par\medskip
\begin{subfigure}{.33\textwidth}\centering
 \caption{Fixed rate under $\pi^*_{FS}$.}
\includegraphics[width=0.97\textwidth]{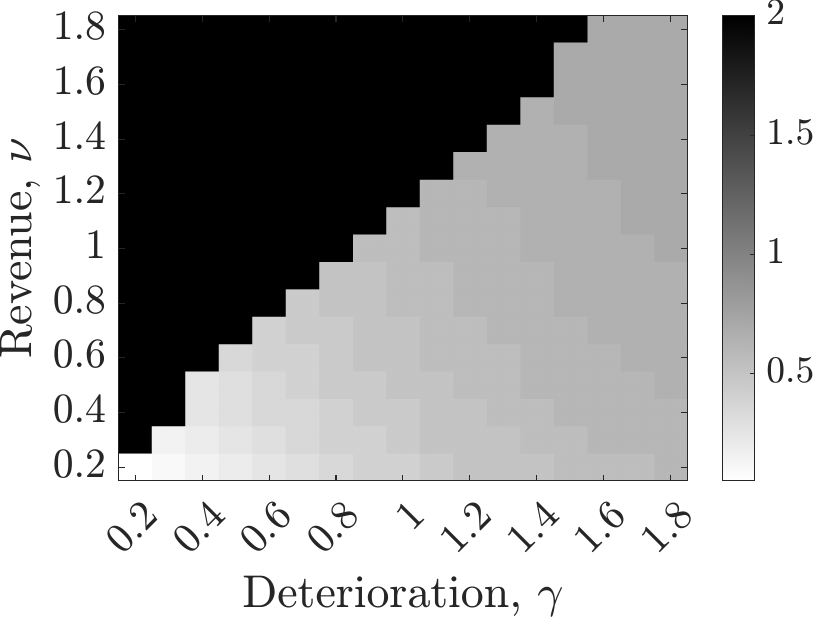}
  \label{fig:pd_rate}
\end{subfigure}%
\begin{subfigure}{.33\textwidth}
\centering
 \caption{Cost under $\pi^*_{CS}$.}
\includegraphics[width=0.94\textwidth]{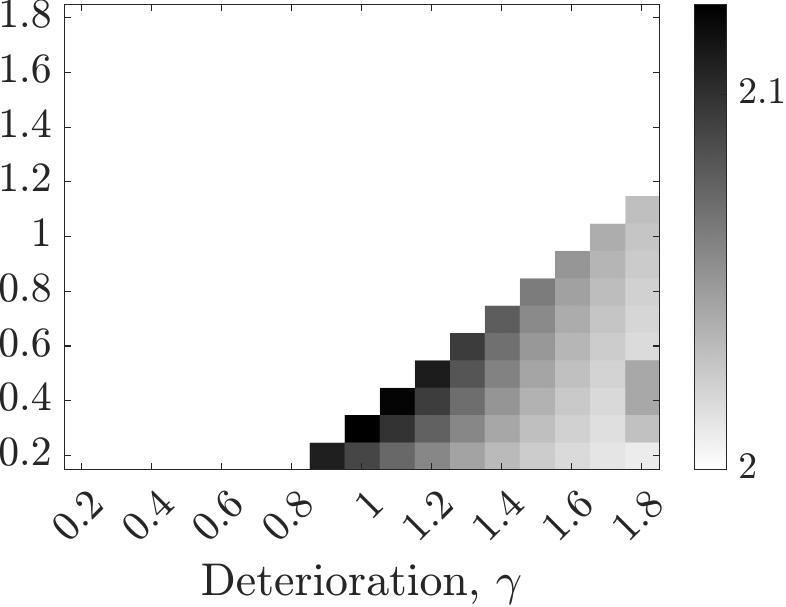}
  \label{fig:pd_costCS}
\end{subfigure}%
\begin{subfigure}{.33\textwidth}
  \centering
    \caption{Cost under $\pi^*_{FS}$.}
\includegraphics[width=0.94\textwidth]{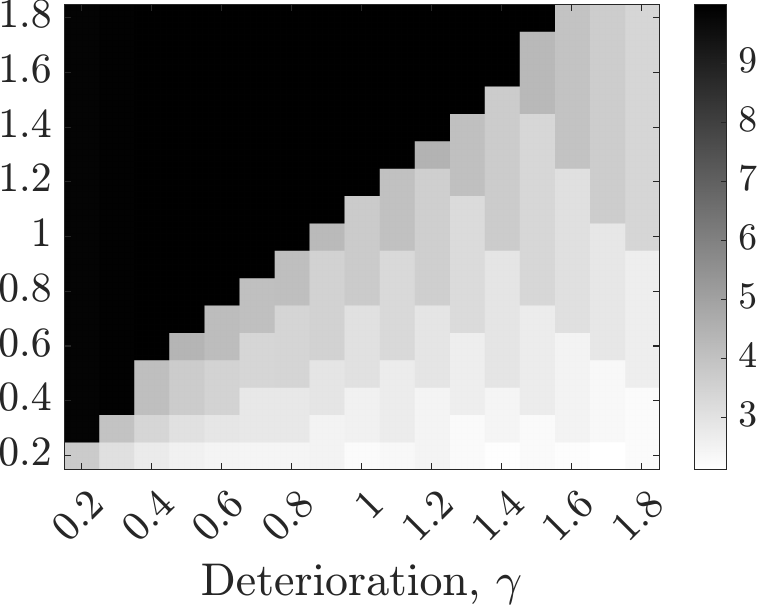}
  \label{fig:pd_costFS}
\end{subfigure}
\caption{Various output measures for different values of the revenue function parameter $\nu$ and the deterioration function parameter $\gamma$ when the other parameters are held fixed at the base case values. Note that the optimal condition-based production policy for the instances on the anti-diagonal and above (i.e. $\nu\geq \gamma$) is of the bang-bang type (cf. Corollary \ref{extremalConvexConcave}).}
\label{pd_parametric}
\end{figure}

Next we investigate the effect of the deterioration and revenue functions on the value of condition-based production. Figure \ref{fig:pd_profitIncrease} indicates that the relative profit increase of optimal condition-based production over static fixed-rate production increases substantially from less than 20\% to more than 120\% as we move closer to the anti-diagonal consisting of instances for which the revenue function and deterioration function parameter are identical. This suggests that condition-based production is most beneficial in settings where production impact on revenue and deterioration are similar. Figure \ref{fig:pd_rate} indicates that the optimal fixed rate equals $\bar{s}$ for most of the instances left of the anti-diagonal where the optimal condition-based production policy is of the bang-bang type (i.e. $\nu\geq \gamma$). This leads to slightly higher revenue accumulation than what the optimal condition-based production policy accumulates. However, different from the condition-based production policy, the fixed-rate policy cannot be turned off when failure is nearing. Consequently, the fixed-rate policy almost always leads to failure in these instances, as indicated by the costs of these instances in Figure \ref{fig:pd_costFS} being roughly equal to the corrective maintenance cost.
For all instances on the right of the anti-diagonal (i.e. $\gamma > \nu$) we observed numerically that the optimal condition-based production policy is not of the bang-bang type. For these instances, we deduce that condition-based production consistently generates higher revenues than fixed-rate production while simultaneously incurring lower maintenance costs, and that absolute revenues and profits of these instances under both production policies are typically lower than the instances for which $\nu\geq \gamma$. 

Finally, Table \ref{tab:resultsTestbed1} in Appendix \ref{subsec:staticCBPapp} indicates that the relative profit increase of condition-based production averaged across a wide range of settings as well as the variability in this profit increase across those settings increase (decrease) in $\nu$ $(\gamma)$. Our discussions above provides explanation for this behavior: As $\nu$ ($\gamma$) increases (decreases), the number of instances for which the optimal condition-based production policy has the bang-bang structure increases (decreases) as well (i.e. $\nu\geq \gamma$), and our observations above become more persistent.

\subsection{Sequential versus Integrated Maintenance Scheduling}
\label{subsec:seqInt}
Theorem \ref{optimalT} enables efficient optimization of the maintenance interval length when the optimal condition-based production policy is followed in between those maintenance moments. This allows DMs to jointly decide upon the two policies (i.e. maintenance-interval and production policies). In practice, these two decisions are often made sequentially: First, the tactical decision on the maintenance interval length is made, followed by the operational decisions on production control based on that interval length.
This is typically approached as follows. The maintenance interval length problem -- without taking into account the production policy and using only the knowledge of the base rate -- is a classical age-based maintenance problem \citep{barlow1960optimum}. The lifetime under the base rate is an Erlang random variable with shape $\xi$ and rate $\lambda$. Hence, the average maintenance cost rate for a given maintenance interval length $t$, denoted $g_{mc}(t)$, equals \begin{align}\label{optmaintlength}
   g_{mc}(t)= \frac{c_p + (c_u-c_p)\cdot \P\big[T_{\xi} \leq  t\big] }{\E\big[\min\{T_{\xi},t\} \big]},
\end{align}
where the numerator is the expected maintenance cost, and the denominator is the expected time to maintenance (either preventive or corrective). Since the hazard rate of an Erlang random variable is increasing, we know that the function $g_{mc}(t)$ has a unique minimizer that is easily computed. Given this optimized maintenance interval length, denoted $t^*_S$, we can then use the optimal condition-based production policy to control the system's deterioration in between those maintenance moments, which yields an expected profit rate per time unit. We call this approach the sequential approach, denoted with $\pi_{S}$. The resulting average profit rate is defined as $\hat{P}_{\pi_{S}} \triangleq J_{\lambda}^{*}(0,t^*_S)/t^*_S$.

Alternatively we can use Theorem \ref{optimalT} and determine the optimal length of the maintenance interval when taking into account that the optimal condition-based condition-based production policy is followed in between those maintenance moments. We call this approach the integrated approach, denoted with $\pi_{I}$. The resulting average profit rate is given by $\hat{P}_{\pi_{I}} \triangleq J_{\lambda}^{*}(0,t^*_I)/t^*_I$, where $t^*_I$ is the optimal length of the maintenance interval, readily obtained by virtue of  Theorem \ref{optimalT}.
The relative value of integrating maintenance and production policies is expressed by
$\mathcal{\hat{R}} = 100\cdot (\hat{P}_{\pi_{I}}-\hat{P}_{\pi_{S}})/\hat{P}_{\pi_{S}}$.

\subsubsection{Full factorial experiment}
As before, we first report on a full factorial experiment to quantify the value of integrating maintenance and condition-based production averaged across a wide range of settings. We use the same parameters as the previous full factorial experiment except the parameter $T$, resulting in 1875 instances. The main results of this full factorial experiment are summarized in Table \ref{tab:resultsTestbedMain2}, which shows the average and standard deviation of $\mathcal{\hat{R}}$, as well as the average value $\hat{P}_{\pi_{I}}$ over all instances. Detailed results can be found in Appendix \ref{subsec:seqIntApp}.

\begin{table}[htbp]
\centering
\caption{Main results sequential versus integrated maintenance scheduling.}
\fontsize{9pt}{9pt}\selectfont
\label{tab:resultsTestbedMain2}
\begin{tabular}{ccc}
\toprule
Mean {$\mathcal{\hat{R}}$}  & SD {$\mathcal{\hat{R}}$} & $\hat{P}_{\pi_{I}}$ \\
\midrule
21.39 & 28.35 & 1.83 \\
    \bottomrule
    \end{tabular}%
\end{table}%
\vspace{-8pt} 

Table \ref{tab:resultsTestbedMain2} shows that integrating the length of maintenance intervals with the optimal condition-based production policy for revenue accumulation in between those maintenance moments leads to significant profit rate increases as opposed to treating them sequentially, with an average increase of 21.39\%. Given that both approaches apply the optimal condition-based production policy on the operational level, we emphasize that this average profit increase can be fully attributed to the integration of maintenance and condition-based production. As with the value of condition-based production explored in the previous section, Table \ref{tab:resultsTestbed2} in Appendix \ref{subsec:seqIntApp} indicates that the value of the integrated approach is significantly influenced by the expected time to failure and how production affects both deterioration and revenue accumulation. We will now investigate this further.

\subsubsection{Comparative statics}
We follow the same approach as in the comparative static analyses of the previous section. The base case remains the same, with the input parameters set such that they coincide with their respective median values used in the full factorial experiment (except $T$ for obvious reasons). We note that for all instances in this section, we observed numerically that the optimal condition-based production policy on the operational level ensures that the system never fails under both the sequential and the integrated approach. Consequently, only the preventive maintenance cost is incurred upon the scheduled maintenance moment for all instances.

\begin{figure}[htb!]
\begin{center}
\begin{subfigure}{.33\textwidth}
    \caption{Relative profit increase $\mathcal{\hat{R}}$.}
\includegraphics[width=0.96\textwidth]{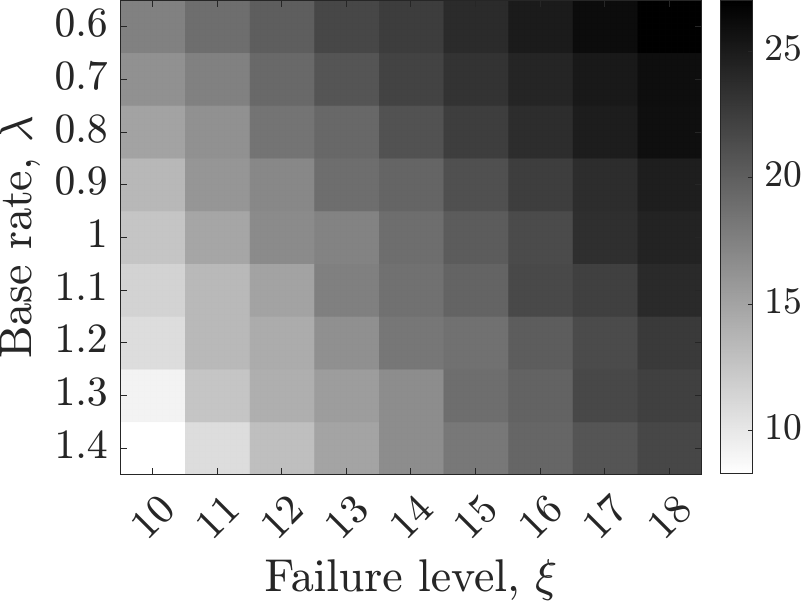}
  \label{fig:lf_T_profitIncrease}
\end{subfigure}
\begin{subfigure}{.33\textwidth}\centering
\caption{Ratio $t^*_{I}/t^*_{S}$.}
\includegraphics[width=0.96\textwidth]{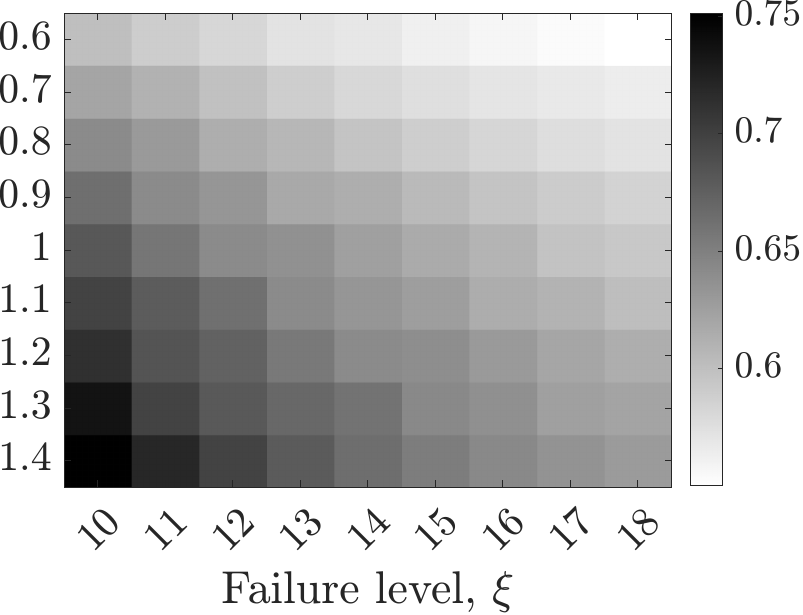}
  \label{fig:lf_T_optT}
\end{subfigure}%
\caption{Output measures for different values of the base rate $\lambda$ and the failure level $\xi$ when the other parameters are held fixed at the base case values.}
\label{lf_T_parametric}
\end{center}
\end{figure}

Figure \ref{fig:lf_T_profitIncrease} indicates that the value of integrating maintenance and condition-based production increases from less then 10\% to more than 25\% as the expected time to failure increases. While both approaches accumulate more profit per time unit as the expected time to failure increases, Figure \ref{fig:lf_T_profitIncrease} indicates that the integrated approach does so at an increasingly higher rate. This latter observation can be explained by means of Figure \ref{fig:lf_T_optT}. As the expected time to failure increases, both approaches increase the maintenance interval length, enabling more revenue accumulation on the operational level. However, as depicted in the figure, the maintenance interval length under the integrated approach grows slower than under the sequential approach. This highlights the advantages of integrating maintenance and condition-based production decisions as opposed to treating them sequentially. The integrated approach sets the maintenance interval length such that revenue accumulation in between those maintenance moments can be maximized while taking into account the maintenance cost at the end of the maintenance interval. By contrast, the sequential approach neglects this revenue accumulation when determining the maintenance interval, leading to longer intervals in which less revenue is accumulated to avoid the cost of failure replacement. The integrated approach's shorter maintenance intervals lead to more frequent maintenance moments and associated preventive maintenance costs. However, these higher costs are significantly offset by the increased revenue accumulation, especially for longer expected times to failure, all while keeping the risk of corrective maintenance small. 
Thus the integrated approach strikes an optimal balance between increased maintenance costs and greater revenue generation. 

\begin{figure}[htb!]
\begin{center}
\begin{subfigure}{.33\textwidth}
    \caption{Relative profit increase $\mathcal{\hat{R}}$.}
\includegraphics[width=0.97\textwidth]{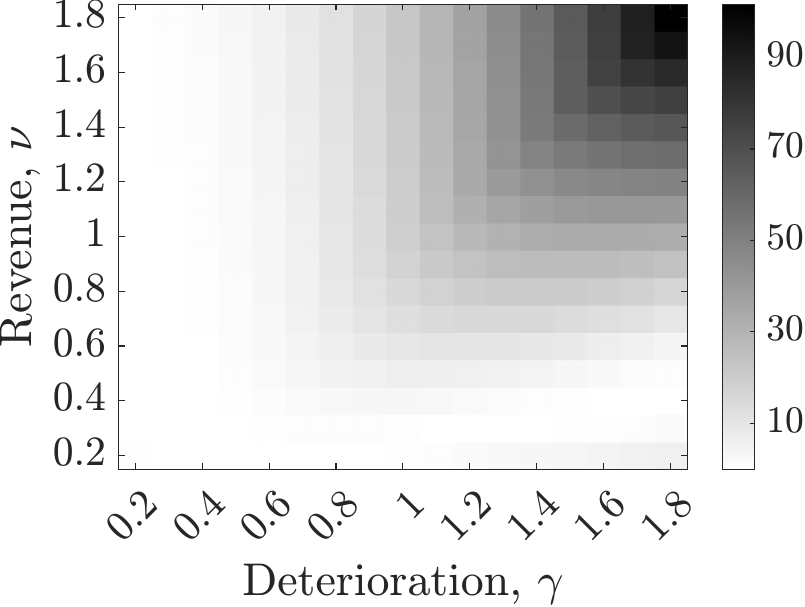}
  \label{fig:pd_T_profitIncrease}
\end{subfigure}
\begin{subfigure}{.33\textwidth}\centering
\caption{Ratio $t^*_{I}/t^*_{S}$.}
\includegraphics[width=0.94\textwidth]{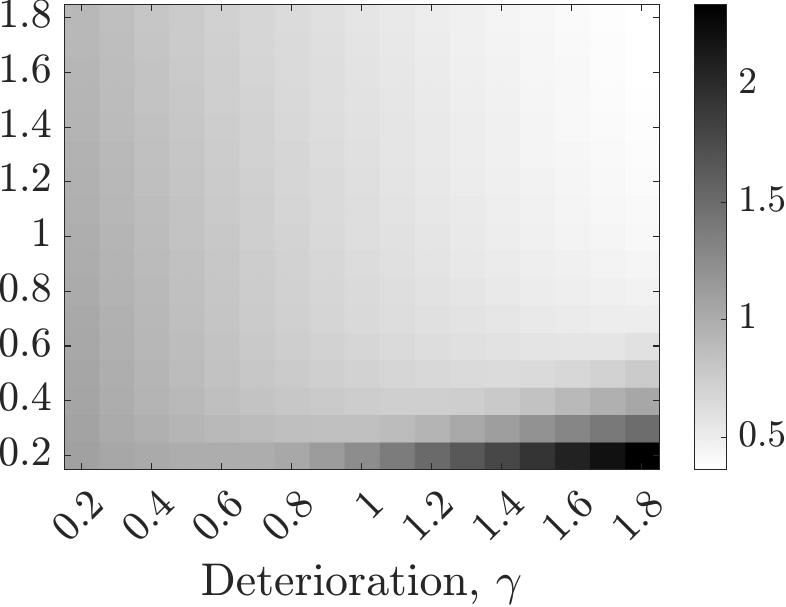}
  \label{fig:pd_T_optT}
\end{subfigure}%
\par \medskip
\begin{subfigure}{.33\textwidth}
    \caption{Revenue rate under $\pi_I$.}
\includegraphics[width=0.97\textwidth]{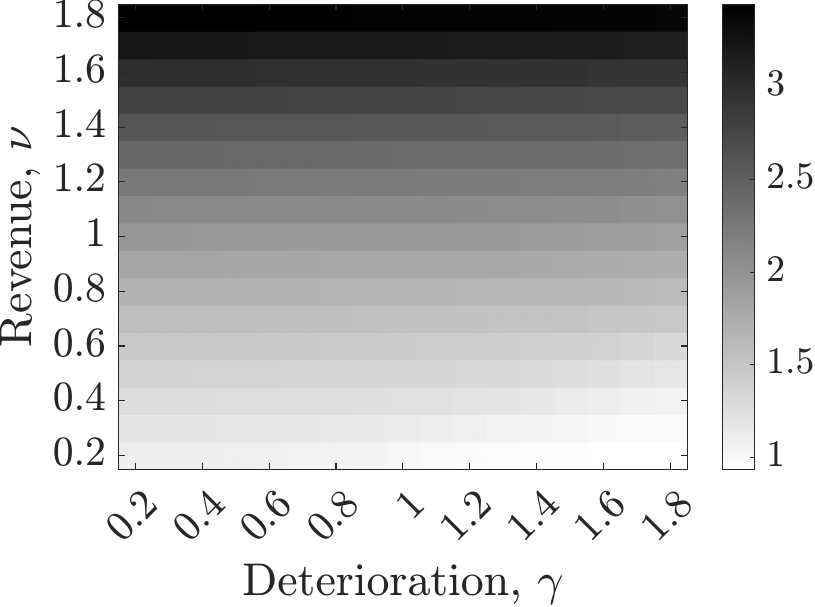}
  \label{fig:pd_T_revInt}
\end{subfigure}
\begin{subfigure}{.33\textwidth}\centering
\caption{Revenue rate under $\pi_S$.}
\includegraphics[width=0.94\textwidth]{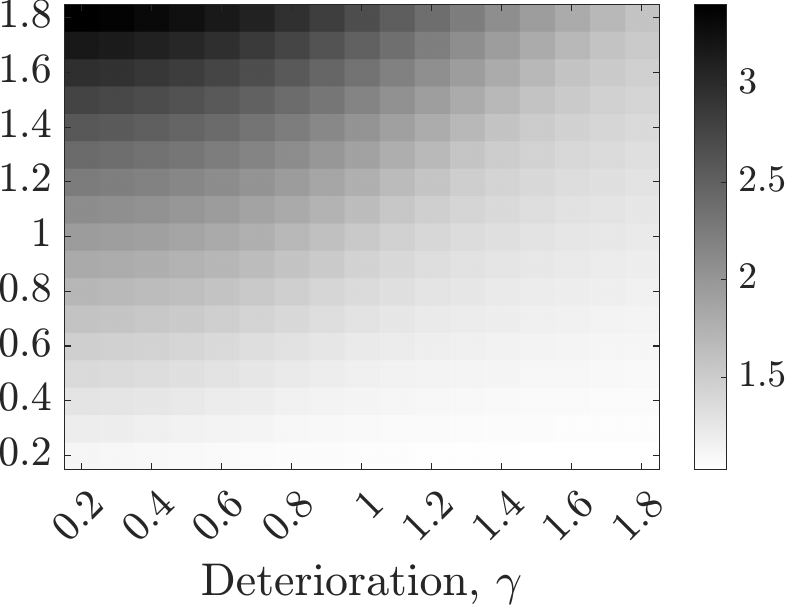}
  \label{fig:pd_T_revSeq}
\end{subfigure}%
\caption{Various output measures for different values of the revenue function parameter $\nu$ and the deterioration function parameter $\gamma$ when the other parameters are held fixed at the base case values.}
\label{pd_T_parametric}
\end{center}
\end{figure}

Next we investigate the impact of the deterioration and revenue functions on the value of integrating maintenance and condition-based production decisions. Figure \ref{fig:pd_T_profitIncrease} indicates that this value increases as the impact of production on revenue accumulation and/or deterioration becomes increasingly higher through increases in $\nu$ and/or $\gamma$. This is expected: The sequential approach neglects operational production decisions in between maintenance moments, and in particular their effects on revenue accumulation and deterioration. As these effects increase, the value of setting the maintenance interval length in conjunction with the operational production decisions increases. Figure \ref{fig:pd_T_optT} depicts the ratio between the optimal maintenance interval length under the integrated approach and the maintenance interval length under the sequential approach, the latter being constant for all instances as it assumes deterioration accumulation under normal operation (cf. Equation \eqref{optmaintlength}).
For the instances left on the anti-diagonal, where the optimal condition-based production policy is of the bang-bang structure, we see that the optimal maintenance interval length decreases as the impact of production deterioration increases. By contrast, the optimal maintenance interval length is robust against changes in the impact of production on revenue accumulation since this impact does not affect the time to failure and only leads to more revenue accumulation (see Figure \ref{fig:pd_T_revInt}). Producing at the maximum rate leads to faster deterioration under higher values of the deterioration function parameter $\gamma$. Hence if the optimal policy has a bang-bang structure, shortening the maintenance interval length allows for more revenue accumulation per time unit, outweighing the more frequently incurred preventive maintenance costs.    
We observe a different structure for the instances on the right of the anti-diagonal.
If all production rates lead to more deterioration, lengthening the maintenance interval length allows for more conservative condition-based production rates such that revenues can be maximized without the risk of failure. 

Figures \ref{fig:pd_T_revInt} and \ref{fig:pd_T_revSeq} show that revenue rates increase as the revenue function parameter increases. By integrating maintenance length intervals and production decisions, Figure \ref{fig:pd_T_revInt} clearly illustrates that this revenue rate is unaffected by increases in the deterioration function parameter. By contrast, Figure \ref{fig:pd_T_revSeq} shows that the revenue rate decreases in the deterioration function parameter under the sequential approach. This is due to the fact that the maintenance interval length is kept constant under the sequential approach, necessitating a condition-based production policy that aims to avoid failure, rather than to accumulate revenue.

\section{Unit-To-Unit Variability: Learning the Base Rate} 
\label{sec:heterogenous}
The previous sections focused on the situation in which the base rate is known and fixed for each system. In this section, we consider the situation in which the base rate is not known and where each newly installed system has a different base rate. In dealing with such unit-to-unit variability in base rates, a DM needs to learn the individual base rate based on the system's deterioration path and the used production policy. Such populations of systems subject to shock deterioration with variability in the base rates are quite prevalent in practice \citep[see, e.g.,][]{drent2021cbm}.

We first propose a Bayesian procedure to learn the unknown base rate and use it to develop a heuristic policy that is inspired by the case where the base rate is known. We then discuss the results of a simulation study in which we compare the performance of the heuristic with an Oracle policy that knows the true base rate a-priori. 

\subsection{Bayesian Learning under Production Policies}
\label{sec:unknownbase}                                                                                 
In line with previous research on unit-to-unit variability in shock deterioration processes \citep[e.g.,][]{drent2021cbm,drent2023pool}, we model unit-to-unit variability in the systems' base rates with a Gamma distributed random variable, denoted by $\Lambda$,  with known shape $\alpha_0$ and rate $\beta_0$; that is, $\Lambda \sim \Gamma (\alpha_0, \beta_0)$. 
Each newly installed system has an a-priori unknown realization of $\Lambda$ as base rate. Hence at $t=0$ ($T$ time units before planned maintenance), upon installation of a new system, the DM's knowledge regarding the unknown base rate is modeled by the prior density function $\varphi_{0}(\lambda|\alpha_0,\beta_0) = \frac{ \beta_0^{\alpha_0} \lambda^{\alpha_0-1} e^{-{\beta_0} \lambda}}{\Gamma(\alpha_0)}$, where the subscript $0$ indicates that this represents $t=0$. Suppose that the DM adopts a production policy $\bm{s}_t\triangleq \{s_u\}_{u\in [0,t)}$ with production rate $s_u\in\mathscr{S}$ at time $u$ over the period $[0, t)$ and the deterioration level at time $t$ is $Y_{t} = y$, the proposition below shows that the posterior distribution at time $t$ is then again a Gamma distribution, but with updated parameters.
\begin{proposition} \label{inferenceUnderPolicy}
The updated posterior distribution of $\Lambda$ at time $t$ given realized deterioration level $Y_t=y$ and an employed production policy $\bm{s}_t$ up to time $t$, is a Gamma distribution  with parameters $\alpha_0+y$ and $\beta_0+ \int_{u=0}^t f(s_u)du$; that is, $\Lambda \sim \Gamma\left(\alpha_0+y,\beta_0+\int_{u=0}^t f(s_u)du\right)$.
\end{proposition}
Proposition \ref{inferenceUnderPolicy} implies that inference of the unknown base rate at each time $t$ only requires the employed production policy $\bm{s}_t$ and the deterioration level $Y_t=y$. 
The mean of this updated Gamma distribution is equal to $\frac{\alpha_0 + y}{\beta_0+ \int_{u=0}^t f(s_u)du}$, which has intuitive appeal. If a system has been controlled by a policy with a higher aggregate amount of production rates captured by $\int_{u=0}^t f(s_u)du$ (recall that $f(\cdot)$ is increasing) than another, but both have the same deterioration, then we expect the base rate of the former to be smaller. Similarly, if a policy has been operated with the same aggregate amount of production rates, but one system has a higher deterioration level than another, then we expect the system with more deterioration to have a larger base rate.  

\subsection{A Certainty-Equivalent Policy}
We now propose an intuitive heuristic policy called the certainty-equivalent (CE) policy to deal with unit-to-unit variability in the base rate that makes use of the Bayesian framework described in the previous section. Under this policy, the DM estimates the base rate based on her current knowledge and assumes that the true base rate is equivalent to her estimation with certainty. The CE policy is parameterized by a positive integer  $N^{\text{opt}}\in \N$, which denotes how often the DM re-optimizes the policy based on the accrued knowledge. In the CE policy, we divide the planning horizon in $N^{\text{opt}}+1$ phases, each of equal length.

\begin{figure}[!htb]
        \centering
        \resizebox{0.6\textwidth}{!}{
\begin{tikzpicture}
\draw[thick, |-|] (0,0) -- (10,0);

\foreach \x in {1,2,3,4,5,6,7,8,9}
\draw (\x cm,1.5pt) -- (\x cm,-1.5pt);


\draw [black, thick ,decorate,decoration={brace,amplitude=5pt}] (0.1,0.2)  -- (4.9,0.25) 
       node [black,midway,above=3pt,xshift=-2pt] {\footnotesize Use policy $u^*_{\hat{\lambda}_0}$ };
       
\draw [black, thick ,decorate,decoration={brace,amplitude=5pt}] (5.1,0.2)  -- (9.9,0.25) 
       node [black,midway,above=3pt,xshift=-2pt] {\footnotesize  Use policy $u^*_{\hat{\lambda}_1}$};
       
\draw[] (0,0) node[below=2pt] {$t=0$} node[above=2pt] {};
\draw[] (5,0) node[below=2pt] {$t_1=\frac{T}{2}$} node[above=2pt] {};
\draw[] (10,0) node[below=2pt] {$t=T$} node[above=2pt] {};

\draw[ latex-, line width=1pt] (5,0.15) -- (5,1);
\draw[] (5,0.7) node[above=3pt,thick] {\footnotesize Compute $\hat{\lambda}_1=\frac{\alpha_0+x_1}{\beta_0+\bm{s}_{t_1}}$, and obtain policy $u^*_{\hat{\lambda}_1}$} node[above=3pt] {};
\end{tikzpicture}
}\caption{Schematic illustration of CE policy with  $N^{\text{opt}}=1$.}
\label{illucep}
\end{figure}
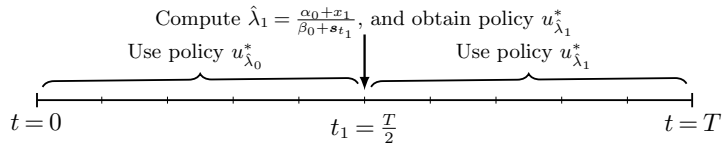

At the start of the planning horizon, upon installation of the system, the DM's estimate for the base rate, denoted with $\hat{\lambda}_0$, is equal to the mean of the initial prior distribution: $\hat{\lambda}_0 = \frac{\alpha_0}{\beta_0}$. Using this estimate, we obtain the optimal policy for the known case with $\hat{\lambda}_0$ as base rate, and implement this policy until the first re-optimization epoch at $t_1=\frac{T}{N^{\text{opt}}+1}$. At this point, the DM updates the estimate of the base rate based on the current deterioration level, denoted with $x_1$, and the policy used $\bm{s}_{t_1}$, to $\hat{\lambda}_1=\frac{\alpha_0 + x_1}{\beta_0+\bm{s}_{t_1}}$ and computes the optimal policy as if $\hat{\lambda}_1$ is the true base rate. This policy is then implemented until the next re-optimization epoch, after which the procedure is repeated until the scheduled maintenance moment. Figure \ref{illucep} provides a schematic illustration of the CE policy with $N^{\text{opt}}=1$. Recall that $u_{{\lambda}}^*$ denotes the optimal policy for a system with known base rate $\lambda$.

The CE policy is easy to implement and its behavior is well-understood from the structural properties established in Theorem \ref{optPolicy}. Specifically, if the DM learns that a system has a high (low) base rate, she will decrease (increase) the production in all states. The policy thus learns when it is possible to increase production to accumulate more revenue, but also when lower production is needed to prevent failures and high associated maintenance costs. 

\subsection{Performance of the Certainty-Equivalent Policy}
\label{subsec:performanceCE}
In the numerical studies in \S \ref{sec:num} we were able to numerically compute the performance measures needed for the comparisons in those studies. The performance of the CE policy, however, can only be evaluated through simulation. In this section and Appendix \ref{subsec:cePolicyApp} we report on the results of a comprehensive simulation study in which we compare the performance of the CE policy with an Oracle policy that knows the true base rate upon installation. 

We conduct a full factorial experiment in which we vary all input parameters across various values. The test bed of this experiment closely resembles the ones used in \S \ref{sec:num}, with the addition of two input parameters that are specific to the setting studied here. First,
in a simulation instance (we shortly describe the simulation procedure), the base rates are sampled from a Gamma distribution in which we fix the mean to 1, and vary its coefficient of variation, denoted with $cv_{\Lambda}$. This coefficient of variation is a measure for the unit-to-unit variability in the base rates of the systems. Second, the CE policy requires as input the parameter $N^{\text{opt}}$, which denotes how often the DM re-optimizes the policy based on the accrued knowledge. 
Permuting all input parameter values lead to a total of 4050 instances, see Table \ref{tab:testbedSim} in Appendix \ref{subsec:cePolicyApp} for a complete overview.

We benchmark the performance of the certainty equivalent policy, denoted with $\pi_{CEP}$, by the notion of relative regret, i.e. the percentage of expected profit loss (over a planning horizon) relative to the Oracle policy $\pi_{O}$. We compute this
relative regret as $\mathcal{\bar{R}} = 100\cdot (\bar{P}_{\pi_{O}}-\bar{P}_{\pi_{CEP}})/\bar{P}_{\pi_{O}}$, where $\bar{P}_\pi$ is the average profit of policy $\pi\in\{\pi_{CEP},\pi_{O}\}$. Since $\pi_{O}$ is unattainable in  practice, $\mathcal{R}$ serves as an upper bound on the relative regret with the best achievable performance.

The simulation procedure is as follows. For each instance, we first sample a base rate, say $\lambda^*$, from the Gamma distribution $\Lambda\sim\Gamma(\alpha_0,\beta_0)$. Then we compute the expected optimal profit of the Oracle policy using base rate $\lambda^*$:  $J^*_{\lambda^*}(0,T)$. We then simulate a sample path of an inhomogeneous Poisson process with base rate $\lambda^*$ that is controlled by the CE policy using the well-known thinning algorithm of \cite{lewis1979simulation} and obtain the resulting profit. We repeat this 2,000 times to obtain estimates for both $\bar{P}_{\pi_{O}}$ and $\bar{P}_{\pi_{CEP}}$, respectively, with the width of the 95\% confidence intervals being less than 0.1\% of the point estimates. We then compute $\mathcal{\bar{R}}$ using $\bar{P}_{\pi_{O}}$ and $\bar{P}_{\pi_{CEP}}$. 
Table \ref{tab:resultsTestbedMain3} summarizes the simulation study's main results. This table presents the average and maximum $\mathcal{\bar{R}}$ and the $\bar{P}_{\pi_{O}}$ across all instances. Detailed results can be found in Appendix \ref{subsec:cePolicyApp}.

Table \ref{tab:resultsTestbedMain3} indicates that the CE policy performs excellently with average profit losses of only 0.72\% compared to the Oracle policy across a wide range of instances. Additionally, with a maximum profit loss of only 11.63\% over all 4050 instances, the CE policy is also robust, providing good performance under varying conditions. 

\begin{table}[htp!]
\centering
\caption{Main results certainty-equivalent policy.}
\fontsize{9pt}{9pt}\selectfont
\label{tab:resultsTestbedMain3}
\begin{tabular}{ccc}
\toprule
 Mean {$\bar{\mathcal{R}}$}  & Max {$\bar{\mathcal{R}}$} & $\bar{P}_{\pi_{O}}$ \\
\midrule
0.72 & 11.63 & 10.33 \\
    \bottomrule
    \end{tabular}%
\end{table}%
\vspace{-8pt} 

It is important to mention that the detailed results in Table \ref{tab:resultsTestbed} in Appendix \ref{subsec:cePolicyApp} indicate that a concave deterioration function and a convex revenue function lead to particularly small profit losses. This is precisely the regime (cf. Corollary \ref{extremalConvexConcave} $(iii)$) where optimal condition-based production policies, and thus the CE and Oracle policy, are bang-bang policies (recall that this policy characterization does not depend on the value of $\lambda$).  Hence, these results suggest that CE policies are robust against unit-to-unit variability in the bang-bang regime; that is, they achieve low relative profit losses in that regime compared to the Oracle policy regardless of the degree of unit-to-unit variability.  In the next section, we report on additional numerical analyses that provide strong further evidence for this finding. 


\subsection{Robustness of Certainty-Equivalent Bang-Bang Policies} \label{subsec:robustCE}
\begin{figure}[!htb]
\begin{center}
\begin{subfigure}{.33\textwidth}
    \caption{}
\includegraphics[width=0.97\textwidth,keepaspectratio]{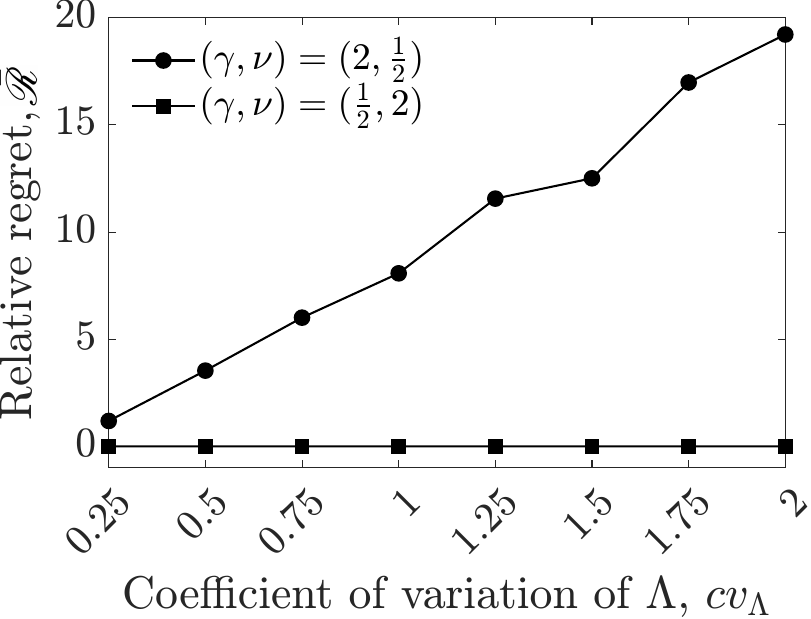}
  \label{fig:robust_2}
\end{subfigure}
\begin{subfigure}{.33\textwidth}\centering
\caption{}
\includegraphics[width=0.93\textwidth,keepaspectratio]{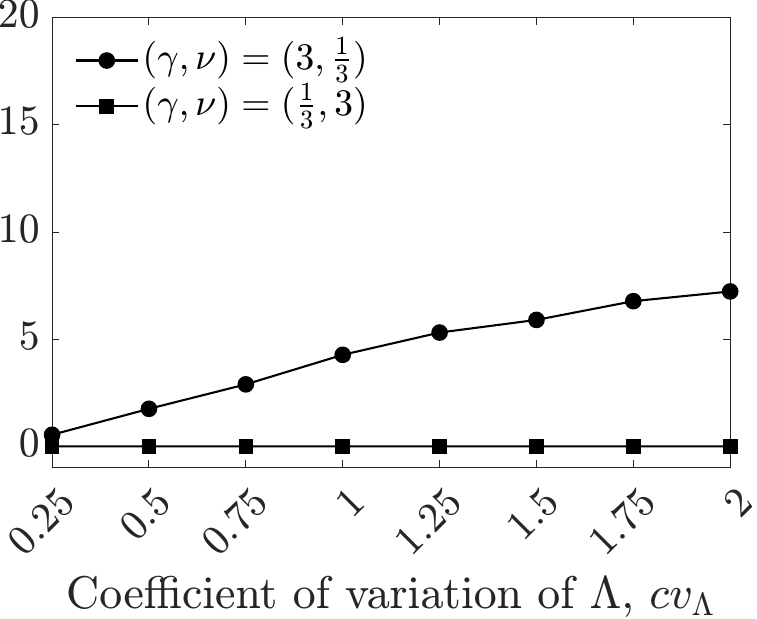}
  \label{fig:robust_3}
\end{subfigure}%
\begin{subfigure}{.33\textwidth}
  \centering
    \caption{}
\includegraphics[width=0.94\textwidth,keepaspectratio]{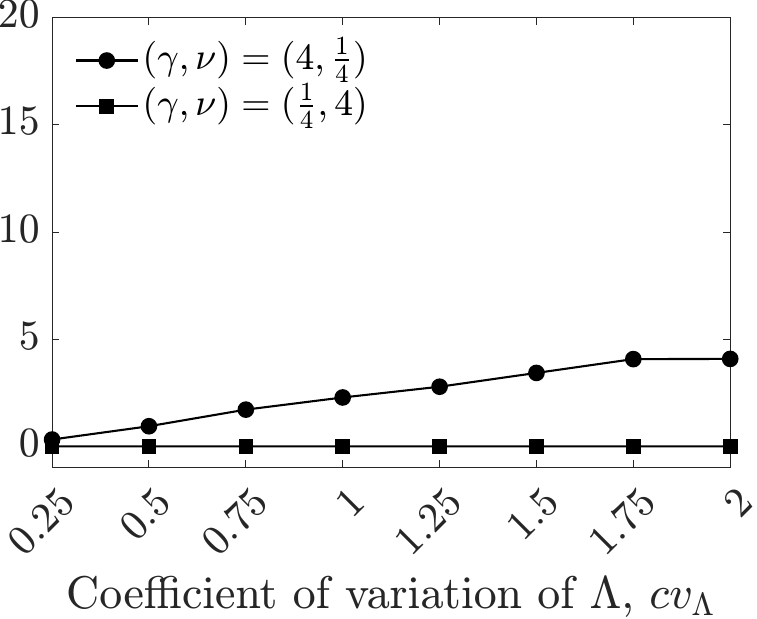}
  \label{fig:robust_4}
\end{subfigure}
\caption{Comparisons of the relative profit loss (in \%) of the CE policy with $N^{\text{opt}}=0$ as a function of $cv_{\Lambda}$ for three pairs of $(\gamma,\nu)$-values with $f=s^{\gamma}$ and $r=s^{\nu}$. In each comparison, the relative profit loss of the CE bang-bang policy is depicted with the solid squares.}
\label{robustCEBang}
\end{center}
\end{figure}

Figure \ref{robustCEBang} contains plots that exemplify behaviour we consistently observed in our numerical experiments. 
In these plots, we use the base case instance from \S \ref{subsec:cbpparametric}, except that we compare pairs of $(\gamma,\nu)$ values, fix $\E[\Lambda]=1$, and vary $cv_{\Lambda}.$ We also keep $N^{\text{opt}}$ fixed at 0 to isolate and study the effect of the bang-bang structure -- and not the effect of re-optimization -- on the relative profit loss $\bar{\mathcal{R}}$ (which we compute using the simulation procedure described in the previous subsection). In each plot, one $(\gamma,\nu)$ value corresponds to a bang-bang policy, while the other does not. The plots clearly indicate that $\bar{\mathcal{R}}$ is close to zero when the CE policy is bang-bang (squares) for all values of $cv_{\Lambda}$ while it increases in $cv_{\Lambda}$ when the CE policy is not bang-bang (circles).  

The robustness of a CE bang-bang policy against unit-to-unit variability can be explained as follows. When the Oracle is bang-bang, the CE policy is as well -- and vice versa -- as this only depends on $f(\cdot)$ and $r(\cdot)$, and not on (the realization of) the base rate $\lambda$. Since a bang-bang policy can only take two values ($0$ or $\bar{s}$) in each state $(x,t)$, it is much more likely that the CE and Oracle policy coincide throughout each sample path realization, than when the policy is not bang-bang in which the policy can take on infinitely many values in each state $(x,t)$. As a result, the expected profits of the Oracle and CE policy in the bang-bang regime will also be of similar magnitude. Indeed, we observe in Figure \ref{robustCEBang} that the relative profit loss of the CE policy is not impacted by unit-to-unit variability in the bang-bang regime.

This finding has an important managerial implication for DMs facing unit-to-unit variability in condition-based production systems. When DMs operate in the regime where a bang-bang policy is optimal (which we established can be verified without knowing the value of $\lambda$, see Proposition \ref{bangbangCondition}), DMs do not have to take variability in deterioration rates into account, but only require the population's average deterioration rate $\E[\Lambda]$. As DMs in practice can often readily obtain an unbiased estimate of $\E[\Lambda]$ through data, implementing a CE policy with $N^{\text{opt}}=0$ is not only computationally less demanding and straightforward, it also ensures superior performance regardless of the realization of $\lambda$. When the policy is not bang-bang, DMs can rely on a CE policy with $N^{\text{opt}}>0$ to learn $\lambda$ and still effectively manage condition-based production systems that have unit-to-unit variability. The numerical results in \S\ref{subsec:performanceCE} suggest that the CE policy is namely capable of achieving good overall performance when $N^{\text{opt}}>0$.  

\section{Concluding Remarks}\label{sec:conclusion}

Building on advancements in sensor technology and Internet-of-Things connectivity, we addressed two interconnected decision problems for stochastically deteriorating production systems, focusing on real-time monitoring and control.
We first studied the operational decision problem of how to adjust the production rate to balance deterioration and revenue accumulation based on the real-time condition signal given that maintenance moments are fixed. We modeled this problem as a continuous-time MDP, and established intuitive monotonic properties of the optimal policy. We also presented easily verifiable conditions that are sufficient for a bang-bang policy to be optimal. 
Next, we considered the tactical decision problem of finding the optimal planned maintenance moments given that the optimal condition-based production policy is followed in between. We showed that this decision problem can be easily optimized, thereby enabling the optimal integration of condition-based production on an operational level and maintenance decisions on a tactical level.
Finally, we extended our analysis to populations of production systems that exhibit unit-to-unit variability in the sense that the production-deterioration relation is not known and needs to be learned on-the-fly. For the operational decision problem in this setting, we proposed a Bayesian procedure to learn the unknown deterioration rate under any production policy in real time. 

Our extensive numerical studies have led to several important insights. Firstly, condition-based production policies enabled by real-time condition monitoring and remote control significantly outperform static policies. Secondly, the integration of tactical maintenance planning with operational condition-based production policies leads to substantial profit increases, as opposed to addressing these two problems sequentially. This underscores the importance of setting maintenance intervals in conjunction with operational condition-based production decisions. Thirdly, in cases where production systems exhibit unit-to-unit variability in the production-deterioration relationship, our Bayesian heuristic policy performs close to an Oracle policy. Notably, our findings highlight a key insight: The heuristic policy demonstrates particularly superior performance in the bang-bang regime, which remains robust even in the presence of high unit-to-unit variability.



\bibliographystyle{informs2014}
\setlength{\bibsep}{-2pt}
\renewcommand*{\bibfont}{\footnotesize}
\bibliography{myrefs}
\newpage
 \begin{APPENDICES}
 
\section{Proofs of Statements}\label{AppendixProofs}
\subsection{Proof of Lemma \ref{lemmaHJB}}
\proof{Proof.}
We use a standard argument in which we consider the dynamics of the deterioration under a chosen setting -- essentially an inhomogeneous Poisson process -- over a very small time interval.  

Consider the dynamics over a small interval of length $\delta t$ (that is, when the time until the next maintenance moment changes from $t$ to $t-\delta t$) when $t>0$. Let $x < \xi$. If  the decision maker selects production rate $s\in\mathscr{S}$, then with probability $\delta t \lambda f(s) + o(\delta t)$, $x$ changes to $x+1$, and with probability $1-\delta t \lambda f(s) +o(\delta t)$, $x$ stays the same. Here, $o(\delta t)$ is a quantity that goes to zero faster than $\delta t$. Secondly, if $x=\xi$, then the only admissible action is $s=0$, so that with probability 1, $x$ stays $\xi$. In all cases, $t$ changes to $t-\delta t$. By the principle of optimality, for $t>0$, we have
\begin{align}\label{heuristicx}
J_{\lambda}^*(x,t) = \max_{s\in \mathscr{S}}\Big[ r(s)\delta t + \delta t \lambda f(s) J_{\lambda}^*(x+1,t-\delta t)  +  \big(1-\delta t \lambda f(s)\big)J_{\lambda}^*(x,t-\delta t) + o(\delta t)\Big], 
\end{align}
when $x< \xi $ and 
\begin{align}\label{heuristicx0}
J_{\lambda}^*(\xi,t) =  J_{\lambda}^*(\xi,t-\delta t), 
\end{align}
when $x=\xi$. Rearranging \eqref{heuristicx} and \eqref{heuristicx0}, dividing by $\delta t$, and taking the limit $\delta t \rightarrow 0$, leads to the desired result. The boundary conditions follow from the maintenance cost that is incurred when the system has deterioration  level $x$ at the scheduled maintenance moment. \hfill \Halmos \endproof

\subsection{Proof of Lemma \ref{lemmaIncreasing}}
\proof{Proof.} We prove each assertion of the statement using a sample path argument. For the monotonicity in $x$, consider a sample path where the system with deterioration level $x$ is governed by the optimal condition-based production policy for the system with deterioration level $x+1$ until either the system with deterioration level $x+1$ breaks down, or until the scheduled maintenance moment, whichever happens first. Note that since the policies are coupled, they will face the same deterioration sample path. If the system with deterioration level $x+1$ breaks down at time $\tau$ before the scheduled maintenance moment, then the system with deterioration level $x$ can still generate revenue over $(T-\tau,0]$,  possibly even with less maintenance costs if the system does not break down. On the other hand, if the scheduled maintenance moment arrives first, the two systems will have generated the same revenue and will have the same maintenance cost. Consequently, the system with deterioration level $x$ generates at least as much revenue and incurs at most the same maintenance costs. Hence, for their expected profits, we have that $J_{\lambda}^*(x,t)\geq J_{\lambda}^*(x+1,t)$. Using similar arguments one can also prove the monotonicity in $t$ and $\lambda$.
\hfill \Halmos \endproof

\subsection{Proof of Proposition \ref{concaveInX}}
\proof{Proof.}
We first introduce some additional notation. Let $\Delta_{\lambda}^2(x,t) \triangleq \Delta_{\lambda}(x,t) - \Delta_{\lambda}(x+1,t)$ for $x \leq \xi-2$. To establish the concavity of $J_{\lambda}^*(x,t)$ in $x$ for all $t\geq 0$ and each base rate,  we need to show that $\Delta_{\lambda}^2(x,t)\leq 0$ for all  $x \leq \xi-2$ and $t,\lambda \geq 0$.

For $t=0$, it is easy to verify that $\Delta_{\lambda}^2(x,0)\leq 0$ for  all  $x \leq \xi-2$ and $\lambda \geq 0$, because $-c_m(x)$ is concave in $x$ by assumption.


We now consider $t>0$.  We first construct some useful inequalities that we will use in an inductive proof to establish the result. Let $x\leq \xi-2$, using \eqref{HJEqn}, we have 
\begin{align}\label{eqn1}
\frac{\partial J_{\lambda}^*(x,t)}{\partial t} = r\big( s_{\lambda}^*(x,t) \big) - \lambda f\big( s_{\lambda}^*(x,t) \big) \Delta_{\lambda}(x,t),
\end{align}
and since $s_{\lambda}^*(x,t)$ is also admissible, but not necessarily optimal, for state $(x+1,t)$ with base rate $\lambda$, we have 
\begin{align}\label{eqn2}
\frac{\partial J_{\lambda}^*(x+1,t)}{\partial t} \geq  r\big( s_{\lambda}^*(x,t) \big) - \lambda f\big( s_{\lambda}^*(x,t) \big) \Delta_{\lambda}(x+1,t).    
\end{align}
Then, subtracting \eqref{eqn2} from \eqref{eqn1} leads to 
\begin{align} \label{eqn3}
\frac{\partial \Delta_{\lambda} (x,t)}{\partial t}  \leq - \lambda f\big( s_{\lambda}^*(x,t) \big) \Delta_{\lambda}^2(x,t).  
\end{align}
Similarly, we have 
\begin{align}\label{eqn5}
\frac{\partial J_{\lambda}^*(x+2,t)}{\partial t} = r\big( s_{\lambda}^*(x+2,t) \big) - \lambda f\big( s_{\lambda}^*(x+2,t) \big) \Delta_{\lambda}(x+2,t),
\end{align}
and since $s_{\lambda}^*(x+2,t)$ is also admissible, but not necessarily optimal, for state $(x+1,t)$ with base rate $\lambda$, we have 
\begin{align}\label{eqn4}
\frac{\partial J_{\lambda}^*(x+1,t)}{\partial t} \geq  r\big( s_{\lambda}^*(x+2,t) \big) - \lambda f\big( s_{\lambda}^*(x+2,t) \big) \Delta_{\lambda}(x+1,t).    
\end{align}
Then, subtracting \eqref{eqn5} from \eqref{eqn4} gives
\begin{align} \label{eqn6}
\frac{\partial \Delta_{\lambda} (x+1,t)}{\partial t}  \geq - \lambda f\big( s_{\lambda}^*(x+2,t) \big) \Delta_{\lambda}^2(x+1,t).  
\end{align}
Next, subtracting \eqref{eqn6} from \eqref{eqn3} leads to 
\begin{align} \label{eqn7a}
\frac{\partial \Delta_{\lambda}^2(x,t)}{\partial t}  \leq - \lambda f\big( s_{\lambda}^*(x,t) \big) \Delta_{\lambda}^2(x,t) + \lambda f\big( s_{\lambda}^*(x+2,t) \big) \Delta_{\lambda}^2(x+1,t).  
\end{align}
We now proceed with our inductive proof. We first prove the base case, i.e. $\Delta_{\lambda}^2(\xi-2,t)\leq 0$ for all $t>0$ and $\lambda \geq 0$. By \eqref{eqn3} we have that 
\begin{align}\label{eqn7}
\frac{\partial \Delta_{\lambda} (\xi-2,t)}{\partial t}  \leq - \lambda f\big( s_{\lambda}^*(\xi-2,t) \big) \Delta_{\lambda}^2(\xi-2,t).
\end{align} 
Note that due to the linearity of the derivative, we have  
\begin{align}\label{eqn8}
\frac{\partial \Delta_{\lambda} (\xi-1,t)}{\partial t} &= \frac{\partial J_{\lambda}^*(\xi-1,t)}{\partial t} - \frac{\partial J_{\lambda}^*(\xi,t)}{\partial t} 
= \frac{\partial J_{\lambda}^*(\xi-1,t)}{\partial t}  \geq 0,
\end{align} 
where the second equality follows from \eqref{HJEqn} for $x=\xi$ and  the inequality follows from Lemma \ref{lemmaIncreasing}. Equation \eqref{eqn7} minus Equation \eqref{eqn8} on both sides leads to 
\begin{align}\label{eqn9}
\frac{\partial \Delta_{\lambda}^2(\xi-2,t)}{\partial t}  \leq - \lambda f\big( s_{\lambda}^*(\xi-2,t) \big) \Delta_{\lambda}^2(\xi-2,t).
\end{align} 
Applying Gr\"onwall's Lemma to Inequality \eqref{eqn9} leads to 
\begin{align}
\Delta_{\lambda}^2(\xi-2,t) \leq  \Delta_{\lambda}^2(\xi-2,0)\cdot \exp \left({\int_0^t-\lambda f( s_{\lambda}^*(\xi-2,u) ) d u }\right).
\end{align} 
Since $\Delta_{\lambda}^2(\xi-2,0) \leq 0$, (recall that $-c_m(x)$ is concave in $x$ by assumption), and because $\exp \left( {\int_0^t-\lambda f( s_{\lambda}^*(\xi-2,u) ) d u }\right) >0$ for all $t> 0$ and $\lambda \geq0$, we have that  $\Delta_{\lambda}^2(\xi-2,t)\leq 0$, which proves the base case. 

Assume now inductively that  $\Delta_{\lambda}^2(x,t)\leq 0$ for an $x\leq \xi - 2$. We will show that this implies that $\Delta_{\lambda}^2(x-1,t)\leq 0$ for all $t>0$ and all $\lambda\geq0$. 

For $(x-1,t)$, we have 
\begin{align}\label{eqn10}
\frac{\partial \Delta_{\lambda}^2(x-1,t)}{\partial t} & \leq - \lambda f\big( s_{\lambda}^*(x-1,t) \big) \Delta_{\lambda}^2(x-1,t) + \lambda f\big( s_{\lambda}^*(x+1,t) \big) \Delta_{\lambda}^2(x,t) \nonumber \\
&\leq  - \lambda f\big( s_{\lambda}^*(x-1,t) \big) \Delta_{\lambda}^2(x-1,t).
\end{align}
The first inequality is Equation \eqref{eqn7a}, and the second inequality is due to the induction hypothesis.
Applying Gr\"onwall's Lemma to Inequality \eqref{eqn10} leads to 
\begin{align}
 \Delta_{\lambda}^2(x-1,t) \leq  \Delta_{\lambda}^2(x-1,0)\cdot \exp \left( {\int_0^t-\lambda f( s_{\lambda}^*(x-1,u) ) d u }\right).  
\end{align}
Since $\Delta_{\lambda}^2(x-1,0) \leq 0$ for all $x\leq \xi-2$ (again by the concavity assumption of $-c_m(x)$ in $x$)  and because $ \exp\left( {\int_0^t-\lambda f( s_{\lambda}^*(x-1,u) ) d u }\right) >0$ for all $t> 0$ and $\lambda \geq0$, 
we have that  $\Delta_{\lambda}^2(x-1,t)\leq 0$ for all $t>0$ and all $\lambda\geq0$. 
\hfill \Halmos \endproof

\subsection{Proof of Proposition \ref{increasingInT}}
\proof{Proof.} In this proof, we employ a commonly used transformation of the continuous model into an analog discrete dynamic program \citep[see e.g.,][]{bitran1997periodic,aydin2008pricing,hu2017optimal}. We are allowed to do so because it is known that this discrete dynamic program is guaranteed to converge uniformly to the continuous model \citep[cf.][]{kleywegt2001dynamic}.

To this end, we divide the entire planning horizon into $\hat{T} = \frac{T}{\delta t}$ periods, each of which is short enough that 
the probability of more than one deterioration increment during an interval of length $\delta t$ is negligible. Indeed, note that since the function $f(\cdot)$ is bounded, we have that for each choice of production rate $s$, the probability that more than 1 arrival will occur in a small interval of length $\delta t$ is of order $o(\delta t)$. 

Let $J_{\lambda}^*(x,k)$ denote the maximum expected revenue when the deterioration level is $x$ and there are $k$ periods to go until the planned maintenance moment, and the base rate is $\lambda$. We then need to show that the following inequality holds
\begin{align}\label{toshow}
J_{\lambda}^*(x-1,k)-J_{\lambda}^*(x,k) \geq J_{\lambda}^*(x-1,k-1) - J_{\lambda}^*(x,k-1).
\end{align}

The dynamic program for the revenue maximization problem over the $\hat{T}$-period horizon is given by
 \begin{align}\label{optDiscrete}
 J_{\lambda}^*(x,k) = \max_{s\in\mathscr{S}}\Big[ r(s)\delta t + \lambda f(s) \delta t  J_{\lambda}^*(x+1,k-1) + \big(1 - \lambda f(s) \delta t \big)J_{\lambda}^*(x,k-1) \Big],
 \end{align}
 with boundary conditions  $J_{\lambda}^*(x,0) = -c_m(x)$ for all $x\in \mathscr{X}$.
 
Let $s_{\lambda}^*(x,k)$ denote the optimal production rate for state $(x,k)$ when the base rate is $\lambda$. We then have, using \eqref{optDiscrete},
\begin{align}\label{optPeriod}
J_{\lambda}^*(x,k) = r\big(s_{\lambda}^*(x,k)\big) \delta t + \lambda f\big(s_{\lambda}^*(x,k)\big)\delta t J_{\lambda}^*(x+1,k-1) + \big(1-\lambda f\big(s_{\lambda}^*(x,k)\big)\delta t\big)J_{\lambda}^*(x,k-1).
\end{align}
Subtracting $J_{\lambda}^*(x,k-1)$ from both sides of \eqref{optPeriod} leads to the following equality
\begin{align}\label{optPeriodSubtracted}
J_{\lambda}^*(x,k) - J_{\lambda}^*(x,k-1) = r\big(s_{\lambda}^*(x,k)\big) \delta t + \lambda f\big(s_{\lambda}^*(x,k)\big)\delta t \big(J_{\lambda}^*(x+1,k-1)-J_{\lambda}^*(x,k-1)\big).
\end{align}
Since $s_{\lambda}^*(x,k)$ is admissible in state $(x-1,k)$ when the base rate is $\lambda$ but not necessarily optimal, we obtain the following inequality
\begin{align}\label{ineq1}
J_{\lambda}^*(x-1,k) & \geq  \nonumber\\ & r\big(s_{\lambda}^*(x,k)\big) \delta t + \lambda f\big(s_{\lambda}^*(x,k)\big)\delta t J_{\lambda}^*(x,k-1) + \big(1-\lambda f\big(s_{\lambda}^*(x,k)\big)\delta t\big)J_{\lambda}^*(x-1,k-1).
\end{align}
Subtracting $J_{\lambda}^*(x-1,k-1)$ from both sides of \eqref{ineq1} and then using the concavity of $J_{\lambda}^*(x,k)$ in $x$ (see Proposition \ref{concaveInX}) yields
\begin{align}
J_{\lambda}^*(x-1,k)-J_{\lambda}^*(x-1,k-1)  &\nonumber\\
\geq\ & r\big(s_{\lambda}^*(x,k)\big) \delta t + \lambda f\big(s_{\lambda}^*(x,k)\big)\delta t  \big(J_{\lambda}^*(x,k-1) - J_{\lambda}^*(x-1,k-1)\big) \nonumber \\
 \geq\ & r\big(s_{\lambda}^*(x,k)\big) \delta t + \lambda f\big(s_{\lambda}^*(x,k)\big)\delta t  \big(J_{\lambda}^*(x+1,k-1) - J_{\lambda}^*(x,k-1)\big). \label{ineq2}
\end{align}
Observe that the right hand side of \eqref{ineq2} can be replaced by \eqref{optPeriodSubtracted}, which yields the following inequality 
\begin{align}
J_{\lambda}^*(x-1,k)-J_{\lambda}^*(x-1,k-1)\geq J_{\lambda}^*(x,k) - J_{\lambda}^*(x,k-1) ,\nonumber
\end{align}
which, after rearranging terms, equals \eqref{toshow}.
%
\hfill \Halmos \endproof

\subsection{Proof of Proposition \ref{concaveInT}}
\proof{Proof.} We need to show that $\frac{\partial J_{\lambda}^*(x,t)}{\partial t}$ is decreasing in $t$. Let $u>t$. We have 
\begin{align}\label{eqn1concave}
\frac{\partial J_{\lambda}^*(x,u)}{\partial u} = r\big( s_{\lambda}^*(x,u) \big) - \lambda f\big( s_{\lambda}^*(x,u) \big) \Delta_{\lambda}(x,u),
\end{align}
and since $s_{\lambda}^*(x,u)$ is also admissible, but not necessarily optimal (recall that $s_{\lambda}^*(x,u)$ is the optimal production rate for state $(x,u)$ when the base rate is $\lambda$), for state $(x,t)$ when the base rate is $\lambda$, we have 
\begin{align}\label{eqn2concave}
\frac{\partial J_{\lambda}^*(x,t)}{\partial t} \geq  r\big( s_{\lambda}^*(x,u) \big) - \lambda f\big( s_{\lambda}^*(x,u) \big) \Delta_{\lambda}(x,t).    
\end{align}
Subtracting \eqref{eqn2concave} from \eqref{eqn1concave} yields 
\begin{align*}
\frac{\partial J_{\lambda}^*(x,t)}{\partial t} - \frac{\partial J_{\lambda}^*(x,u)}{\partial u} \geq \lambda f\big( s_{\lambda}^*(x,u) \big) \Big(\Delta_{\lambda}(x,u) - \Delta_{\lambda}(x,t) \Big) \geq 0,
\end{align*}
where the last inequality is because $\Delta_{\lambda}(x,u) \geq \Delta_{\lambda}(x,t)$ for $u>t$ (see Proposition \ref{increasingInT}).
\hfill \Halmos \endproof

\subsection{Proof of Lemma \ref{submodExpect}}
\proof{Proof.} 
Let $x_1\leq x_2$ and $\theta_1\leq \theta_2$. Since $Z(\theta_1) \leqst{} Z(\theta_2)$, there exist two random variables, $\tilde{Z}(\theta_1)$ and  $\tilde{Z}(\theta_2)$, on the same probability space, such that  $\tilde{Z}(\theta_1) \eqst{} Z(\theta_1)$,  $\tilde{Z}(\theta_2) \eqst{} Z(\theta_2)$, and $\tilde{Z}(\theta_1) \leq \tilde{Z}(\theta_2)$ almost surely (here, $\eqst{}$ denotes equality in law) \citep[see, e.g.,][Theorem 1.A.1.]{shaked2007stochastic}. We therefore have: 
\begin{align*}
 \E[f(x_1+Z(\theta_1),\theta_1)] - \E[f(x_2+Z(\theta_1),\theta_1)] 
&=  \E[f(x_1+\tilde{Z}(\theta_1),\theta_1) - f(x_2+\tilde{Z}(\theta_1),\theta_1)] &\\
 &\leq \E[f(x_1+\tilde{Z}(\theta_1),\theta_2) - f(x_2+\tilde{Z}(\theta_1),\theta_2)]&\\
 &\leq \E[f(x_1+\tilde{Z}(\theta_2),\theta_2) - f(x_2+\tilde{Z}(\theta_2),\theta_2)]&\\
 &=  \E[f(x_1+Z(\theta_2),\theta_2)] - \E[f(x_2+Z(\theta_2),\theta_2)],&
\end{align*}
where the first inequality is from the submodularity of $f(x,\theta)$, and the second inequality is from the concavity in $x$ of $f(x,\theta)$. Hence,  $\E[f(x+Z(\theta),\theta)]$ is submodular in $(x,\theta)$.

The preservation of concavity under expectation is a well-known result. For completeness, we provide a short proof here. Let $x_1,x_2\in \R$ and $\alpha\in[0,1]$. We then have 
\begin{align*}
 \E[f(\alpha x_1 + (1-\alpha)x_2 +Z(\theta),\theta)] &= \E[f(\alpha x_1 + (1-\alpha)x_2 + \alpha Z(\theta) +(1-\alpha)  Z(\theta),\theta)]\\
 &= \E[f(\alpha\big(x_1+Z(\theta)\big) + (1-\alpha)\big(x_2+Z(\theta)\big),\theta)]\\
 &\geq \E[\alpha f\big(x_1+Z(\theta), \theta\big) +  (1-\alpha)f\big(x_2+Z(\theta),\theta \big)]\\
 &= \alpha \E[f\big(x_1+Z(\theta), \theta\big)] +  (1-\alpha)\E[f\big(x_2+Z(\theta),\theta \big)],
\end{align*}
where the inequality follows from the concavity of $f(x,\theta)$ in $x$.
\hfill \Halmos \endproof

\subsection{Proof of Lemma \ref{submodFunction}}
\proof{Proof.} We need to show that for $x_1 \leq x_2$ and $y_1 \leq y_2$ the following inequality holds: $$g(y_1)f(x_1,y_1)-g(y_2)f(x_1,y_2)\leq g(y_1)f(x_2,y_1)-g(y_2)f(x_2,y_2).$$  To this end, let $x_1 \leq x_2$ and $y_1 \leq y_2$, we then have
\begin{align*}
&g(y_1)f(x_1,y_1)-g(y_2)f(x_1,y_2)-g(y_1)f(x_2,y_1)+g(y_2)f(x_2,y_2) & \\
&= g(y_1) \big( f(x_1,y_1) -f(x_2,y_1)\big) - g(y_2) \big( f(x_1,y_2)-f(x_2,y_2)\big)& \\
& \leq g(y_2) \big( f(x_1,y_1) -f(x_2,y_1) - f(x_1,y_2) + f(x_2,y_2) \big) \leq 0. 
\end{align*}
The first inequality is because $g(y_1)\leq g(y_2)$ and $f(x_1,y_1)\geq f(x_2,y_1)$, and the second inequality is due to the submodularity of $f(x,y)$ in $(x,y)$. 
\hfill \Halmos \endproof

\subsection{Proof of Proposition \ref{submodvalue}}
\proof{Proof.} Again, we use the analog discrete dynamic program for the revenue maximization problem over the $\hat{T}$-period horizon and rewrite it in expectation form:
 \begin{align}
 J_{\lambda}^*(x,k) = \max_{s\in\mathscr{S}}\Big[ r(s)\delta t +  \E[J_{\lambda}^*(x+Z(\lambda,s),k-1)] \Big],
 \end{align}
 with boundary conditions  $J_{\lambda}^*(x,0) = -c_m(x)$ for all $x\in \mathscr{X}$, and where $Z(\lambda,s)$ is a Bernoulli random variable with success probability $p= \delta t\lambda f(s)$. 
 
 We need to show that for all $k\geq 0$, $J_{\lambda}^*(x,k)$ is submodular in $(x,\lambda)$. We prove this by induction on $k$. Note that for the base case $k=0$,  $J_{\lambda}^*(x,0) = -c_m(x)$ is evidently submodular in $(x,\lambda)$. Assume that for some $k-1$, $J_{\lambda}(x,k-1)$ is submodular in $(x,\lambda)$, we will now show that $J_{\lambda}(x,k)$ is submodular in  $(x,\lambda)$. 
 
Let $g_{\lambda}(s,x,k) \triangleq  r(s)\delta t +  \E[J_{\lambda}^*(x+Z(\lambda,s),k-1)]$. Since $Z(\lambda,s)$ is stochastically increasing in $\lambda$, the induction hypothesis and the fact that $J_{\lambda}^*(x,k-1)$ is concave in $x$ (see Proposition \ref{concaveInX}), we can employ Lemma \ref{submodExpect} to conclude that $g_{\lambda}(s,x,k)$ is submodular in $(x,\lambda)$ for all $s$. 
 
 Now, let $s_{\lambda}^*(x,k)\in \arg\max_{s\in\mathscr{S}}\{g_{\lambda}(s,x,k)\}$, and let $x_1<x_2$, and $\lambda_1<\lambda_2$, we then have: 
 \begin{align*}
 J_{\lambda_1}^*(x_1,k) + J_{\lambda_2}^*(x_2,k) &= g_{\lambda_1}\big(s_{\lambda_1}^*(x_1,k),x_1,k\big)+g_{\lambda_2}\big(s_{\lambda_2}^*(x_2,k),x_2,k\big)\\
 &\leq g_{\lambda_2}\big(s_{\lambda_1}^*(x_1,k),x_1,k\big)+g_{\lambda_1}\big(s_{\lambda_2}^*(x_2,k),x_2,k\big)\\
 &\leq g_{\lambda_2}\big(s_{\lambda_2}^*(x_1,k),x_1,k\big)+g_{\lambda_1}\big(s_{\lambda_1}^*(x_2,k),x_2,k\big)\\
 &=  J_{\lambda_2}^*(x_1,k) + J_{\lambda_1}^*(x_2,k),
 \end{align*}
 where the first inequality is due to the submodularity of $g_{\lambda}(s,x,k)$ in $(x,\lambda)$ for all $s$, and the second inequality is due to the suboptimality of $s_{\lambda_1}^*(x_1,k)$ for $(x_1,\lambda_2)$ and  $s_{\lambda_2}^*(x_2,k)$ for $(x_2,\lambda_1)$, respectively. Hence $J_{\lambda}^*(x,k)$ is submodular in $(x,\lambda)$.
\hfill \Halmos \endproof

\subsection{Proof of Theorem \ref{optPolicy}}
\proof{Proof.} The proof is divided into three parts, focused on the monotonicity in $x$, $t$, and $\lambda$, respectively. In each part we use the notation $g_{\lambda}(s,x,t) = r(s) - \lambda f(s) \Delta_{\lambda}(x,t)$, that we previously introduced.
 \begin{description}
   \item[Proof of monotonicity in $x$:] For notational clarity we drop the subscript $\lambda$. Let $\hat{s}\triangleq s^*(x,t)$ and consider $s>\hat{s}$, we have that 
\begin{align}\label{ineqOptX}
&g(\hat{s},x+1,t) - g(\hat{s},x,t) + g(s,x,t) - g(s,x+1,t) \nonumber \\
&= -\lambda f(\hat{s})\Delta(x+1,t) + \lambda f(\hat{s})\Delta(x,t) - \lambda f(s)\Delta(x,t) + \lambda f(s)\Delta(x+1,t) \nonumber \\ 
&= \lambda f(\hat{s})\Delta^2(x,t) - \lambda f(s)\Delta^2(x,t) \nonumber \\
&= \lambda \big(f(\hat{s})-f(s)\big) \Delta^2(x,t) \geq 0.
\end{align}
The inequality follows from the increasing property of $f(\cdot)$ with $s>\hat{s}$, and since $\Delta^2(x,t)\leq 0$ (see Proposition \ref{concaveInX}). Rearranging Inequality \eqref{ineqOptX} leads to 
\begin{align}\label{ineqOptX2}
g(\hat{s},x+1,t) - g(s,x+1,t) \geq g(\hat{s},x,t) -g(s,x,t)\geq 0,
\end{align}
where the second inequality is because $\hat{s}\triangleq s^*(x,t)$. We can therefore conclude from \eqref{ineqOptX2}, that any $s>\hat{s}$ cannot be optimal for $(x+1,t)$, hence $s^*(x+1,t)\leq s^*(x,t)$.
\item[Proof of monotonicity in $t$:] For notational clarity we again drop the subscript $\lambda$. Let $u<t$ and $\hat{s}\triangleq s^*(x,u)$, and consider $s>\hat{s}$. We have that \begin{align}\label{ineqOptT}
&g(\hat{s},x,t) - g(\hat{s},x,u) + g(s,x,u) - g(s,x,t) \nonumber \\
&= -\lambda f(\hat{s})\Delta(x,t) + \lambda f(\hat{s})\Delta(x,u) - \lambda f(s)\Delta(x,u) + \lambda f(s)\Delta(x,t) \nonumber \\ 
&= \lambda\big(f(\hat{s})- f(s)\big)\big(\Delta(x,u)-\Delta(x,t)\big) \geq 0.
\end{align}
The inequality follows from the increasing property of $f(\cdot)$ and $s>\hat{s}$, and since $\Delta(x,u)\leq \Delta(x,t)$ for $u<t$ (see Proposition \ref{increasingInT}). Rearranging Inequality \eqref{ineqOptT} leads to 
\begin{align}\label{ineqOptT2}
g(\hat{s},x,t) - g(s,x,t) \geq g(\hat{s},x,u) -g(s,x,u)\geq 0,
\end{align}
where the second inequality is because $\hat{s}\triangleq s^*(x,u)$. We can therefore conclude from \eqref{ineqOptT2}, that any $s>\hat{s}$ cannot be optimal for $(x,t)$, hence $s^*(x,t)\leq s^*(x,u)$ for $u<t$.
\item[Proof of monotonicity in $\lambda$:] In this part, we use the analog discrete dynamic program. In light of the well-known Topkis' Theorem \citep[see, for instance,][Theorem 2.8.1]{topkis2011supermodularity}, we need to show that $g_{\lambda}(s,x,k) =  r(s)\delta t + \lambda f(s) \delta t J_{\lambda}^*(x+1,k-1) + (1-\lambda f(s) \delta t)J_{\lambda}^*(x,k-1)$ is submodular in $(s,\lambda)$ to conclude that  $s_{\lambda}^*(x,k)$ is decreasing in $\lambda$.   

Let $s_1\leq s_2$ and $\lambda_1 \leq \lambda_2$. We then need to show that $$g_{\lambda_1}(s_1,x,k) - g_{\lambda_2}(s_1,x,k)-g_{\lambda_1}(s_2,x,k) + g_{\lambda_2}(s_2,x,k) \leq 0.$$
We have (after some algebraic manipulations):
\begin{align*}
   & g_{\lambda_1}(s_1,x,k) - g_{\lambda_2}(s_1,x,k)-g_{\lambda_1}(s_2,x,k) + g_{\lambda_2}(s_2,x,k) &  \\
    &= \lambda_1 f(s_1)  J_{\lambda_1}^*(x+1,k-1) - \lambda_1 f(s_1) J_{\lambda_1}^*(x,k-1) &\\
    &- \lambda_2 f(s_1) J_{\lambda_2}^*(x+1,k-1) + \lambda_2 f(s_1)  J_{\lambda_2}^*(x,k-1) &\\
    &- \lambda_1 f(s_2) J_{\lambda_1}^*(x+1,k-1) + \lambda_1 f(s_2) J_{\lambda_1}^*(x,k-1) &\\ 
    &+ \lambda_2 f(s_2) J_{\lambda_2}^*(x+1,k-1)  - \lambda_2 f(s_2)  J_{\lambda_2}^*(x,k-1) &\\
    &= \big(f(s_1)-f(s_2)\big)\big( \lambda_1  J_{\lambda_1}^*(x+1,k-1) - \lambda_1J_{\lambda_1}^*(x,k-1) -\lambda_2  J_{\lambda_2}^*(x+1,k-1) + \lambda_2 J_{\lambda_2}^*(x,k-1)  \big). & 
\end{align*}
Since $f(s_1)-f(s_2)\leq 0$ ($f(\cdot)$ is increasing and $s_2\geq s_1$), $g_{\lambda}(s,x,k)$ is submodular in $(s,\lambda)$ if  $$\lambda_1  J_{\lambda_1}^*(x+1,k-1) - \lambda_1J_{\lambda_1}^*(x,k-1) -\lambda_2  J_{\lambda_2}^*(x+1,k-1) + \lambda_2 J_{\lambda_2}^*(x,k-1) \geq 0,$$
that is, if $\lambda J_{\lambda}^*(x,k-1)$ is submodular in $(x,\lambda)$.

Note that $J^*_{\lambda}(x,k-1)$ is submodular in $(x,\lambda)$ (Proposition \ref{submodvalue}) and decreasing in $x$ (Lemma \ref{lemmaIncreasing}), hence by applying Lemma \ref{submodFunction} with $g(\lambda)\triangleq \lambda$ and $f(x,\lambda)\triangleq J^*_{\lambda}(x,k-1)$, we can conclude that $\lambda J_{\lambda}^*(x,k-1)$ is submodular in $(x,\lambda)$.
\end{description}
\hfill \Halmos \endproof

\subsection{Proof of Proposition \ref{bangbangCondition}}
\proof{Proof.} 
Recall that for a state $(x,t)$, we are interested in the following optimization problem: 
\begin{align} \max_{s\in [0,\bar{s}]} \Big[r(s) -  \lambda f(s)\Delta_{\lambda}(x,t) \Big]. \label{maxproblem}
\end{align} 
Suppose the contrary of the proposition; that is, if $\frac{r(s)}{f(s)} \leq  \frac{r(\bar{s})}{f(\bar{s})}$ for all $s\in (0,\bar{s})$ then the optimal policy has not a bang-bang structure. Let $s\in (0,\bar{s})$ be the optimal production rate. We will now show that this leads to a contradiction. 

Since $s\in (0,\bar{s})$ is optimal, we have by Equation \eqref{maxproblem}: 
\begin{align}\nonumber
r(s) -  \lambda f(s)\Delta_{\lambda}(x,t) & > r(0) -  \lambda f(0)\Delta_{\lambda}(x,t), \text{ and} \\
r(s) -  \lambda f(s)\Delta_{\lambda}(x,t) & > r(\bar{s}) -  \lambda f(\bar{s})\Delta_{\lambda}(x,t).\nonumber
\end{align}
These inequalities are equivalent to, respectively, 
\begin{align} \label{eqn1Bang}
r(s) & > r(0) +  \lambda \Delta_{\lambda}(x,t)\left(f(s)-f(0)\right), \text{ and}\\\label{eqn2Bang}
r(s) & > r(\bar{s}) +  \lambda \Delta_{\lambda}(x,t)\left(f(s)-f(\bar{s})\right).
\end{align}
Multiplying \eqref{eqn1Bang} with $\left(1-\frac{f(s)}{f(\bar{s})}\right)$ and  \eqref{eqn2Bang} with $\frac{f(s)}{f(\bar{s})}$ on both sides, and adding them together leads to
\begin{align} 
r(s) & > \left(1-\frac{f(s)}{f(\bar{s})}\right)\cdot \left( r(0) +  \lambda \Delta_{\lambda}(x,t)\left(f(s)-f(0)\right)\right) + \left(\frac{f(s)}{f(\bar{s})} \right)\cdot 
 \left( r(\bar{s}) +  \lambda \Delta_{\lambda}(x,t)\left(f(s)-f(\bar{s})\right)\right) \nonumber \\
 &= \left(1-\frac{f(s)}{f(\bar{s})}\right)\cdot \lambda \Delta_{\lambda}(x,t) f(s) + \frac{f(s)}{f(\bar{s})}r(\bar{s}) +   \left(\frac{f(s)}{f(\bar{s})} \right)\cdot 
 \lambda \Delta_{\lambda}(x,t)\left(f(s)-f(\bar{s})\right) \nonumber \\
 &= \lambda \Delta_{\lambda}(x,t) f(s) -  \left(\frac{f(s)}{f(\bar{s})} \right)\cdot \lambda \Delta_{\lambda}(x,t) f(\bar{s})  + \left(\frac{f(s)}{f(\bar{s})} \right)\cdot r(\bar{s}) \nonumber \\
 &= \frac{f(s)}{f(\bar{s})} \cdot r(\bar{s}). \label{eqn3Bang}
\end{align}
Rewriting Inequality \eqref{eqn3Bang} leads to  $\frac{r(s)}{f(s)} > \frac{r(\bar{s})}{f(\bar{s})}$, which is a contradiction. Hence, $s$ cannot be optimal.
\hfill \Halmos \endproof

\subsection{Proof of Corollary \ref{extremalConvexConcave}}
\proof{Proof.} 
In light of Proposition \ref{bangbangCondition}, we need to show that  the inequality $\frac{r(s)}{f(s)} \leq \frac{r(\bar{s})}{f(\bar{s})}$ holds for all $s\in (0,\bar{s})$ when either condition $(i)$, $(ii)$, or $(iii)$ holds.

It is easy to verify that the above inequality holds when $r(s)=f(s)$ for all $s\in\mathcal{S}$, hence yielding condition $(i)$. For condition $(ii)$, we first note that $f(\cdot)$ and $r(\cdot)$ are assumed to be increasing functions of $s$ (see Assumption \ref{assumptionFunctions}), and thus $a>0$ and $b>0$. We write $
    \frac{r(s)}{f(s)} = \frac{a}{b}s^{\alpha-\beta},$
which is non-decreasing in $s$ for $ab\geq 0$ and $\alpha\geq\beta$. This implies that $\frac{a}{b}s^{\alpha-\beta} \leq  \frac{a}{b}\bar{s}^{\alpha-\beta}$ for all $s\in(0,\bar{s})$ under condition $(ii)$ and Assumption \ref{assumptionFunctions}. 
We now focus on condition $(iii)$.
For any $0<s<\bar{s}$, let $s=t\cdot\bar{s}$, so that $0<t<1$. The convexity of $r(\cdot)$ yields $t\cdot r(\bar{s}) + (1-t)\cdot r(0) \geq  r(t\cdot \bar{s})= r(s),$ which after using $r(0)\leq 0$, leads to 
\begin{align}\label{eqFs}
t\cdot r(\bar{s}) \geq r(s).
\end{align} 
Likewise, using the concavity of $f(\cdot)$ and $f(0)=0$, we obtain
\begin{align}\label{eqRs}
t\cdot f(\bar{s}) \leq f(s).
\end{align}
Combining Inequality \eqref{eqFs} and Inequality \eqref{eqRs} leads to $\frac{r(\bar{s})}{f(\bar{s})} \geq  \frac{r(s)}{f(s)}$ for all $s\in (0,\bar{s})$. 
Note that this result can also be obtained by applying Bauer's maximum principle \citep[see, e.g., Theorem 7.69 of][]{aliprantis1994}. When $r(\cdot)$ is convex and $f(\cdot)$ is concave, our objective function  $r(s) -  \lambda f(s)\Delta_{\lambda}(x,t)$ is convex, and since  $[0,\bar{s}]$ is a compact and convex set, its maximum must be attained at the boundary of the set $[0,\bar{s}]$. 
\hfill \Halmos \endproof

\subsection{Proof of Theorem \ref{optimalT}}
\proof{Proof.} 
We use the shorthand notation $J(T) \triangleq J_{\lambda}^*(0,T)$. Taking the first and second derivative of $g(T)$ with respect to $T$ yields
\begin{align}\label{firstDer}
    g'(T) =\frac{TJ'(T)-J(T)}{T^2},
\end{align}
and, after some algebraic simplifications,
\begin{align}\label{secondDer}
    g''(T) = \frac{T^2J''(T) -2\cdot\big(TJ'(T) - J(T)\big)}{T^3}.
\end{align} 
Now observe that for any $t$ with $g'(t)=0$, we have that $tJ'(t)-J(t)=0$ (see Equation \eqref{firstDer}), and hence that $g''(t) = \frac{t^2 J''(t)-2\cdot0}{t^3} = \frac{J''(t)}{t}$ (see Equation \eqref{secondDer}). We further know that $J(t)$ is concave in $t$, i.e. $J''(t)\leq0$ for all $t\geq 0$  (see Proposition \ref{concaveInT}), so that the following condition holds
\begin{align}\label{contra}
 g''(t)\leq 0, \text{ for all } t\geq 0 \text{ such that }g'(t)=0.
\end{align}

We now proceed with the proof. Let $T^*$ be a strict local maximizer of $g(T)$. Let us assume the contrary of the result, i.e. that $T^*$ is not a strict global maximizer. We will show that this leads to a contradiction. 


Let $T'$ denote the strict global maximizer and consider the case that $T'>T^*$ (the proof follows verbatim for the case that $T'<T^*$). By construction, $g(T)$ attains it maximum over the closed interval $[T^*,T']$ in $T'$. By the extreme value theorem, $g(T)$ also attains a minimum over this closed interval $[T^*,T']$ at least once. Let $t_0$ be the smallest number in $[T^*,T']$ such that $g(t_0)$ is a local minimum.
Observe that since $T^*$ is a strict local maximizer, there exists some $\epsilon_0>0$ such that $g(T^*)>g(T^*+\epsilon)$ for all $0<\epsilon<\epsilon_0$, implying that $t_0\in [T^*+\epsilon,T')$ for some $\epsilon>0$.
As $g(T)$ is differentiable and continuous in $[T^*,T']$, we have $g'(t_0)=0$. 
But by \eqref{contra} this implies that $g''(t_0)\leq 0$, which contradicts that $t_0$ is a local minimizer. Thus $T^*$ must be a strict global maximizer.
\hfill \Halmos \endproof


\subsection{Proof of Proposition \ref{inferenceUnderPolicy}}
\proof{Proof.} 
Let $f_t(\bm{s}_t) \triangleq \int_{u=0}^t f(s_u)du$, so that the mean of the deterioration process at time $t$ under an employed production policy $\bm{s}_t$ and given $\lambda$ is equal to $\lambda f_t(\bm{s}_t)$ (this follows from the fact that the deterioration process is a non-homogenous Poisson process with rate $\lambda f(s_t)$ at time $t$). We can then obtain the posterior distribution at time $t$ using Bayes' theorem. That is, the posterior distribution at time $t$ is proportional to the likelihood function, which is a Poisson distribution with mean $\lambda f_t(\bm{s}_t)$, times the prior distribution: 
\begin{align*}
\varphi_{u}(\lambda|\bm{s},y,\alpha_0,\beta_0) &\propto \P[Y_u= y |\lambda,\bm{s}]\cdot \varphi_{0}(\lambda|\alpha_0,\beta_0) =
 \frac{(\lambda f_t(\bm{s}))^y e^{- \lambda f_t(\bm{s}) }}{y!} \frac{ \beta_0^{\alpha_0}\lambda^{\alpha_0-1} e^{-\beta_0 \lambda}}{\Gamma(\alpha_0)}\\
&=  \frac{\big( f_t(\bm{s})  \big)^y}{y!} \frac{ \beta_0^{\alpha_0} \lambda^{y+\alpha_0-1} e^{-\lambda(\beta_0 + f_t(\bm{s}) ) }}{\Gamma(\alpha_0)} \propto  \lambda^{\alpha_0+y-1} e^{-\lambda(\beta_0+f_t(\bm{s})) }.
\end{align*}
Note that this is again a Gamma distribution with parameters $\alpha_0+y$ and $\beta_0 + f_t(\bm{s}_t)$.
\hfill \Halmos \endproof

\section{Extension to Multiple Production Systems}\label{extensionMulti} In Remark \ref{remarkRev2}, we explained that by appropriately choosing the revenue function $r$, one can model situations where a DM aims to satisfy a (continuous) demand rate $D$ with production rate $s$ (by penalizing the downfall $(D-s)^+$ and rewarding the surplus $(s-D)^+$). One natural extension to this case is when this demand rate $D$ can be satisfied by multiple production systems together. In this section, we discuss this extension.   

\subsection{Problem Formulation}
We have $N$ production systems. Each individual system $i\in\{1,2,\ldots, N\}$ has base rate  $\lambda_i$ and deterioration function $f_i$. For ease of exposition we assume that all systems have the same production rate set $\mathcal{S}$, failure threshold $\xi$, and maintenance cost function $c_m$. (These can also depend on $i$ leading only to additional notation). Let $\bm{s}\triangleq (s_1,s_2,\ldots,s_n)\in \mathscr{S}^N$ denote a vector (bold notation indicates vectors) of production rates, where $s_i$ is the production rate of system $i$. As the DM aims to achieve the common demand rate $D$ with the cumulative production rate from all systems, we let the revenue function $r: \mathscr{S}^N \rightarrow \R_+$ be of the form: 
\begin{align}
    r(\bm{s}) \triangleq  b \cdot \left(\sum_{i=1}^n s_i - D\right)^+ - p\cdot \left(D-\sum_{i=1}^n s_i\right)^+, \label{eq:objectiveMultiple}
\end{align}
where $p$ is a penalty cost rate per unit of shortfall, while $b$ is a bonus rate per unit of surplus (we discuss this in Remark \ref{remarkRev2} for the single system case). In this section, we model the setting that the bonus (penalty) is linear in the surplus (shortfall), but the revenue function easily generalizes to any increasing, continuous function. Since we have a common demand rate $D$, but differing production systems in terms of their base rates and deterioration functions, the DM is interested in determining how to allocate production among these various systems in a way that best meets demand during the planning horizon while also minimizing maintenance costs incurred at the end of the planning horizon.

Let $\bm{x}\triangleq(x_1,x_2,\ldots,x_N)\in \mathscr{X}^N$ be the vector of deterioration levels for all systems, where $x_i$ is the deterioration level of system $i$. Again, we count time backwards, so that at the start of the planning horizon $t=T$, while at the planned maintenance moment we have $t=0$. The feasible set of production decisions at any time $t<T$ given $\bm{x}$ is denoted by $\mathscr{S}(\bm{x})\triangleq \{\bm{s} \in\mathscr{S}^N : s_i = 0 \mbox{ if } x_i=\xi \}$, which means that if system $i$ is failed, i.e. $x_i=\xi$, we can only turn it off. For any time to planned-maintenance $t$, let $J^*_{\bm{\lambda}}(\bm{x},t)$ be the optimal expected profit when the systems have deterioration levels $\bm{x}$ and base rates $\bm{\lambda}\triangleq(\lambda_1,\lambda_2,\ldots,\lambda_N)$. We shall now derive the HJB equation for $J^*_{\bm{\lambda}}(\bm{x},t)$ by following a similar argument as in the proof of Lemma \ref{lemmaHJB}. 

Assume that the state is $(\bm{x},t)\in \mathscr{X}^N \times [T,0)$  and we choose production rates $\bm{s}\in \mathscr{S}(\bm{x})$. The deterioration process is now an $N$-dimensional non-homogeneous Poisson process with rate vector $\left(\lambda_i f_i(s_i)\right)_{i\in \{1,2,\ldots, N\}}$, so we can consider a small interval of length $\delta t$ such that the probability of a deterioration increment at system $i$ is approximately $\lambda_i f_i(s_i) \delta t$. We then have  by the principle of optimality:
\begin{align}\label{heuristicMultiple}
J_{\bm{\lambda}}^*(\bm{x},t) = \max_{\bm{s}\in \mathscr{S}(\bm{x})}\Bigg[ r(\bm{s})\delta t + \sum_{i=1}^N \underbrace{\delta t \lambda_i f_i(s_i) J_{\bm{\lambda}}^*(\bm{x}+\bm{e}_i,t-\delta t)}_{\text{system $i$ deteriorates}}  +  \overbrace{\left(1- \sum_{i=1}^N \delta t \lambda_i f_i(s_i)\right)J_{\bm{\lambda}}^*(\bm{x},t-\delta t)}^{\text{no deterioration}} + o(\delta t)\Bigg], 
\end{align}
where $o(\delta t)$ is a quantity that goes to zero faster than $\delta t$, and where $\bm{e}_i$ denotes the $i$-th standard basis vector in $\R^N$. Subtracting $J_{\bm{\lambda}}^*(\bm{x},t-\delta t)$ from both sides of \eqref{heuristicMultiple}, dividing by $\delta t$ and taking the limit $\delta t \rightarrow 0$ leads to: 
\begin{align}\label{HJBMultiple}
\frac{ \delta J_{\bm{\lambda}}^*(\bm{x},t)}{\delta t} = \max_{\bm{s}\in \mathscr{S}(\bm{x})}\Bigg[ r(\bm{s}) - \sum_{i=1}^N  \lambda_i f_i(s_i) \Big( J_{\bm{\lambda}}^*(\bm{x},t)-J_{\bm{\lambda}}^*(\bm{x}+\bm{e}_i,t)\Big)\Bigg], 
\end{align}
with boundary condition $J_{\bm{\lambda}}^*(\bm{x},0) = -\sum_{i=1}^N c_m(x_i)$ for all $\bm{x}\in\mathscr{X}^N$. 

To see the reasoning behind \eqref{heuristicMultiple}, consider $N=2$, $(\bm{x},t)$ with $x_1,x_2<\xi$ and $t>0$, and $\bm{s}=(s_1,s_2)$. There are four events that can occur: 
\begin{enumerate}
\item Both system 1 and 2 deteriorate and we move to $J^*_{\bm{\lambda}}((x_1+1,x_2+1),t-\delta t)$. This happens with probability $\lambda_1 f_1(s_1) \delta t \cdot \lambda_2 f_2(s_2) \delta t=o(\delta t)$;  
\item Both system 1 and 2 do not deteriorate and we move to  $J^*_{\bm{\lambda}}((x_1,x_2),t-\delta t)$. This happens with probability $(1-\lambda_1 f_1(s_1) \delta t)\cdot(1-\lambda_2 f_2(s_2) \delta t) = 1-\lambda_1 f_1(s_1) \delta t - \lambda_2 f_2(s_2) \delta t - \lambda_1 f_1(s_1) \delta t \cdot \lambda_2 f_2(s_2) \delta t =    1-\big((\lambda_1 f_1(s_1) + \lambda_2 f_2(s_2)\big) \delta t + o(\delta t)$; 
\item System 1 deteriorates but system 2 does not and we move to $J^*_{\bm{\lambda}}((x_1+1,x_2),t-\delta t)$. This happens with probability $\lambda_1 f_1(s_1) \delta t \cdot (1-\lambda_2 f_2(s_2) \delta t) = \lambda_1 f_1(s_1) \delta t - \lambda_1 f_1(s_1) \delta t\cdot \lambda_2 f_2(s_2) \delta t)  = \lambda_1 f_1(s_1) \delta t  + o(\delta t)$;
\item System 2 deteriorates but system 1 does not and we move to $J^*_{\bm{\lambda}}((x_1,x_2+1),t-\delta t)$. This happens with probability $\lambda_2 f_2(s_2) \delta t \cdot (1-\lambda_1 f_1(s_1) \delta t) = \lambda_2 f_2(s_2) \delta t - \lambda_2 f_2(s_2) \delta t\cdot \lambda_1 f_1(s_1) \delta t)  = \lambda_2 f_2(s_2) \delta t  + o(\delta t)$.
\end{enumerate}
Event 1 is $o(\delta t)$, event 2 is captured in the ``no deterioration'' term in in $\eqref{HJBMultiple}$, and event 3 and 4 are captured in the ``system $i$ deteriorates'' term in $\eqref{HJBMultiple}$. For general $N$, the derivation of \eqref{heuristicMultiple} is along similar lines. Note that for $N=1$, Equation \eqref{HJBMultiple} reduces to the single system HJB equation that we formulated in Lemma \ref{lemmaHJB}.

Solving \eqref{HJBMultiple} can be done numerically via (a multivariate version of) the solution approach outlined in Appendix \ref{appendixProcedure}. However, even for moderate values of $N$, $\xi$, and $\bar{s}$ (the maximum production rate),  the size of the state and action space becomes very large; consequently, obtaining the solution
of the HJB recursion \eqref{HJBMultiple} is computationally demanding and can only be done for small instances. Hence, in the following section, we will confine our numerical experiments to settings that involve only two production systems.
\subsection{Numerical Experiments with Two Production Systems}
We first consider two symmetric production systems. Both production systems have the same values as the example in Figure \ref{optimalPolicyIllus}. Thus, both systems have $\xi=10$, $c_m(\xi)=5$, and $c_m(\cdot)=1$ for all non-failed states. Further, $T$ equals 15, and the rate space for both systems equals $\mathscr{S}=[0,1]$.  For the deterioration function, we use $f_1(s)=f_2(s)=s^{2}$, and we set the base rates $\lambda_1$ and $\lambda_2$ to 1. We use $D=1.5$, $p=1$, and $b=0$ for the parameters of the revenue function \eqref{eq:objectiveMultiple}. Thus the cumulative production output of both systems must meet a demand rate of $1.5$, and any shortfall of that required production rate is penalized at rate $1$.
Figure \ref{optimalPolicySym} presents joint sample paths of the optimal production policy and the corresponding controlled deterioration for system 1 (top) and system 2 (bottom). 

\begin{figure}[htb!]
\begin{center}
\includegraphics[width=0.55\textwidth]{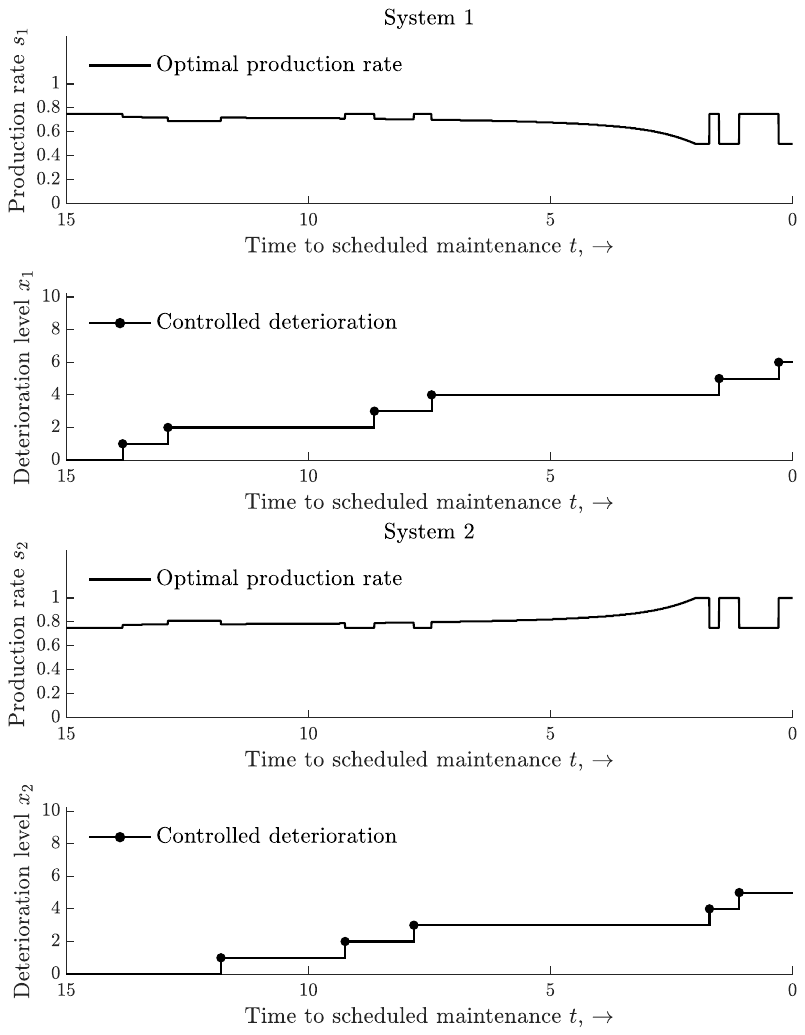}
\caption{Illustration of the optimal condition-based production policy for two symmetric production systems.}
\label{optimalPolicySym}
\end{center}
\end{figure}

Based on Figure \ref{optimalPolicySym}, we can state three important observations. First, when both production systems have equal deterioration levels, they share the demand rate requirement equally (i.e. $s_1=s_2$ when $x_1=x_2$). Second, if one system's condition worsens, its production rate decreases instantly, while the other system instantly compensates by an equal increase, thereby still meeting the demand rate requirement. Third, as the time to the planned maintenance moment decreases, the distribution of the required demand rate across both production system changes when the deterioration levels differ. For example, we clearly see that the optimal production rate of system 1 for fixed deterioration level $x_1=4$ decreases as $t$ decreases, while the optimal production rate of system 2, which has deterioration level $x_2=3$ at that time, increases with the same amount. This latter observation also breaks the monontonic behavior of the optimal production policy in $t$ for a single system (cf. Theorem \ref{optPolicy}). These three observations together indicate that the optimal condition-based production policy for symmetric production systems leads to synchronized controlled deterioration.

Next, we consider two production systems that are asymmetric in their base rates. We use the same values as the symmetric setting, except for the base rates which are set at $\lambda_1=1.5$ and $\lambda_2=0.5$, so that production system 1 deteriorates faster under normal operating conditions than production system 2. Figure \ref{optimalPolicyAsym} presents joint sample paths of the optimal production policy and the corresponding controlled deterioration for system 1 (top) and system 2 (bottom). 
\begin{figure}[htb!]
\begin{center}
\includegraphics[width=0.55\textwidth]{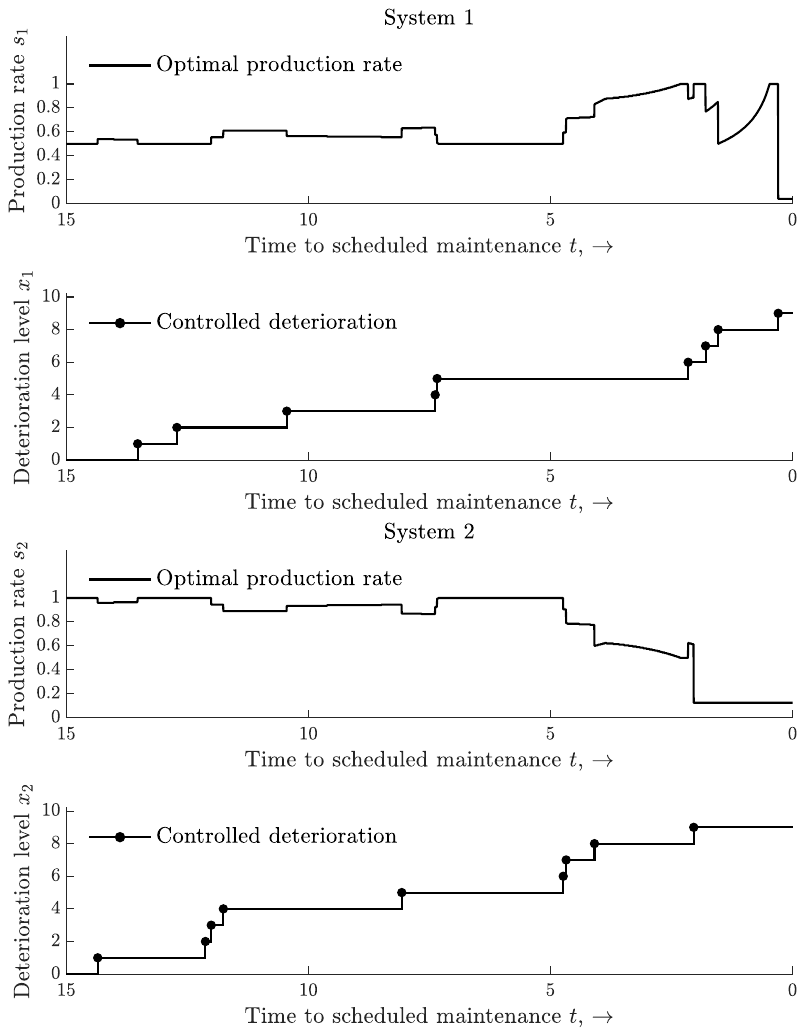}
\caption{Illustration of the optimal condition-based production policy for two asymmetric production systems. System 1 (top) has base rate $\lambda_1=1.5$ and system 2 (bottom) has base rate $\lambda_2=0.5$.}
\label{optimalPolicyAsym}
\end{center}
\end{figure}

Figure \ref{optimalPolicyAsym} suggests that the optimal condition-based production policy for asymmetric production systems also results in synchronized controlled deterioration. However, due to differing base rates, the distribution of the demand rate requirement across the two systems varies and is inversely proportional to these base rates. At the start of the planning horizon, we still see that if one system's condition worsens, its production rate decreases instantly, while the other system instantly compensates by an equal increase, thereby still meeting the demand rate requirement. However, as the time to planned maintenance shortens and both systems approach their failure levels, this compensatory effect becomes less pronounced. Notably, close to the end of the planning horizon, the production systems do not even meet the demand rate requirement any more, indicating that the expected maintenance costs have exceeded those associated with the production shortfall.

\section{Solving Continuous Time Dynamic Programs}\label{appendixProcedure}
 To compute the optimal policy we discretize the planning horizon $[T, 0]$ into $N\triangleq T/ \Delta t$ time intervals of length $\Delta t$ and use a finite difference equation to approximate the optimality equation in \eqref{HJEqn}. Specifically, we replace the partial derivative in Equation \eqref{HJEqn} by the finite difference $\big(J_{\lambda}^*(x,t+\Delta t)-J_{\lambda}^*(x,t)\big) / \Delta t$, after which we obtain the finite difference equation:
\begin{align*}
& J_{\lambda}^*(x,t+\Delta t) = & \\
& \begin{cases}
		J_{\lambda}^*(x,t) + \Delta t \cdot \max_{s\in \mathscr{S}} \Big[r(s) -  \lambda f(s)\big( J_{\lambda}^*(x,t)-J_{\lambda}^*(x+1,t)\big)\Big], & \text{if $x \in \{0,1,\ldots, \xi-1\},$}\\
            J_{\lambda}^*(x,t), & \text{if $x=\xi$},
		 \end{cases}&
\end{align*}
which can be solved by a backward recursion on the discrete time set $\big\{n\Delta t \ |\  n\in\{0,1,\ldots, N\} \big\}$ starting from the boundary condition $J_{\lambda}^*(x,0) = -c_m(x)$ for all $x\in \mathscr{X}$. In all our numerical experiments and illustrations we use a sufficiently small $\Delta t$, such that the approximation is sufficiently accurate.

\section{Detailed Results Numerical Experiments}
\label{sec:fullFactorial}
This section presents detailed results from our numerical experiments. We first discuss the results from the full factorial experiments conducted on systems with known, fixed base rates. Subsequently, we discuss the results of the full factorial experiment involving systems exhibiting unit-to-unit variability in the base rates.

\subsection{Static versus Condition-Based Production}
\label{subsec:staticCBPapp}
To assess the value of condition-based production over static production policies, we perform an extensive numerical study. This numerical study has the following parameter settings. We use base rates $\lambda \in \{0.5, 0.75, 1, 1.25, 1.5 \}$, failure thresholds $\xi\in \{10,12,14,18,20\}$, and planning horizons $T\in \{ 10, 20\}$. We use the functions $f=s^{\gamma}$ and $r=s^{\nu}$, for the degradation and revenue function, respectively, and use the interval $\mathscr{S}=[0,2]$ for the production rates to choose from. By using the values $\{0.5, 0.75, 1, 1.33, 2\}$ for $\gamma$ and $\nu$, we model concave, linear, and convex functions, respectively. We vary the preventive maintenance cost $c_p$ between $\{1,2,3\}$, and choose $c_m(\xi)=10$ as corrective maintenance cost for the failed state.  The complete set of parameter values is summarized in Table \ref{tab:testbed1}. We perform a full factorial test bed, which results in 3750 instances. 

\begin{table*}[!htbp]
  \centering
              \fontsize{8pt}{8pt}\selectfont
  \caption{Input parameter values for numerical study}
  \label{tab:testbed1}
  \begin{tabularx}{0.8\textwidth}{c X c l}
    \toprule
    & Input parameter& No. of choices&Values\\
                             \midrule
1 & Base rate, $\lambda $ & 5 &  0.5, 0.75, 1, 1.25, 1.5      \\
2 & Preventive maintenance cost, $c_p$					 	            & 3 &  1, 2, 3  \\
3 & Failure threshold, $\xi$ & 5 & 10, 12, 14, 18, 20\\
4 & Length of planning horizon, $T$							        & 2 &  10, 20 \\
5 & Parameter of deterioration function, $\gamma$							        & 5 & 0.5, 0.75, 1, 1.33, 2 \\
6 & Parameter of revenue function, $\nu$							        & 5 & 0.5, 0.75, 1, 1.33, 2 \\
   \bottomrule
  \end{tabularx}
\end{table*}

The results of the numerical study are summarized in Table \ref{tab:resultsTestbed1}. 
In this table, we present the average and standard deviation of $\mathcal{R}$, and the $P_{\pi_{CS}}$,  where we first distinguish between subsets of instances with the same value for a specific input parameter of Table \ref{tab:testbed1}, and then present the results for all instances.

\begin{table}[htbp]
\centering
\caption{Relative profit increase (in \%) of condition-based production compared to static production.}
\fontsize{8pt}{8pt}\selectfont
\label{tab:resultsTestbed1}
\begin{tabular}{lccccc}
\toprule
  Input parameter & \multicolumn{1}{l}{Value} &       & Mean $\mathcal{R}$  & SD $\mathcal{R}$  & $P_{\pi_{CS}}$  \\
\midrule
\multirow{5}[1]{*}{Base rate, $\lambda $} & \multicolumn{1}{l}{0.5}                                                            &         & 15.18	&	11.46	&	33.74	\\
& \multicolumn{1}{l}{0.75}                                                                                                    &          & 26.82	&	18.26	&	24.89	\\
& \multicolumn{1}{l}{1}                                                                                                      &           & 41.58	&	34.97	&	19.39	\\
& \multicolumn{1}{l}{1.25}                                                                                                    &          & 72.31	&	133.78	&	15.73	\\
& \multicolumn{1}{l}{1.5}                                                                                                      &         & 96.23	&	425.59	&	13.16	\\
\addlinespace
\multirow{3}[1]{*}{Preventive maintenance cost, $c_p$ } & \multicolumn{1}{l}{1}                                              &       & 55.94	&	262.30	&	22.38	\\
& \multicolumn{1}{l}{2}                                                                                                     &       & 46.11	&	188.32	&	21.38	\\
& \multicolumn{1}{l}{3}                                                                                                     &       & 49.23	&	136.83	&	20.39	\\
\addlinespace
\multirow{5}[1]{*}{Failure threshold, $\xi$} & \multicolumn{1}{l}{10}                                                        &       & 92.22	&	427.09	&	14.65	\\
& \multicolumn{1}{l}{12}                                                                                                    &        & 66.07	&	129.02	&	17.70	\\
& \multicolumn{1}{l}{14}                                                                                                      &      & 42.56	&	44.98	&	20.59	\\
& \multicolumn{1}{l}{18}                                                                                                    &        & 27.65	&	21.93	&	25.81	\\
& \multicolumn{1}{l}{20}                                                                                                      &      & 23.63	&	17.87	&	28.16	\\
\addlinespace
\multirow{2}[1]{*}{Length of planning horizon, $T$ } & \multicolumn{1}{l}{10}                                                &       & 47.94	&	151.18	&	19.47	\\
& \multicolumn{1}{l}{20}                                                                                                     &       & 52.91	&	243.17	&	23.29	\\
\addlinespace
\multirow{5}[1]{*}{Parameter of deterioration function, $\gamma$ } & \multicolumn{1}{l}{$0.5$}                        &       & 67.08	&	376.91	&	31.59	\\
& \multicolumn{1}{l}{$0.75$}                                                                                         &       & 61.19	&	173.84	&	25.54	\\
& \multicolumn{1}{l}{1}                                                                                                     &       & 48.59	&	168.19	&	20.75	\\
& \multicolumn{1}{l}{$1.33$}                                                                                         &       & 41.02	&	45.08	&	16.42	\\
& \multicolumn{1}{l}{2}                                                                                                     &       & 34.25	&	42.68	&	12.60	\\
\addlinespace
\multirow{5}[1]{*}{Parameter of revenue function, $\nu$}& \multicolumn{1}{l}{$0.5$}                                 &       & 28.95	&	48.01	&	12.14	\\
& \multicolumn{1}{l}{$0.75$}                                                                                         &       & 32.31	&	43.44	&	13.35	\\
& \multicolumn{1}{l}{1}                                                                                                     &       & 53.79	&	169.41	&	15.83	\\
& \multicolumn{1}{l}{$1.33$}                                                                                         &       & 67.49	&	175.60	&	21.49	\\
& \multicolumn{1}{l}{2}                                                                                                     &       & 69.58	&	374.14	&	44.09	\\
\midrule
Total &                                                                                                                     &       & 50.42 & 202.46 & 21.38 \\
    \bottomrule
    \end{tabular}%
\end{table}%

\subsection{Sequential versus Integrated Maintenance Scheduling}
\label{subsec:seqIntApp}
To assess the value of integrating maintenance and production policies, we consider a full factorial test bed using the same parameters of the previous full factorial experiment as displayed in Table \ref{tab:testbed1} (excluding the parameter $T$ for obvious reasons), which results in 1875 instances. We summarize the results of this numerical study in Table \ref{tab:resultsTestbed2}. Again, we present the average and standard deviation of $\mathcal{\hat{R}}$, and the $\hat{P}_{\pi_{I}}$,  where we first distinguish between subsets of instances with the same value for a specific input parameter of Table \ref{tab:testbed1}, and then present the results across all instances.

\begin{table}[htbp]
\centering
\caption{Relative profit rate increase (in \%) of integrating maintenance and production decisions compared to treating them sequentially.}
\fontsize{8pt}{8pt}\selectfont
\label{tab:resultsTestbed2}
\begin{tabular}{lccccc}
\toprule
Input parameter & \multicolumn{1}{l}{Value} &       & Mean $\mathcal{\hat{R}}$   & SD $\mathcal{\hat{R}}$  & $\hat{P}_{\pi_{I}}$  \\
\midrule
\multirow{5}[1]{*}{Base rate, $\lambda $} & \multicolumn{1}{l}{0.5}                                                            &       &	26.18	&	30.00 &	2.04	\\
& \multicolumn{1}{l}{0.75}                                                                                                    &       & 23.22	&	28.99	&	1.92	\\
& \multicolumn{1}{l}{1}                                                                                                      &       &	20.80	&	28.09	&	1.82	\\
& \multicolumn{1}{l}{1.25}                                                                                                    &      & 19.04	&	27.29	&	1.73	\\
& \multicolumn{1}{l}{1.5}                                                                                                      &       & 17.70	&	26.59	&	1.64	\\
\addlinespace
\multirow{3}[1]{*}{Preventive maintenance cost, $c_p$ } & \multicolumn{1}{l}{1}                                              &       &	18.74	&	24.84	&	2.04	\\
& \multicolumn{1}{l}{2}                                                                                                     &       & 21.44	&	28.38	&	1.82	\\
& \multicolumn{1}{l}{3}                                                                                                     &       & 23.99	&	31.25	&	1.63	\\
\addlinespace
\multirow{5}[1]{*}{Failure threshold, $\xi$} & \multicolumn{1}{l}{10}                                                       &       &	15.95	&	23.72	&	1.63	\\
& \multicolumn{1}{l}{12}                                                                                                    &       & 18.39	&	25.87	&	1.75	\\
& \multicolumn{1}{l}{14}                                                                                                    &       & 20.80	&	27.75	&	1.83	\\
& \multicolumn{1}{l}{18}                                                                                                    &       &	25.05	&	30.62	&	1.95	\\
& \multicolumn{1}{l}{20}                                                                                                    &       &	26.76	&	31.70	&	1.99	\\
\addlinespace
\multirow{5}[1]{*}{Parameter of deterioration function, $\gamma$ } & \multicolumn{1}{l}{$0.5$}                        &       &	2.73	&	1.98	&	1.98	\\
& \multicolumn{1}{l}{$0.75$}                                                                                         &       &	8.46	&	4.56	&	1.92	\\
& \multicolumn{1}{l}{1}                                                                                                     &       &   17.40	&	8.99	&	1.87	\\
& \multicolumn{1}{l}{$1.33$}                                                                                         &       &	30.52	&	18.88	&	1.78	\\
& \multicolumn{1}{l}{2}                                                                                                     &       &	47.84	&	47.40	&	1.60	\\
\addlinespace
\multirow{5}[1]{*}{Parameter of revenue function, $\nu$}& \multicolumn{1}{l}{$0.5$}                                 &       &	5.70	&	6.25	&	1.00	\\
& \multicolumn{1}{l}{$0.75$}                                                                                         &       &   11.25	&	8.79	&	1.22	\\
& \multicolumn{1}{l}{1}                                                                                                     &       &   18.23	&	14.24	&	1.50	\\
& \multicolumn{1}{l}{$1.33$}                                                                                         &      &	27.73	&	23.56	&	1.99	\\
& \multicolumn{1}{l}{2}                                                                                                     &       &	44.04	&	47.30	&	3.43	\\
\midrule
Total &                                                                                                                     &       & 21.39	&	28.35	&	1.83	\\
    \bottomrule
    \end{tabular}%
\end{table}%

\subsection{Performance of the Certainty-Equivalent Policy}
\label{subsec:cePolicyApp}
To assess the performance of the certainty-equivalent policy, we use a full-factorial design in which we vary five different input parameters (see Table \ref{tab:testbedSim} for a summary). We use failure thresholds $\xi\in \{10, 15, 20\}$ and planning horizons $T\in \{ 10, 15, 20\}$. Again, we use a single corrective maintenance cost equal to 10 and a single preventive maintenance cost for all non-failed states, denoted with $c_p$, which we vary between two values $1$ and $5$. We use the functions $f=s^{\gamma}$ and $r=s^{\nu}$, for the degradation and revenue function, respectively, and use the interval $\mathscr{S}=[0,2]$ to choose  the production rates from. By using the values $\{0.5, 1, 2\}$ for $\gamma$ and $\nu$, we model concave, linear, and convex functions, respectively. The certainty-equivalent policy is characterized by the parameter $N^{\text{opt}}$, which we vary between 5 values. In a simulation instance (we shortly describe the simulation procedure), the base rates are sampled from a Gamma distribution in which we fix the mean to 1, and vary its coefficient of variation, denoted with $cv_{\Lambda}$.  This coefficient of variation is thus a measure for the population heterogeneity of the systems. These choices lead to a total of 4050 instances.

\begin{table*}[!htbp]
  \centering
              \fontsize{8pt}{8pt}\selectfont
  \caption{Input parameter values for simulation study}
  \label{tab:testbedSim}
  \begin{tabularx}{0.8\textwidth}{c X c l}
    \toprule
    & Input parameter& No. of choices&Values\\
                             \midrule
1 & Coefficient of variation of prior Gamma distribution, $cv_\Lambda $ & 3 &  $\sqrt{\frac{1}{2}}$, 1, $\sqrt{2}$      \\
2 & Preventive maintenance cost, $c_p$					 	            & 2 &  1, 3  \\
3 & Length of planning horizon, $T$							        & 3 &  10, 15, 20 \\
4 & Failure threshold, $\xi$					 	             & 5 &  10, 12, 15, 18, 20  \\
5 & Parameter of deterioration function, $\gamma$							        & 3 &  0.5, 1, 2 \\
6 & Parameter of revenue function, $\nu$							        & 3 &  0.5, 1, 2 \\
7 & Number of re-optimizations, $N^{\text{opt}}$ & 5 &  1, 2, 4, 6, 8 \\
   \bottomrule
  \end{tabularx}
\end{table*}

The results of the simulation study are summarized in Table \ref{tab:resultsTestbed}. 
In this table, we present the average and maximum $\mathcal{\bar{R}}$, and the $\bar{P}_{\pi_{O}}$,  where we first distinguish between subsets of instances with the same value for a specific input parameter of Table \ref{tab:testbedSim}, and then present the results across all instances.

\begin{table}[htbp]
\centering
\caption{Results of simulation study: Relative regret (profit loss in \%) compared with  the Oracle performance.}
\fontsize{8pt}{8pt}\selectfont
\label{tab:resultsTestbed}
\begin{tabular}{lccccc}
\toprule
  Input parameter & \multicolumn{1}{l}{Value} &       & Mean $\mathcal{\bar{R}}$   & Max $\mathcal{\bar{R}}$ & $\bar{P}_{\pi_{O}}$  \\
\midrule
\multirow{3}[1]{*}{Coefficient of variation of prior Gamma distribution, $cv_\Lambda $} & \multicolumn{1}{l}{ $\sqrt{\frac{1}{2}}$}      &       & 0.59	&	7.29	&	10.32	\\
& \multicolumn{1}{l}{1}                                                                                                                  &       & 0.72	&	9.75	&	10.33 \\
& \multicolumn{1}{l}{ $\sqrt{2}$}                                                                                                        &       & 0.84	&	11.63	&	10.33 \\
\addlinespace
\multirow{2}[1]{*}{Preventive maintenance cost, $c_p$} & \multicolumn{1}{l}{1}                                                            &       & 0.65	&	10.06	&	11.33	\\
& \multicolumn{1}{l}{5}                                                                                                         &       &  0.78	&	11.63	&	9.33	\\

\addlinespace
\multirow{3}[1]{*}{Length of planning horizon, $T$		} & \multicolumn{1}{l}{10}                                               &       & 0.38	&	3.10	&	7.03	\\
& \multicolumn{1}{l}{15}                                                                                                           &       & 0.70	&	7.16	&	10.51	\\
& \multicolumn{1}{l}{20}                                                                                                         &       & 1.08	&	11.63	&	13.46	\\
\addlinespace
\multirow{5}[1]{*}{Failure threshold, $\xi$	} & \multicolumn{1}{l}{10}                                                          &       & 1.37	&	11.63	&	9.03	\\
& \multicolumn{1}{l}{12}                                                                                                       &       & 0.94	&	7.98	&	9.75	\\
& \multicolumn{1}{l}{15}                                                                                                         &       & 0.64	&	4.98	&	10.32	\\
& \multicolumn{1}{l}{18}                                                                                                       &       & 0.34	&	2.33	&	11.14	\\
& \multicolumn{1}{l}{20}                                                                                                         &       & 0.28	&	2.09	&	11.43	\\
\addlinespace
\multirow{3}[1]{*}{Parameter of deterioration function, $\gamma$	} & \multicolumn{1}{l}{0.5}                 &       & 0.01	&	0.07	&	9.93	\\
& \multicolumn{1}{l}{1}                                                                                                       &       & 0.78	&	11.63	&	10.24	\\
& \multicolumn{1}{l}{2}                                                                                                         &       & 1.39	&	11.63	&	10.83	\\
\addlinespace
\multirow{3}[1]{*}{Parameter of revenue function, $\nu$	} & \multicolumn{1}{l}{0.5}                 &               & 1.39	&	11.63	&	10.83	\\
& \multicolumn{1}{l}{1}                                                                                                       &       & 0.77	&	11.63	&	10.24	\\
& \multicolumn{1}{l}{2}                                                                                                         &       & 0.01	&	0.05	&	9.93	\\
\addlinespace
\multirow{5}[1]{*}{Number of re-optimizations, $N^{\text{opt}}$} & \multicolumn{1}{l}{1}                 &       & 0.91	&	11.63	&	10.32	\\
& \multicolumn{1}{l}{2}                                                                                                       &       & 0.79	&	10.21	&	10.33	\\
& \multicolumn{1}{l}{4}                                                                                                         &       & 0.68	&	8.64	&	10.33	\\
& \multicolumn{1}{l}{6}                                                                                                       &       & 0.62	&	7.90	&	10.33	\\
& \multicolumn{1}{l}{8}                                                                                                         &       & 0.59	&	7.55	&	10.33	\\
\midrule
Total &                                                                                                                         &       & 0.72	&	11.63	&	10.33	\\
    \bottomrule
    \end{tabular}%
\end{table}%

Table \ref{tab:resultsTestbed} indicates that the certainty-equivalent policy performs excellently with average profit losses of only 0.72\% on average and a maximum profit loss of only 11.63\% over all 4050 instances. The instance with $cv_{\Lambda}=\sqrt{2}$, $c_p=5$, $T=20$, $\xi = 10$, $\gamma = 1$ or $2$, $\nu=0.5$ or $\nu=1$ and $N^{\text{opt}}=1$ leads to this maximum profit loss of 11.63\%. As such, the value 11.63 occurs as the max $\mathcal{\bar{R}}$ for all subsets of instances containing this poorly performing instance in Table \ref{tab:resultsTestbed}. Indeed, we observe that if you deviate from these values and hence exclude the instance with the highest max $\mathcal{\bar{R}}$, the max $\mathcal{\bar{R}}$ of the resulting subset of instances is smaller. This is the case, when $cv_{\lambda}<\sqrt{2}$, $c_p = 5$, $T<20$, $\xi>10$, $\gamma = 0.5$, $\nu = 2$, or $ N^{\text{opt}}>1$. 

Some insightful observations can be drawn from Table \ref{tab:resultsTestbed}. We find that the profit loss is inversely related to the number of re-optimizations and the failure threshold. The decrease in the number of re-optimizations is intuitive; the certainty equivalent policy adapts more frequently to the accrued knowledge of the current system's base rate yielding a lower profit loss with respect to the optimal policy for that base rate. When the failure threshold increases there is less need for controlled deterioration as the system can, on average, be set to its maximum production rate more often. This holds for both the Oracle policy and the certainty-equivalent policy resulting in comparable performance.  Finally, we see that the form of both the degradation and revenue function has a huge impact on the relative profit loss; a concave deterioration function and/or a convex revenue function lead to negligible profit losses. In light of Corollary \ref{extremalConvexConcave}, these results suggest that the certainty-equivalent policy performs comparable to the Oracle policy in the bang-bang regime, regardless of the degree of unit-to-unit variability. In \S\ref{subsec:robustCE}, we numerically confirm this intuition. 
\end{APPENDICES}
\end{document}